\documentclass[10pt,a4paper]{amsart}
\usepackage{graphicx}
\usepackage{amsmath}
\usepackage{amsfonts}
\usepackage{amssymb}
\usepackage{amsthm}
\usepackage{color}
\begin{document}

\def\R{\mathbb{R}}
\def\N{\mathbb{N}}
\def\H{\mathcal{H}}
\def\d{\textrm{div}}
\def\v{\textbf{v}}
\def\I{\hat{I}}
\def\B{\hat{B}}
\def\x{\hat{x}}
\def\y{\hat{y}}
\def\p{\hat{\phi}}
\def\P{\mathcal{P}}
\def\r{\hat{r}}
\def\w{\textbf{w}}
\def\u{\textbf{u}}
\def\K{\mathcal{K}}
\def\Reg{\textrm{Reg}}
\def\s{\textrm{Sing}}
\def\sgn{\textrm{sgn}}

\numberwithin{equation}{section} 

\newtheorem{lem}{Lemma}[section]
\newtheorem{rem}{Remark}[section]
\newtheorem{cor}{Corollary}[section]
\newtheorem{thm}{Theorem}[section]
\newtheorem{prop}{Proposition}[section]
\newtheorem{definition}{Definition}[section]
\newtheorem{con}{Conjecture}[section]
\newtheorem{Main}{Main Result}

\title{Regularity up to the Crack-Tip for the Mumford-Shah problem}

\author{John Andersson and Hayk Mikayelyan}

\thanks{Department of Mathematics, Royal Institute of Technology,
Lindstedtsv\"agen 25,  100 44 Stockholm, Sweden  \hfill johnan@kth.se} 

\thanks{Mathematical Sciences, University of Nottingham Ningbo, 199 Taikang East Road, Ningbo 315100, PR China   \hfill Hayk.Mikayelyan@nottingham.edu.cn}

\begin{abstract}
 We will prove that if $(u,\Gamma)$ is a minimizer of the functional
 $$
 J(u,\Gamma)=\int_{B_1(0)\setminus \Gamma}|\nabla u|^2dx +\H^1(\Gamma)
 $$
 and $\Gamma$ connects $\partial B_1(0)$ to a point in the interior, then $\Gamma$ satisfies a point-wise $C^{2,\alpha}$-estimate at the crack-tip.

 This means that the Mumford-Shah functional satisfies an additional,
 and previously unknown, Euler-Lagrange condition.
\end{abstract}

\maketitle

\tableofcontents

\newpage
\section{Introduction.}

The Mumford-Shah functional 
\begin{equation*}
J(u,\Gamma)=\int_{\Omega\setminus \Gamma}|\nabla u|^2+\alpha (u(x)-h(x))^2dx +\beta\H^{n-1}(\Gamma\cap \Omega)
\end{equation*}
was introduced by David Mumford and Jayant Shah in \cite{MS} in
the context of image processing problems. The idea is to find $u$ as 
the ``piecewise smooth'' approximation of the 
given raw image data $h(x)$.

In this paper we will investigate the regularity of the free discontinuity set $\Gamma$ of minimizers to
the following simplified version of the Mumford-Shah functional 
\begin{equation}\label{MumfordShah}
J(u,\Gamma)=\int_{\Omega\setminus \Gamma}|\nabla u|^2dx +\H^1(\Gamma\cap \Omega)
\end{equation}
where $\Omega \subset \R^2$ is a given set, $u\in W^{1,2}(\Omega\setminus \Gamma)$, $u\big\lfloor_{\partial \Omega}=g$ and $\Gamma$ is a one dimensional
set that is not apriori determined. In particular,
minimizing (\ref{MumfordShah}) involves finding a pair $(u,\Gamma)$ where
the function $u$ is allowed to be discontinuous across $\Gamma$ but the
set $\Gamma$ can not have too large one-dimensional
Hausdorff measure $\H^1(\Gamma)$. We will call the set $\Gamma$ the free discontinuity set. Notice that the 
absence 
of the term $\alpha(u(x)-h(x))$ with the image data requires imposing boundary conditions. 

It should also be remarked 
that the a model, suggested by the British engineer Alan~Arnold Griffith (see \cite{G}), of brittle fracture is based on the
balance between gain in surface energy and strain energy release involves minimization of an integral very similar to (\ref{MumfordShah}). We believe that our results are 
more relevant for this interpretation of the minimizers rather than to image processing.

The existence and regularity of the Mumford-Shah minimizers started by pioneering works of
Ennio De Georgi, Michele Carriero, Antonio Leaci,
Luigi Ambrosio, Guy David, Alexis Bonnet, Nicola Fusco and Diego Pallara
(see \cite{DCL}, \cite{AmbExist}, \cite{DCone}, \cite{B}, \cite{AP}, \cite{AFP}, \cite{BD}).
Most of the known results can be found in the following two monographs \cite{AFPBook}, \cite{DBook}.
We would also like to mention some recent publications, such as \cite{KLM}, \cite{DFi}, \cite{CL}, \cite{DFo} and \cite{LemMik}.

The regularity analysis near the crack-tip is of particular difficulty and interest,
since the crack-tip is the only singularity where the bulk energy
and the surface energy in a ball $B_r$ scale of the same order as $r\to 0$. As a result one cannot exploit the
domination of the surface term over the bulk term which essentially means that the regularity of $\Gamma$ is 
determined by the minimal surface equation.

In the paper we apply linearization techniques to derive the full asymptotic of both, analytic and geometric, 
components of the minimizers near the crack-tip. These classical
techniques has been recently successfully exploited by the first author, Henrik Shahgholian and Georg Weiss in 
regularity analysis of several free boundary problems (see \cite{ASW}, \cite{ASign}). In 
this article we adapt those free boundary theory methods to free discontinuity context. The exact asymptotic of the
minimizers allows to carry out variations of the discontinuity set in the orthogonal direction
near the crack-tip and derive a previously unknown Euler-Lagrange condition for the Mumford-Shah functional
in Section \ref{OrthVarSec}.

We believe that our results are of importance for the regularity analysis of the quasi-static crack-propagation
model of Gilles A. Francfort and Jean-Jacques Marigo (see \cite{FM}, \cite{DT}). 
In particular, \cite{DT} showed existence of solutions for Griffith's model of quasi-static propagation
in brittle materials. The models involve minimizing $J(u(x,t),\Gamma_{u(x,t)})$ for each time $t\in [0,T]$ where $u(x,t)=g(x,t)$ on $\partial \Omega$
and under the extra constraint $\Gamma_{u(x,s)}\subset \Gamma_{u(x,t)}$ for $s\le t$. Their method of analysis is by time discretization.
That is, they minimize the problem at discrete times $t_k=k\delta$ and find solutions to the Mumford-Shah problem $(u(x,t_k), \Gamma_{u(x,t_k)})$
with the extra condition $\Gamma_{u(x,t_k)}\subset \Gamma_{u(x,t_{k+1})}$. By sending $\delta \to 0$ they recover a solution to the original
problem. In order to analyze this time discretized problem we need to analyze the free discontinuity set
$\Gamma_{u(x,t_{k})}\setminus \Gamma_{u(x,t_{k-1})}$ for each $k$. Heuristically, if $\delta$ is small then the set 
$\Gamma_{u(x,t_{k})}\setminus \Gamma_{u(x,t_{k-1})}$ will be close to a crack-tip in a small ball. It will, in particular, be of great importance to be able to make variations {\sl ``in the orthogonal direction''}
as we do in Section \ref{OrthVarSec} in order to analyze the growth of a fracture. We plan to address this problem in a future article.


In the next subsection we will describe, in more detail, the problem setting. Then we will briefly summarize
the background and the relevant known results for the Mumford-Shah problem. We end the introduction with stating our main results.

\subsection{Problem Setting.}

The Mumford-Shah problem consists, in the context of this article, in finding the pair $(u,\Gamma_u)$ that minimizes the energy (\ref{MumfordShah}) with some prescribed boundary
values $u=g$ on $\partial \Omega$. In general, the set $\Gamma$ can be very complicated (disconnected etc.) as in the left figure below.

\begin{center}
\includegraphics[height=5cm,width=10cm]{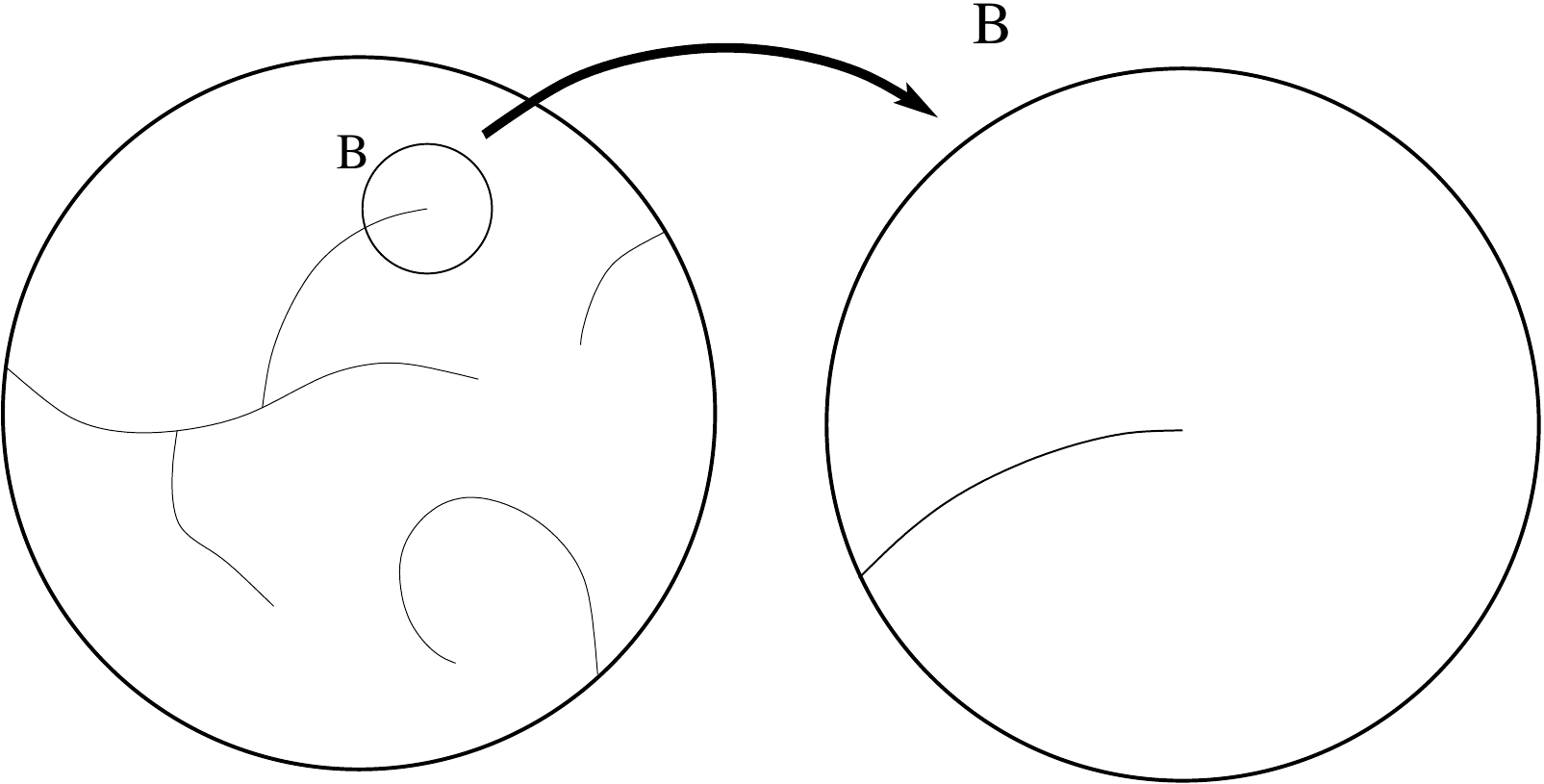}
\end{center}

\vspace{2mm}

{\bf Figure 1:} The geometry of a blow-up.

\vspace{2mm}

It has been shown (see \cite{AP}, \cite{AFP}, \cite{AFPBook} and \cite{DBook}) that the free discontinuity set $\Gamma_u$ is a $C^{1,\alpha}-$graph in
some coordinate system around $\H^1-$almost every point $x^0\in \Gamma_u$. It is also known that certain singularities exist \cite{DBook}
such as {\sl spider points}, consisting of three arcs meeting at $120^\circ$ angles in a point, or crack-tips, points
where the free discontinuity set $\Gamma_u$ ends at some point. The right figure above shows a typical crack-tip.

Not much is known about the regularity properties of the free discontinuity set close to the crack-tip. In particular, it is not known
if the free discontinuity set can spiral around a crack-tip point infinitely many times or if the blow-up of the solution is unique at the crack-tip.
Our main goal is to analyze the behavior of $\Gamma_u$ close to a crack-tip. In this article we will provide an analysis of the crack-tip and
exclude spiraling behaviors as well as providing good regularity estimates at the crack-tip.

In order to specify what we mean by a crack-tip we make the following definition.

\begin{definition}\label{epsCloseDef}
 We say that $(u,\Gamma)$ is $\epsilon-$close to a crack-tip if the following holds:
 \begin{enumerate}
  \item\label{SomeShit} $(u,\Gamma)$ is a minimizer of (\ref{MumfordShah}) in $B_1(0)\,(=\Omega)$ with some specified boundary conditions.
  \item\label{connectedCurve} $\Gamma$ consists of a connected rectifiable curve that starts at the origin and connects the origin 
  to $\partial B_1(0)$. I.e. there exists a Lipschitz mapping $\tau:[0,1]\mapsto \overline{B_1(0)}$ such that $\overline{\Gamma}=\{\tau(s);\; s\in [0,1]\}$, $\tau(0)=0$ and $\tau(1)\in \partial B_1(0)$.
  \item\label{BranchCut} For some $\lambda\in \R$
  $$
  \left(\int_{B_1\setminus\Gamma} \left|\nabla \left( u-\lambda \sqrt{\frac{2}{\pi}}r^{1/2}\sin(\phi/2)\right)\right|^2 dx\right)^{\frac{1}{2}}\le\epsilon.
  $$
  where $(r,\phi)$ are the standard polar coordinates of $\R^2$.
  \item and $u(0,0)=0$.
 \end{enumerate}
\end{definition}

\vspace{1mm}

{\bf Remark:} {\sl As we point out in the remark after Theorem \ref{Ambrosio} there is no loss of generality to assume that 
$\Gamma_u$ is a $C^{1,1/4}$ graph in $B_1\setminus B_{1/2}$. We will, for simplicity, make this assumption throughout the paper.}

\vspace{1mm}

Let us briefly indicate this definition is not vacuous and that solutions that are $\epsilon-$close to a crack-tip indeed do exist.

In \cite{AmbExist} L. Ambrosio showed existence of minimizers in the space of special functions of
bounded variation SBV, in that case the free discontinuity set $\Gamma_u$ is considered to be the singular support of the measure $\nabla u$.
Dal Maso, Morel and Solimini \cite{DMSActa} showed existence of minimizers under the restriction that $\Gamma_u$ should consist
of at most $k$ components, each consisting of the image of a Lipschitz curve. In \cite{DMSActa} it was also showed that the minimizer with $\Gamma_u^k$ consisting of at most $k$
components will converge to the SBV minimizer as $k\to \infty$, see \cite{DMSActa} for exact statements and details. 
For definiteness we can think of minimizers in the sense of \cite{DMSActa}.

To see that the assumption (\ref{connectedCurve}) in Definition \ref{epsCloseDef} is not very restrictive we consider a minimizer constructed in
\cite{DMSActa} with $\Gamma_u$ consisting of a finite number of components. Then if $x^0$ is at a crack-tip, that is if $x_0$ is contained in exactly one 
of the pieces of $\Gamma_u$ and that $x_0$ is the end point of that piece.
To be specific we may assume that $\Gamma$ is parametrized by the Lipschitz 
curves $\gamma_i: I_i\mapsto B_1(0)$, for $i=1,2,...,n$, where $I_i=[a_i,b_i]$ and that $x^0=\gamma_1(a_1)$ and that 
\begin{equation}\label{Archimedes}
x^0\notin \cup_{i=2}^n \{\gamma_i(x); \; x\in I_i\}\bigcup \{\gamma_1(x);\; (a_1,b_1]\}.
\end{equation}
Then, since each of the sets $\{\gamma_i(x); \; x\in I_i\}$ are closed, there is 
 a small ball
$B_r(x^0)$ such that $\Gamma_u$ connects $\partial B_r(x^0)$ to $x^0$ and 
$$
\Gamma_u\cap \{\gamma_i;\; x\in I_i\}=\emptyset \textrm{ for }i=2,3,...,n.
$$
That is $\gamma_1$ parametrizes $\Gamma_u$ in $B_r(x^0)$. We may then define
\begin{equation}\label{Rescaled}
u_r(x)=\frac{u(rx+x^0)}{\sqrt{r}}.
\end{equation}
By scaling invariance of the functional $u_r$ is a minimizer in $B_1(0)$; thus $u_r$ satisfies \ref{SomeShit} of Definition \ref{epsCloseDef}.
Also, since $\Gamma_u$ have finitely many components, $\Gamma_u$ will satisfy \ref{connectedCurve} of Definition \ref{epsCloseDef} if $r$
is small enough. It should be remarked here that if $r$ is small enough we may,
in view of Theorem \ref{Ambrosio}, which we may apply (possibly after decreasing $r$) by Corollary \ref{BlowUpAtCarackTipCor} (which clearly applies since 
$\Gamma_u$ is parametrized by only one curve in small balls centered at $x^0$), 
assume that $\Gamma_u \cap (B_{2r}(x^0)\setminus B_{r/2}(x^0))$ is the graph of a single $C^{1,1/4}$ function which implies that the free discontinuity set does not 
cross $\partial B_r(x^0)$ more than once. The argument is standard and therefore left to the reader.

The geometry
of the situation is indicated in in Figure 1 where we have tried to depict that in a small ball around a crack-tip the free discontinuity set
is a curve connecting the origin to the boundary. The right picture in Figure 1 shows the free discontinuity set of the rescaled function $u_r$.

If we let $r$ be small enough in (\ref{Rescaled}) then \ref{BranchCut} of Definition \ref{epsCloseDef} will be satisfied by $u_r$; this was shown in \cite{B}
see also Corollary \ref{BlowUpAtCarackTipCor}.
In \ref{BranchCut} of Definition \ref{epsCloseDef}, as well as later in this article, we place the branch-cut of $r^{1/2}\sin(\phi/2)$
along $\Gamma_u$. We need to have a branch-cut since $r^{1/2}\sin(\phi/2)$ is two-valued in the plane and therefore we need to choose a branch 
of $r^{1/2}\sin(\phi/2)$, we make this choice of branch so that $r^{1/2}\sin(\phi/2)$ is continuous in $B_1(0)\setminus \Gamma_u$.

Since it is well known that minimizers of the Mumford-Shah functional are $C(B_1\setminus \Gamma_u)$ the final assumption in
Definition \ref{epsCloseDef} can be achieved by adding a constant to $u$.

It is therefore clear that we can construct minimizers $(u,\Gamma_u)$, as in \cite{DMSActa}, such that up to a rescaling, as in (\ref{Rescaled}),
and an additive constant the solution is $\epsilon-$close to a crack-tip at every end point of the set $\Gamma_u$. All our Theorems are valid for 
SBV-minimizers, as considered in \cite{AmbExist}, given that the conditions of Definition \ref{epsCloseDef} are satisfied. 

\vspace{2mm}

{\bf Variations:} Before we continue we need to talk about variations and Euler-Lagrange equations. We have talked about two kinds of minimizers,
and for definiteness mentioned that we are interested in the minimizers where $\Gamma_u$ has finitely many components.
In practice we are interested in any kind of minimizer that is $\epsilon-$close to a crack-tip that satisfies the three types of variations that we
will discuss presently.

The first type of variations are variations in $u$. Let $(u,\Gamma_u)$ be a minimizer in $\Omega$ and let $\Sigma \subset \Omega$.
Then for any function $v\in W^{1,2}(\Sigma)$ such that $v=0$ on $\partial \Sigma \setminus \Gamma_u$ we have that $(u+\epsilon v,\Gamma_u)$ is a competitor for
minimality in the subset $\Sigma$. Extending $v$ by zero in $\Omega$ we can calculate
$$
0=\frac{d J(u+\epsilon v,\Gamma_u)}{d\epsilon}\bigg\lfloor_{\epsilon=0}dx=2\int_{\Sigma\setminus \Gamma_u} \nabla v\cdot \nabla u dx=
2\int_{\partial (\Sigma\setminus \Gamma_u)}v\frac{\partial u}{\partial \nu}d\H^1.
$$
Thus minimizers $u$ satisfy a Neumann boundary condition on $\Gamma_u$.

The second type of variations are the domain variations. These variations we will only do away form the crack-tip. For any function
$\eta\in C^{\infty}(\Omega; \R^2)$ with compact support it follows that
$v_\epsilon(x)=u(x+\epsilon \eta(x))$ and $\Gamma_{v_\epsilon}=\{x;\; x+\epsilon \eta(x)\in \Gamma_u\}$
is a competitor for minimality. A standard calculation, \cite{AFPBook}, leads to
\begin{equation}\label{DomVar}
0=\lim_{\epsilon\to 0}\frac{J(u,\Gamma_u)-J(v_\epsilon,\Gamma_{v_\epsilon})}{\epsilon}=
\end{equation}
$$
=\int_\Omega \left[|\nabla u|^2\textrm{div}(\eta)-2\langle\nabla u, \nabla\eta\cdot\nabla u\rangle \right]dx+\int_{\Gamma_u}\textrm{div}_{\Gamma_u}\eta d\H^1,
$$
where $\textrm{div}_{\Gamma_u} (\eta)$ is the tangential divergence on 
$\Gamma_u$. For us, the most important situation will be when
$\Gamma_u$ is a graph of a $C^1$ function on the support of $\eta$. 

Choosing $\eta=\xi(x_1,x_2)\frac{(-f'(x_1),1)}{\sqrt{1+|f'|^2}}$ and converting the volume integral in (\ref{DomVar}) to a boundary integral by means of an
integration by parts, assuming that everything is smooth, shows that
\begin{equation}\label{PointwiseCurvW}
\int_{\Gamma_u}\left[|\nabla u(x_1,f(x_1)^+)|^2-|\nabla u(x_1,f(x_1)^-)|^2\right]\xi(x_1,f(x_1))d\H^1=
\end{equation}
$$
-\int_{\Gamma_u}\frac{\partial}{\partial x_1}
\left( \frac{f'(x_1)}{\sqrt{1+(f'(x_1))^2}}\right)\xi(x_1,f(x_1))d\H^1.
$$

We interpret equation (\ref{PointwiseCurvW}) as
\begin{equation}\label{PointwiseCurv}
 \frac{\partial}{\partial x_1}\left( \frac{f'(x_1)}{\sqrt{1+|f'|^2}}\right)= -\left[ |\nabla u(x_1,f(x_1))|^2\right]^\pm
\end{equation}
in a weak sense where
$$
\left[ |\nabla u|^2\right]^\pm=|\nabla u(x_1,f(x_1)^+)|^2-|\nabla u(x_1,f(x_1)^-)|^2,
$$
where 
$$
\nabla u(x_1, f(x_1)^\pm)=\lim_{x_2\to f(x_1)^\pm}\nabla u(x_1,x_2).
$$

The final type of variation we will do is to change the position of the crack-tip. We assume that $(u,\Gamma_u)$
is $\epsilon-$close to a crack-tip - this includes the assumptions that the crack-tip of $(u,\Gamma_u)$ is located at the origin
and that $(u,\Gamma_u)$ is a minimizer. By minimality it follows that
$$
J(u,\Gamma_u)\le J(v,\Gamma_v)
$$
for any pair $(v,\Gamma_v)$ where $\Gamma_v$ is some arc in $B_1$ and the crack tip of $\Gamma_v$ is at some point $x^1\ne 0$. This third type of
variation has previously been done in the tangential direction of the crack when the crack consists of a half line. We will do such tangential variations 
of the crack-tip in section \ref{PrelEstSec}.
In section \ref{OrthVarSec} we will make comparisons with functions
that has the crack-tip slightly moved in the in the direction orthogonal the crack-tip, see the right figure below. These variations are, to our knowledge, entirely new and
needed to show that the curvature vanishes at the crack-tip. These variations require much calculation but are in principle trivial
once one has good enough asymptotic information about the solution in a neighborhood of the crack-tip. 

It should be remarked that the notion of ``tangential'' and ``orthogonal'' variations are somewhat imprecise since we do not 
know that the crack-tip has a well defined tangent in section \ref{PrelEstSec} and section \ref{OrthVarSec}. In practice we will make these variations 
only for minimizers $(u,\Gamma_u)$ that are $\epsilon-$close to a crack-tip. We will interpret the tangent of the straight approximating crack as
an approximating tangent of $\Gamma_u$ in the ball $B_1$ and make the variation in the approximate tangent's direction. It would be more appropriate to 
talk about ``variations in the direction of the approximating tangent in the ball $B_1$'' and similarly for orthogonal variations; but we have opted for 
the shorter phrase ``tangential variations.''

\vspace{3mm}

\begin{center}
\includegraphics[height=4.5cm,width=9cm]{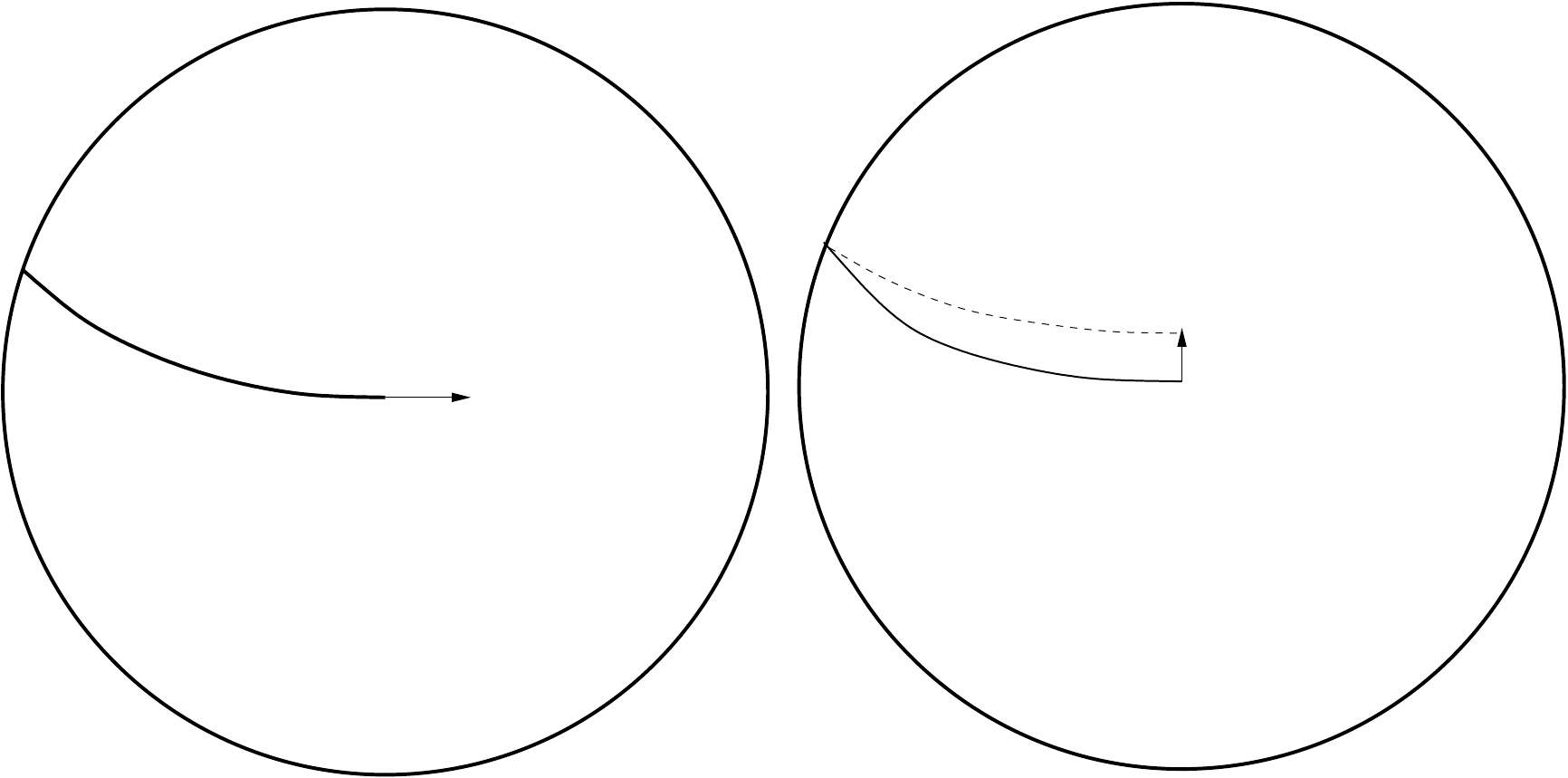}
\end{center}

\vspace{3mm}

{\bf Figure 2.} {\sl In the ``tangential variations'' of the crack-tip we extend the crack as in the left picture. In the ``orthogonal varations'' of the
crack-tip we compare the energy of the crack to the energy of a crack that has the crack-tip moved slightly in the orthogonal direction of the
crack tip. This is shown in the right picture with the comparison crack represented by the dashed line.}

\vspace{2mm}

{\bf Previous versions of the article and a thanks.} A manuscript claiming the same results where published by the authors
on arxive.org. That version contained a sign mistake in the linearized equations. With the correct sign in the equation the 
linearized system has a solution $(\mathfrak{z},\mathfrak{h})$ that is not Lipschitz in $\mathfrak{h}$ (see equations (\ref{DefOfW}) and (\ref{DefOfh})). We would like to thank an anonymous referee for pointing this out to 
us. This paper corrects this mistake. However, the methods published on arxive.org goes through with minor changes even with the right sign. 
The main new development in the arxive.org version of the paper was the new orthogonal variations that allows us to 
get rid of the first term in the expansion of the linearized system. That part of the paper goes through without changes, see Section \ref{OrthVarSec}. 

The first author was supported by Swedish Research 
Council grant (grant number: VR-2016-0363). The second author has been partially supported by National Science Foundation of China (grant number: 11650110437). 

\subsection{Background.}\label{Sec:BackG}

In this sub-section we gather some known results and also try to situate the present research within the field. As we already remarked, the
question of existence of minimizers is already established in different settings \cite{AmbExist} or \cite{DMSActa}. The next step is
usually to prove that the minimizer is more regular that the generic function in the space we minimize in. For the Mumford-Shah problem
this means to show two things, that the minimizing function $u$ is smooth in $B_1\setminus \Gamma_u$ and that the free discontinuity set
$\Gamma_u$ is smooth (except at some controllable set of singular points). The more important, and difficult, problem is to establish that $\Gamma_u$ is smooth. A typical
result in this direction is.

\begin{thm}\label{Ambrosio}
There exists a $\gamma_0>0$ such that if $(u,\Gamma)$ be a minimizing pair of the Mumford-Shah functional and, for some $\gamma<\gamma_0$,
\begin{equation}\label{nablausmall}
\frac{1}{s}\int_{B_s\setminus \Gamma} |\nabla u|^2 dx< \gamma
\end{equation}
and
$$
\Gamma \cap B_1(0)\subset \{(x_1,x_2);\; |x_2|< s^2\gamma\}
$$
then $\Gamma\cap B_{s/2}(0)$ is a $C^{1,1/4}-$graph, of a function $f\in C^{1,1/4}$,
$$
\Gamma=\{(x_1,f(x_1)); \; x_1\in (-s/2, s/2)\}
$$
and there is a universal constant $C$ such that 
\begin{equation}\label{NastyyRHS}
\|f\|_{C^{1,1/4}(B_{s/2})(0)}\le C s^{\frac{1}{4}}\left(1+ \frac{\gamma}{s}\right)^{\frac{1}{2}}.
\end{equation}

\end{thm}

\vspace{1mm}

{\bf Remark:} {\sl Since we are interested in the regularity at the crack-tip the above theorem allows us to assume that the free discontinuity 
$\Gamma_u$ is a $C^{1,1/4}$ graph in $B_1\setminus B_{1/2}$ with small norm. In particular, if $(u,\Gamma_u)$ is $\epsilon-$close to a crack tip with small
enough $\epsilon$ then, by Theorem \ref{Ambrosio}, $(\Gamma_u\cap B_{1/2})\setminus B_{1/4}$ will be a $C^{1,1/4}-$graph. By rescaling 
$u(x)\mapsto \sqrt{2}u(x/2)$ we get a new solution $(u_{1/2},\Gamma_{u_{1/2}})$ in $B_1$ that is $C\epsilon-$close to a crack-tip and $\Gamma_{u_{1/2}}$
is a $C^{1,1/4}$ graph in $B_1\setminus B_{1/2}$.} 

\vspace{1mm}

Theorem \ref{Ambrosio} is a slightly weaker version of the regularity result in \cite{AFPBook} and is close to the state of the art
regularity theory for the Mumford-Shah problem in $\R^2$. The authors of \cite{AFPBook} shows a stronger result in $\R^n$, but
the above theorem is good enough for our purposes, and easier to formulate.

The first thing to notice about Theorem \ref{Ambrosio} is that the condition (\ref{nablausmall}) excludes the crack-tip. The conclusion of the
theorem includes the statement that $B_{s/2}\setminus \Gamma_u$ is disconnected which isn't true for the crack-tip. The problem is that at a
crack-tip $x^0$
\begin{equation}\label{IvarwouldKnow}
u_r(x)=\frac{u(rx+x^0)}{\sqrt{r}}\to \sqrt{\frac{2}{\pi}}r^{1/2}\sin((\phi+\phi_0)/2)
\end{equation}
as $r\to 0$ through some sub-sequence, see \cite{B}) or Corollary \ref{BlowUpAtCarackTipCor}. This means that, for small $r$, $u_r$ will not 
satisfy (\ref{nablausmall}) for small $\gamma$.

We would like to remark that the constant $\sqrt{\frac{2}{\pi}}$ in the limit solution in (\ref{IvarwouldKnow}) expresses the right
balance between the surface and Dirichlet energy in the functional and is uniquely determined by the functional, see Lemma \ref{VariationsInTang}.

Also the estimate (\ref{NastyyRHS}) is not good enough to analyze points close to the crack-tip. In particular, the right hand side in the
estimate (\ref{NastyyRHS}) includes the term $\gamma$ which measures the size of $\|\nabla u\|_{L^2}^2$. However, from the
Euler-Lagrange equations (\ref{PointwiseCurv}) we know that the geometry of $\Gamma_u$ is determined by the difference
$[|\nabla u|^2]^\pm$. If we have symmetry and cancellation in $[|\nabla u|^2]^\pm$ then the curvature of $\Gamma_u$ might be very small even
though $\|\nabla u\|_{L^2}^2$ is not smaller than $\gamma$. At the crack-tip we know that, after a rotation of the coordinate system,
$$
u\approx \sqrt{\frac{2}{\pi}}r^{1/2}\sin(\phi/2).
$$
In particular, we almost have symmetry of the minimizer $u$ close to the crack-tip. This opens up for a more refined analysis of the
regularity of $\Gamma_u$ close to the crack-tip.

We will indeed use this symmetry effect in our proof - but it is hidden away after the obscure equation (\ref{expressionNablaPi}) in the middle of a
long and technical calculation - but that is one of the main ideas of the paper.

There are other regularity proofs for the crack-tip in \cite{B} and \cite{DBook}. However none of them are using the symmetry close to the crack tip
and their proofs does not include the regularity of the free discontinuity sets up to the crack-tip.

\subsection{Main Results.}

Our main result is:

\begin{Main}\label{MainResult1}{\sc [See Theorem \ref{C2alpha}.]} 
There exists an $\epsilon_0>0$ such that if $(u,\Gamma)$ is $\epsilon$-close to a crack-tip solution for some $\epsilon\le \epsilon_0$
 then $\Gamma_u$ is $C^{2,\alpha}$ at the crack-tip for every $\alpha< \alpha_2-3/2$, $\alpha_2\approx 1.889$.

 This in the sense that the tangent at the crack-tip is a well defined line, which we may assume (after rotating the coordinate system) 
 to be $\{(x_1,0);\; x_1\in \R\}$,
 and there exists a constant $C_\alpha$ such that
 $$
 \Gamma_u\subset \left\{(x_1,x_2);\; |x_2|<C_\alpha \epsilon |x_1|^{2+\alpha},\; x_1<0\right\}.
 $$
 Here $C_\alpha$ may depend on $\alpha$ but not on $\epsilon<\epsilon_0$.
\end{Main}

We interpret this as a point-wise $C^{2,\alpha}$ result at the crack-tip. To actually show that $\Gamma$ is a $C^{2,\alpha}$ manifold up to the crack-tip
would require an argument, similar to the argument in \cite{AFPBook} or the main argument here, taking the symmetry into consideration
at points close to the crack-tip. We do not expect any difficulties in proving such a result. However, that would add to the length of an
already long and technical paper. We also believe that the point-wise regularity at the crack-tip, and the techniques leading up to this result,
is of greater importance in order to analyze the quasi-static crack growth.

The proof of the main result consists of linearizing the Euler-Lagrange equations close to a crack-tip. Solutions of the linearized equations can
be expressed as a series of simple solutions all of which are homogeneous except one. The possible homogeneities of the solutions are $\alpha_1=1/2$ and 
the solutions  $\alpha_k\ge 1$, for $k=2,3,4...$, of the following equation
\begin{equation}\label{wobbly}
\tan\left( \alpha \pi\right)=-\frac{2}{\pi}\frac{\alpha}{\alpha^2-\frac{1}{4}}.
\end{equation}
The linearized equations also have a non-homogeneous solution $(\mathfrak{z},\mathfrak{h})$ (see equations (\ref{DefOfW}) and (\ref{DefOfh})). We analyze the linearized system in the Appendix.

As always when one uses linearization it follows that the solution $(u,\Gamma)$ is almost as regular as the solution $(v,f)$
to the linearized system. However, since the linearized problem to the Mumford-Shah problem allows the solution $(\mathfrak{z},\mathfrak{h})$
where $\mathfrak{z}\notin C^{1/2}$ and $\mathfrak{h}\notin C^{0,1}$ one cannot directly draw any conclusions regarding the $C^{1,\alpha}$
regularity of $(u,\Gamma)$ from the linearization - as a matter of fact one cannot even state that $\Gamma$ is a graph close to 
the crack tip from the linearization procedure. 

The natural way to get around this difficulty is to prove that if we linearize a sequence of solutions $(u^j,\Gamma_j)$ to the Euler-Lagrange equations for 
the Mumford-Shah problem we get a solution to the linearized system (see Theorem \ref{theLinearSyst}) that does not contain  the ``bad'' solution 
$(\mathfrak{z},\mathfrak{h})$ in its series expansion. That implies that the solutions we get to the linearized system through a linearization of the Mumford-Shah problem is more 
regular than the generic solution to the linearized system.

The argument is somewhat delicate. We will first, in section \ref{Sec:Lin} and section \ref{SceStrongConv}, show that the Mumford-Shah
problem linearize to the linear system in Theorem \ref{theLinearSyst}. In particular that means that if a minimizer $(u^j,\Gamma_j)$ is $\epsilon_j$-close to a 
crack-tip then, in polar coordinates,
\begin{equation}\label{ujExp}
u^j(r,\phi)= \sqrt{\frac{2}{\pi}}r^{1/2}\sin(\phi/2)+\epsilon_j v(r,\phi)+ R_j(r,\phi)
\end{equation}
where $\|\nabla R_j\|_{L^2(B_1\setminus \Gamma_j)}< < \epsilon_j$ and $v$ solves the linearized system. That $v$ solves the linearized system
implies (Theorem \ref{theLinearSyst}) that
\begin{equation}\label{vjExp}
v(r,\phi)=a_0 \mathfrak{z}(r,\phi)+\sum_{k=1}^\infty a_k r^{\alpha_k}\cos(\alpha_k \phi)+\sum_{k=1}^\infty b_k r^{k-1/2}\sin((k-1/2)\phi).
\end{equation}
The main step in proving $C^{2,\alpha}$ regularity at the crack tip is to show that the $\mathfrak{z}-$term does not appear in the expansion of $v$. 
This has to follow from the assumption that $(u^j,\Gamma_j)$ minimizes the Mumford-Shah energy. We are therefore led to making variations 
of $(u^j,\Gamma_j)$. But knowing the exact form of $u^j$ in (\ref{ujExp}) and (\ref{vjExp}) this is a matter of calculation. As it turns out, 
the right variations are a new type of variations that we call orthogonal to the crack-tip. We will show, in Section \ref{OrthVarSec}, that if $(u^j,\Gamma_j)$
minimizes the Mumford-Shah energy then the constant $a_0=0$ in (\ref{vjExp}). Once we have shown that the Mumford-Shah energy linearizes to 
a ``nice'' function it is a matter of standard technique to show that minimizers are (almost) as regular as the solution to the linearized problem.

A slight complication arises in that the regularity of the solution to the linearized problem is determined to be $C^{1,1/2}$
by a $r^{3/2}\sin(3\phi/2)$ term in its asymptotic development. However, the $r^{3/2}\sin(3\phi/2)$ term will not really affect the 
regularity of the free discontinuity set since it is an odd function in $x_2$. This means that we need a further investigation in Section \ref{SecC2alpha}
in order to prove $C^{2,\alpha}$ regularity at the crack-tip. The analysis is rather subtle and we will need to estimate the 
third order term in the asymptotic expansion in order to prove full $C^{2,\alpha}$ regularity. From a technical point of view the argument in 
Section \ref{SecC2alpha} will be very close to the argument in the previous sections - however some details change.


\subsection{Notation and Conventions.}

We will use several notational conventions in this article. We will freely switch from Cartesian coordinates $(x_1,x_2)$ to polar coordinates $(r,\phi)$.
The one dimensional Hausdorff measure will be denoted $\H^1$. We will often use $\Gamma_0=\{(x_1,0);\; x_1\in (-1,0)\}$. The projection operator
$\Pi(u,s)$ is defined in Definition \ref{ProjDef}. In situations when we have some minimizer $(u,\Gamma_u)$ and consider some multi-valued function,
almost always $\lambda r^{1/2}\sin(\phi/2)$, so that we need to define a branch cut then we assume that the branch cut is along $\Gamma_u$. We will also
use $\nu$ for the normal of a given set. At times we will denote the upper and lower normal of $\Gamma_u$ by $\nu^\pm$.

We will at time refer to the space $H^{1/2}(\partial B_1(0)\setminus \{y\})$ where $y\in \partial B_1(0)$. When we do this we mean the space of all traces on
$\partial B_1(0)$ of Sobolev functions $v\in W^{1,2}(B_1(0)\setminus \Gamma_0)$. We also remark that this space is equivalent to the space of traces of
Sobolev functions $v\in W^{1,2}(B_1(0)\setminus \Gamma_u)$ for any one dimensional set $\Gamma_u$ that is $C^1$ close to $\partial B_1(0)$, non-tangential to
$\partial B_1(0)$ and $\Gamma_u \cap \partial B_1(0)=\{y\}$. This implies, using the remark after Theorem \ref{Ambrosio}, that we can talk about $H^{1/2}(\partial B_1(0)\setminus \Gamma_u)$
for any minimizer $(u,\Gamma_u)$ that is $\epsilon-$close to a crack-tip. We will also use the notation $\P(u,\Gamma)$ for the projection 
of $u$ onto solutions of the linearized system, see Definition \ref{DefProjOnLin}.

We will also use $H^{-1/2}(\Gamma)$, the dual space to $H^{1/2}(\Gamma)$, where $\Gamma$ is a Lipschitz curve. We will identify a function 
$v\in H^{-1/2}(\Gamma)$ with a divergence free vector field $\eta$. In particular, if $w\in H^{1/2}(\Gamma)$ then $w$ is the trace of some 
$W^{1,2}$ function that we still denote $w$. We may identify the pairing $\langle w,v \rangle_{(H^{1/2}, H^{-1/2})}$ as the integral
$\int_{B_1\setminus \Gamma} \nabla w \cdot \eta= \int_{\partial (B_1(0)\setminus \Gamma)}w \eta\cdot \nu$ where $\nu$ is the normal of $B_1\setminus \Gamma$.

\section{Preliminary Estimates.}\label{PrelEstSec}

In this section we gather various preliminary estimates about minimizers of the Mumford-Shah problem. Most of the results are well known or simple wherefore we 
have chosen to only sketch some of the proofs. The reader is advised to look for further details in the references if any result in this section is
unfamiliar.

We begin with a quantitative version of the result that the constant in the blow-up limit (\ref{IvarwouldKnow}) is $\sqrt{\frac{2}{\pi}}$.
A similar calculation is found in \cite{DBook}.

\begin{lem}\label{VariationsInTang}
 There exist constants $\epsilon_0>0$ and $C_0>0$ such that if $0<\epsilon\le \epsilon_0$
 then there is a pair $(v,\Gamma_v)$ such that
 \begin{equation}\label{HomogBdryCond}
 v(x)=\sqrt{\frac{2}{\pi}}(1+\epsilon)r^{1/2}\sin(\phi/2) \textrm{ on }\partial B_1(0),
 \end{equation}
 and
 $$
 J(v,\Gamma_v)\le J\left(\sqrt{\frac{2}{\pi}}(1+\epsilon)r^{1/2}\sin(\phi/2),\Gamma_0\right)-C_0\epsilon^2,
 $$
 where 
 \begin{equation}\label{DefoFGamma0}
  \Gamma_0=\{(x_1,x_2);\; x_1\le 0,\; x_2=0\}.
 \end{equation}
\end{lem}
{\bf Remark:} Notice that the Lemma states that a minimizer in $B_1(0)$ with boundary conditions of the form (\ref{HomogBdryCond})
can only be a homogeneous extension of the boundary data if $\epsilon=0$. In particular, if $(u,\Gamma_u)$ is $\epsilon-$close
to a crack-tip with $\epsilon=0$ then
$$
v(x)=\sqrt{\frac{2}{\pi}}r^{1/2}\sin(\phi/2).
$$

\vspace{2mm}

\textsl{Proof of Lemma \ref{VariationsInTang}:} In this proof we let $\mu=c\epsilon$ where $c$ is a small constant. We will also use the notation
$O(\mu^2)$ for some function whose absolute value is smaller than $C_1\mu^2$ for some fixed constant $C_1$ - that is $O(\mu^2)$ will be a uniform
{\sl ``big Oh notation''}. From the proof we will see that $C_1$ will only depend on the derivatives of $r^{1/2}\sin(\phi/2)$ and is thus universal.

We will assume that $\epsilon>0$ for definiteness. 
The proof consists of showing that, when $\epsilon>0$, there exists a competitor for minimality
with strictly larger (under inclusion) free discontinuity set. If $\epsilon<0$ then an analogous argument works where one shortens the
length of the crack instead.

We consider the ball $\hat{B}=B_{1-\mu}(\mu e_1)\subset B_1(0)$ and use the notation $\Gamma_0=\{(x_1,x_2);\; x_2=0 \textrm{ and }x_1\le 0\}$ and
$$
p(r,\phi)=\sqrt{\frac{2}{\pi}}(1+\epsilon)r^{1/2}\sin(\phi/2).
$$
We will show that 
$$
\int_{\hat{B}}|\nabla p|^2+\H^1(\Gamma_0\cap \hat{B})=\int_{B_{1-\mu}(0)}|\nabla p|^2+\H^1(\Gamma_0\cap B_{1-\mu}(0))-\mu+ O(\mu^2)=
$$
\begin{equation}\label{StarInSomeCrap}
=(1+\epsilon)^2(1-\mu)+(1-2\mu)+O(\mu^2).
\end{equation}
To derive (\ref{StarInSomeCrap}) we calculate the derivative
\begin{equation}\label{Star19thMay}
\frac{\partial }{\partial t}\int_{B_{1-\mu}(t\mu e_1)}|\nabla p|^2\Big\lfloor_{t=0}=\int_{B_{1-\mu}(t)}\cos(\phi)|\nabla p|^2=0.
\end{equation}
Using a second order Taylor expansion in $t$ of $\int_{B_{1-\mu}(t\mu e_1)}|\nabla p|^2$ and using that the linear term is zero by (\ref{Star19thMay}) it 
follows that
\begin{equation}\label{GrisOrvar}
\int_{\hat{B}}|\nabla p|^2=\int_{B_{1-\mu}(0)}|\nabla p|^2+O(\mu^2).
\end{equation}
The calculation in (\ref{StarInSomeCrap}) is just (\ref{GrisOrvar}) together with an explicit calculation $\H^1(\Gamma_0\cap \hat{B})=1-\mu$.

Next we observe that at $\left( r\cos(\phi)+\mu, r\sin(\phi)\right)\in\partial \hat{B}$
$$
p\left( r\cos(\phi)+\mu, r\sin(\phi)\right)=p(r,\phi)+\frac{\partial p(r,\phi)}{\partial x_1}\mu +O(\mu^2)=
$$
$$
=\sqrt{\frac{2}{\pi}}(1+\epsilon -\frac{\mu}{2})r^{1/2}\sin(\phi/2)+O(\mu^2),
$$
where we used a Taylor expansion again.

Letting $\hat{r}$ and $\hat{\phi}$ be polar coordinates centered at the center of $\hat{B}$ we may define the harmonic function in $\hat{B}\setminus \{(\hat{r},\pi);\; \hat{r}>0\}$
$$
w(\hat{r},\hat{\phi})=\sqrt{\frac{2}{\pi}}(1+\epsilon -\frac{\mu}{2})\hat{r}^{1/2}\sin(\hat{\phi}/2)+O(\mu^2)
$$
such that $w=p$ on $\partial \hat{B}$. Moreover we choose $\Gamma_w=\{(\hat{r},\pi);\; \hat{r}>0\}$. Then
$$
\int_{\hat{B}\setminus \Gamma_w}|\nabla w|^2+\H^1(\Gamma_w)=\left(1+\epsilon-\frac{\mu}{2}\right)^2(1-\mu)+(1-\mu)+O(\mu^2).
$$

We can conclude that
$$
\int_{\hat{B}\setminus \Gamma_0}|\nabla p|^2+\H(\Gamma_0\cap \hat{B})-\left(\int_{\hat{B}\setminus \Gamma_w}|\nabla w|^2+\H^1(\Gamma_w) \right)=
$$
\begin{equation}\label{IveLostControlAgain}
=(1+\epsilon)^2(1-\mu)+(1-2\mu)-\left(1+\epsilon-\frac{\mu}{2}\right)^2(1-\mu)-(1-\mu)+O(\mu^2)=
\end{equation}
$$
=\mu\epsilon -\frac{5\mu^2}{4}+O(\mu^2)>0
$$
if $0<\mu<c\epsilon$ and $\epsilon$ is small enough and $c$ is a small constant. The calculation (\ref{IveLostControlAgain}) shows that the function
$$
v(x)=\left\{
\begin{array}{ll}
 p(x) & \textrm{in }B_1(0)\setminus \hat{B} \\
 w(x) & \textrm{in } \hat{B}
\end{array}
\right.
$$
satisfies the conclusion of the lemma.
\qed

A simple perturbation result shows that any minimizer, not just minimizers with the limited boundary conditions
specified by (\ref{HomogBdryCond}), satisfies a similar estimate. The formulation of this corollary differs somewhat from the formulation of the
preceding Lemma - but the proofs are similar.

\begin{cor}\label{CPiLem}
 Let $(u,\Gamma_u)$ be a minimizer to the Mumford-Shah problem and assume that
 \begin{equation}\label{CloseToHomog}
 \int_{B_1\setminus \Gamma_u}\left|\nabla \left(u-\lambda\sqrt{\frac{2}{\pi}} r^{1/2}\sin(\phi/2)\right)\right|^2\le \epsilon^2
 \end{equation}
 for some small $\epsilon>0$.

 Then there exists a constant $C_\Pi$ such that
 \begin{equation}\label{PiCoeffEst}
  1-C_\Pi\sqrt{\lambda \epsilon} \le |\lambda|\le 1+C_\Pi \sqrt{\lambda \epsilon}.
 \end{equation}
\end{cor}
\textsl{Proof:} It is enough to show this for small $\epsilon$. We may also, for the same reasons as in the remark after Theorem \ref{Ambrosio}, 
assume that $\Gamma_u$ is a $C^{1,1/4}$ graph in $B_1\setminus B_{1/2}$ with norm less than $1/4$ if $\epsilon$ is small enough.

If $u$ satisfies (\ref{CloseToHomog}) then, by the trace theorem,
$$
\left\| u\lfloor_{\partial B_1}-\lambda \sqrt{\frac{2}{\pi}}r^{1/2}\sin(\phi/2)\right\|_{H^{1/2}(\partial B_1(0)\setminus\Gamma_u)}\le C\epsilon,
$$
where the constant $C$ can be assumed to be universal if $\epsilon$ is so small that 
$$
\|\Gamma_u \|_{C^{1,1/4}(B_1\setminus B_{1/2})}\le 1/4.
$$

Let us denote
$$
w=u\lfloor_{\partial B_1}-\lambda \sqrt{\frac{2}{\pi}}r^{1/2}\sin(\phi/2) \quad \textrm{ on }\quad \partial B_1(0).
$$

We know, from Lemma \ref{VariationsInTang}, that there exists a function $v$ such that
$$
v= \lambda\sqrt{\frac{2}{\pi}} r^{1/2}\sin(\phi/2) \quad\textrm{ on }\partial B_1(0)
$$
and
$$
J(v,\Gamma_v)\le J\left(\lambda\sqrt{\frac{2}{\pi}} r^{1/2}\sin(\phi/2),\Gamma_0\right)-C_0 (\lambda-1)^2.
$$

If we extend $w$ to a harmonic function solving
$$
\begin{array}{ll}
 \Delta w=0 & \textrm{ in }B_1(0)\setminus \Gamma_v \\
 \frac{\partial w}{\partial \nu}=0 & \textrm{ on }\Gamma_v.
\end{array}
$$

Then
$$
J(v+w,\Gamma_0)\le J(v,\Gamma_v)+C\epsilon^2 +2\int_{B_1(0)\setminus \Gamma_v}\nabla v\cdot \nabla w\le
$$
$$
\le J\left(\lambda\sqrt{\frac{2}{\pi}} r^{1/2}\sin(\phi/2),\Gamma_0\right)-C_0 (\lambda-1)^2+C_1\lambda \epsilon.
$$

We may also use (\ref{CloseToHomog}) to calculate
$$
J(u,\Gamma_u)\ge J\left(\lambda\sqrt{\frac{2}{\pi}} r^{1/2}\sin(\phi/2),\Gamma_0\right)- C_2\lambda \epsilon,
$$
where we also used that $\H^1(\Gamma_0)\le \H^1(\Gamma_u)$.

Since $J(u,\Gamma_u)\le J(v+w,\Gamma_0)$ we can conclude that
$$
C_0(\lambda-1)^2\le \left(C_1 +C_2\right)\lambda \epsilon,
$$
which implies the Corollary.
\qed

Next we define an operator that projects minimizers on the space of homogeneous functions. The right projection operator is
$\Pi$ defined in the following definition. $\Pi$ projects the solutions onto the manifold of all rotations and multiples of the $1/2-$homogeneous 
minimizers.

\begin{definition}\label{ProjDef}
 Let $(u,\Gamma)$ be a minimizer to the Mumford-Shah problem in $B_s$, then we define the
 projection operator $\Pi(u,\Gamma,s)=\Pi(u,s)$, or at times just $\Pi(u)$ if 
 $s=1$, to be the minimizer of the following expression
 \begin{equation}\label{THeDefofPi21}
 \min_{p\in P(\Gamma)}\int_{B_s\setminus \Gamma}\left|\nabla \left(u-p\right)\right|^2
 \end{equation}
 where
 $$
 P(\Gamma)=\left\{ p;\, p(r,\phi)=\sqrt{\frac{2}{\pi}} \lambda r^{\frac{1}{2}}\sin\left(\frac{\phi}{2}+\phi_0\right);\, \phi_0\in (-\pi,\pi],\; \lambda\in \R \right\},
 $$
 the $\Gamma$ in $P(\Gamma)$ indicates the the branch-cut of the functions $p\in P(\Gamma)$ are made at $\Gamma$.
\end{definition}

{\bf Remark:} {\sl Since the minimization in 
(\ref{THeDefofPi21}) is a continuous minimization problem in two variables 
$\phi_0$ and $\lambda$. Since the minimization problem is clearly coercive in 
$\lambda$ and $\phi_0$ may be considered to contained 
compact unit circle it follows that a minimizer exist.}

We think of $\Pi$ as a projection operator since
$$
\int_{B_1(0)\setminus \Gamma_u}\left( \nabla \left(u-\Pi(u,1) \right)\right) \cdot \nabla v=0
$$
for any $v\in \mathcal{P}$. In particular we may conclude, by using minimality and make a variation in $\lambda$, that
\begin{equation}\label{OrthtoSin}
\int_{B_1(0)\setminus \Gamma_u}\left( \nabla \left(u-\Pi(u,1) \right)\right) \cdot \nabla \left( r^{1/2}\sin(\phi/2)\right)=0
\end{equation}
and similarly, by making a variation in $\phi_0$, that
$$
\int_{B_1(0)\setminus \Gamma_u}\left( \nabla \left(u-\Pi(u,1) \right)\right) \cdot \nabla \left( r^{1/2}\cos(\phi/2)\right)=0.
$$

\vspace{2mm}

Next we formulate a lemma to control the Hausdorff measure of the free discontinuity set.

\begin{lem}
 Assume that $(u,\Gamma)$ is $\epsilon-$close to a crack-tip. Then there is a constant $C_\H$ such that
 \begin{equation}\label{HausdorffEst}
  \H^1(\Gamma\cap B_{1})\le 1+C_\H\epsilon.
 \end{equation}
\end{lem}
\textsl{Proof:} The proof is very simple. In particular, we may rotate the coordinate system so that $\Gamma_u\cap \partial B_1(0)=(-1,0)$ and 
compare the energy of $(u,\Gamma)$ in $B_1(0)$ to the competing pair $(w,\Gamma_0)$ where $\Gamma_0=\{(x_1,0);\; x_1\in (-1,0)\}$ and $w$ is 
harmonic in $B_1(0)\setminus \Gamma_0$, has zero Neumann boundary data on $\Gamma_0$ and equals $u$ on the boundary $\partial B_1(0)$. 
We may write 
$$
w(x)=h(x)+\lambda \sqrt{\frac{2}{\pi}}r^{1/2}\sin\left(\frac{\phi}{2}\right),
$$
where $h$ is harmonic in $B_1\setminus \Gamma_0$ and $\|\nabla h\|_{L^2(B_1\setminus \Gamma_0)}\le C\epsilon$.
We may also choose $\lambda$ so that $\Pi(u)=\lambda\sqrt{\frac{2}{\pi}}r^{1/2} \sin(\phi/2)$ and therefore 
$u-\lambda\sqrt{\frac{2}{\pi}}r^{1/2} \sin(\phi/2)$ is orthogonal to 
$r^{1/2} \sin(\phi/2)$ as in (\ref{OrthtoSin}).
We would like to remind the reader that we will consider the function 
$\lambda \sqrt{\frac{2}{\pi}}r^{1/2}\sin(\phi/2)$ to have different branch cuts in different formulas, for instance when $\lambda \sqrt{\frac{2}{\pi}}r^{1/2}\sin(\phi/2)$
appears in an integral such as (\ref{OrthtoSin}) the branch-cut is chosen 
along $\Gamma_u$ but if it appears in an integral over 
$B_1(0)\setminus \Gamma_0$ then the branch-cut will be considered to be over 
$\Gamma_0$.

From minimality of $(u,\Gamma_u)$ we can conclude that
$$
\int_{B_1(0)\setminus \Gamma_u}|\nabla u|^2+\H^1(\Gamma_u \cap B_1)\le \int_{B_1\setminus \Gamma_0}|\nabla w|^2+\H^1(\Gamma_0\cap B_1)\le
$$
\begin{equation}\label{bergs}
\le \int_{B_1\setminus \Gamma_0}\left|\nabla \left(\lambda\sqrt{\frac{2}{\pi}}r^{1/2}\sin(\phi/2) \right)\right|^2+ 
2\int_{B_1\setminus \Gamma_0}\nabla h\cdot \nabla \left(\lambda\sqrt{\frac{2}{\pi}}r^{1/2}\sin(\phi/2) \right)+
\end{equation}
$$
+\int_{B_1\setminus \Gamma_0}|\nabla h|^2 +1\le \int_{B_1\setminus \Gamma_0}\left|\nabla \left(\lambda\sqrt{\frac{2}{\pi}}r^{1/2}\sin(\phi/2) \right)\right|^2+ C\epsilon+1
$$
where we used $\H^1(\Gamma_0)=1$ in the last estimate.
Since $(u,\Gamma)$ is $\epsilon-$close to a crack-tip it follows that
$$
J(u,\Gamma_u)=\int_{B_1\setminus\Gamma_u}\left| \nabla \left(u-\Pi(u)+\Pi(u)\right)\right|^2+\H^1(\Gamma_u)=
$$
$$
=\int_{B_1\setminus\Gamma_u}\left| \nabla \left(u-\Pi(u)\right)\right|^2+
2\int_{B_1\setminus\Gamma_u}\nabla (u-\Pi(u))\cdot \nabla \Pi(u)+\int_{B_1\setminus\Gamma_u}\left| \nabla \left(\Pi(u)\right)\right|^2 +
$$
\begin{equation}\label{rummen}
+\H^1(\Gamma_u)\ge\int_{B_1\setminus\Gamma_u}\left| \nabla \left(\lambda\sqrt{\frac{2}{\pi}}r^{1/2}\sin(\phi/2)\right)\right|^2 +\H^1(\Gamma_u),
\end{equation}
where we used that $\Pi(u)=\lambda\sqrt{\frac{2}{\pi}}r^{1/2}\sin(\phi/2)$ and that $\nabla(u-\Pi(u))$ is orthogonal to $\Pi(u)$ as in (\ref{OrthtoSin}).
%
The estimate (\ref{bergs}) and the estimate ending in (\ref{rummen}) clearly implies that
$$
\H^1(\Gamma_u \cap B_\mu)\le \H^1(\Gamma_0\cap B_1)+C\epsilon.
$$
\qed

\begin{cor}\label{L2boundsCor}
 Assume that $(u,\Gamma)$ is $\epsilon-$close to a crack-tip. Then there is a line, $l$, from the origin to $\partial B_1(0)$
 such that
 $$
 \sup_{x\in \Gamma_u} \textrm{dist}(x,l)\le C \sqrt{\epsilon}.
 $$
\end{cor}
\textsl{Proof:} We may rotate the coordinate system so that $(-1,0)\in \Gamma_u$. Let $(x_1,x_2)\in \Gamma_u$ then
\begin{equation}\label{sdfg}
\sqrt{(1+x_1)^2+x_2^2}+\sqrt{x_1^2+x_2^2}\le \H^1(\Gamma_u\cap B_{1})\le 1+C\epsilon
\end{equation}
the left inequality is obviously true for any 
$(x_1,x_2)\in \Gamma_u\cap B_1(0)$ since the left side is the length
of the two straight line segments form $(-1,0)$ to $(x_1,x_2)$ and then from 
$(x_1,x_2)$ to the origin whereas $\Gamma_u$ connects the same points without necessarily consisting of straight line segments.
The left side in (\ref{sdfg}) attains its minimal value for fixed $x_2$ when $x_1=-1/2$. Thus
$$
\sqrt{1+4x_2^2}\le \H^1(\Gamma\cap B_{1})\le 1+C\epsilon.
$$
This implies that $|x_2|\le C\sqrt{\epsilon}$. Similarly, for fixed $x_1$ the minimal value of (\ref{sdfg})
is attained for $x_2=0$ which implies that $x_1<C\sqrt{\epsilon}$. This implies the corollary. \qed

\begin{lem}\label{FirstBadC1}
 For each $\tau>0$ there exists an $\epsilon_\tau>0$ such that if $(u,\Gamma)$ is $\epsilon-$close to a crack-tip for some
 $\epsilon<\epsilon_\tau$ then $\Gamma$ is a graph in $B_1\setminus B_\tau$. That is
 $$
 \Gamma\cap B_1=\{(x_1,f(x_1));\; x_1\in (-1,-\tau)\}\textrm{ in the set } -1<x_1\le -\tau.
 $$

 Furthermore there exists a modulus of continuity $\sigma$ such that 
 \begin{equation}\label{BadC14Estf}
 \left[f'(x_1)\right]_{C^{1/4}(-1,-\tau)}\le \sigma(\epsilon)
 \end{equation}
 and 
 \begin{equation}\label{BadLipEstu}
 \sup_{B_{1-\tau}(0)\setminus (B_\tau(0)\cup \Gamma_u)} |\nabla (u-\Pi(u))|\le \sigma(\epsilon).
 \end{equation}
\end{lem}
\textsl{Proof:} The proof is a direct application of Theorem \ref{Ambrosio}. We let $\delta>0$ be a small constant, in particular we assume that
$0<\delta<\gamma_0$, where $\gamma_0$ is as in Theorem \ref{Ambrosio}. We also consider a small ball
$B_s(-te_1)$ with $0<s<t$ and $s$ small enough so that
\begin{equation}\label{Dsmall}
\frac{1}{s}\int_{B_s(-t e_1)\setminus \Gamma_0}\left|\nabla \left( \sqrt{\frac{2}{\pi}}r^{1/2}\sin(\phi/2)\right)\right|^2\le \frac{\delta}{2}.
\end{equation}
Notice that we may choose $s = c \delta t$ for some fixed small constant $c$.

Clearly since $u$ is $\epsilon-$close to a cracktip it follows that  if $\epsilon$ is small enough, 
depending only on $s$ and $\delta$, then
\begin{equation}\label{WhyWhy}
\frac{1}{s}\int_{B_s(-t e_1)\setminus \Gamma_u}\left|\nabla u\right|^2\le \delta<\gamma_0.
\end{equation}
Furthermore by Corollary \ref{L2boundsCor}, if $\epsilon$ is small enough,
\begin{equation}\label{Trapped}
\Gamma_u\subset \{(x_1,x_2);\; |x_2|< s^2\delta\}.
\end{equation}
We may conclude that for any $0<\delta<\gamma_0$ then if $\epsilon$ is small enough so that (\ref{Trapped}) and (\ref{WhyWhy}) are satisfied then it follows,
from Theorem \ref{Ambrosio}, that $\Gamma_u$ is a $C^{1,1/4}-$graph in $B_{s/2}(-te_1)$ and
\begin{equation}\label{ElPasoHellHole}
 \left[f'(x_1)\right]_{C^{1/4}(-t-s/2,-t+s/2)}\le Cs^{1/4}\left( 1+\frac{\delta}{s}\right)^{1/2}=
\end{equation}
$$
=
 Ct^{1/4}\delta^{1/4}\sqrt{1+\frac{1}{ct}}\le C\sqrt{1+\frac{1}{c\tau}}\delta^{1/4}
$$
on the set $-t\le \tau<0$.

But $\delta$ is any constant s.t. $0<\delta<\gamma_0$ and therefore (\ref{ElPasoHellHole}) implies that
$|f'(-t)|$ can be made arbitrarily small on the set $\{-t\le -\tau\}$ if $\epsilon$ is small enough.

The estimate (\ref{BadLipEstu}) follows from the fact that $u-\Pi(u)$ is harmonic in $B_{1-\tau}\setminus \Gamma_u$
with $\|\nabla (u-\Pi(u))\|_{L^2}\le \epsilon$. Also the normal derivative of $u-\Pi(u)$ on $\Gamma_u$ will equal the normal 
derivative of $\Pi(u)$ on $\Gamma_u$ which will be less than $C\sigma(\epsilon)$ since $\Gamma_u$ is a $C^{1,1/4}-$graph with 
norm at most $\sigma(\epsilon)$. The estimate (\ref{BadLipEstu}) follows from this and standard regularity for harmonic functions.
\qed

The following monotonicity formula was first proved in \cite{B}.

\begin{lem}\label{BonnetMonFormula}
 Let $(u,\Gamma_u)$ be a minimizer of the Mumford-Shah energy and assume that $\Gamma_u$ is connected in $B_1(0)$. Then
 \begin{equation}\label{MonFunk}
 \frac{1}{r}\int_{B_r(0)\setminus \Gamma_u}|\nabla u|^2
 \end{equation}
 is non-decreasing in $r$. Furthermore the functional in (\ref{MonFunk}) is constant only if $u(x)-u(0)$ is $r^{1/2}-$homogeneous.

 In particular, if $(u,\Gamma_u)$ is $\epsilon-$close to a crack-tip then (\ref{MonFunk}) is non decreasing.
\end{lem}
\textsl{Proof:} For the proof of the first part of the lemma see \cite{B}. The second part follows since, by assumption,
a minimizer that is $\epsilon-$close to a crack-tip is a minimizer and $\Gamma_u\cap B_1(0)$ has one component. \qed

A slight difficulty in dealing with convergence issues for sequences $(u^j,\Gamma_{u^j})$ of minimizers of the Mumford-Shah problem is that since 
$u^j\in W^{1,2}(B_1\setminus \Gamma_{u^j})$ and the free discontinuity set $\Gamma_{u^j}$ is different for each $u^j$
every function $u^j$ in the sequence belongs to a different Sobolev space. This makes it rather subtle to define convergence for the sequence $u^j$.
In order to handle the problems with $\Gamma_{u^j}$ being different for each $j$ we will show convergence for the regular part of the gradient.

\begin{definition}
We will use the notation $\Reg(\nabla v^j)$ to denote the regular part of the measure $\nabla v^j$, when $v^j$ is considered to be
a function of bounded variation in $B_1$. In particular, $\Reg(\nabla v^j)$ will be an $L^2-$function agreeing with $\nabla v^j$ on $B_1\setminus \Gamma_j$.
\end{definition}

\begin{cor}\label{BlowUpAtCarackTipCor}
 If $(u,\Gamma)$ is $\epsilon-$close to a crack-tip then, for any sequence $s_j\to 0$ such that $\frac{u(s_j)}{\sqrt{s_j}}$ converges in $L^2(B_1)$, 
 \begin{equation}\label{BlowupIsCT}
 \lim_{s_j\to 0}\frac{u(s_jx)}{\sqrt{s_j}}=\sqrt{\frac{2}{\pi}}r^{1/2}\sin\left(\frac{\phi+\phi_0}{2}\right),
 \end{equation}
 where the convergence in (\ref{BlowupIsCT}) is also strong in $L^2(B_1)$. Furthermore if $s_j\to 0$ is a sequence such that (\ref{BlowupIsCT})
 holds in $L^2(B_1)$ then $\Reg\left(\nabla \frac{u(s_j x)}{\sqrt{s_j}}\right)$ converges strongly.
\end{cor}
{\bf Remark:} {\sl  In the conclusion of the corollary we do not exclude that the constant $\phi_0$ may depend on the particular sequence $s_j\to 0$.
Later we will see that the rotation $\phi_0$ is actually independent of the sequence $s_j\to 0$.}

\vspace{2mm}

\textsl{Sketch  of the proof of Corollary \ref{BlowUpAtCarackTipCor}:} The proof is essentially contained in \cite{B}. 

Using the notation $u_{s_j}(x)=\frac{u(s_jx)}{\sqrt{s_j}}$ it follows from Lemma \ref{BonnetMonFormula} that 
the integral
$$
\frac{1}{s_j}\int_{B_{s_j}(0)\setminus \Gamma_{u}}|\nabla u|^2=\int_{B_1(0)\setminus \Gamma_{u_{s_j}}}|\nabla u_{s_j}|^2
$$
is non-decreasing, here we also used a simple change of variables in the equality and indicate by a subscript in $\Gamma$
if we consider the free discontinuity set of $u$ or of $u_{s_j}$.

Also by Lemma \ref{BonnetMonFormula} it follows that, for any $t\in (0,1]$,
\begin{equation}\label{100thTime}
\int_{B_1\setminus \Gamma_{u_{s_j}}}|\nabla u_{s_j}|^2 \ge \frac{1}{t}\int_{B_t \setminus \Gamma_{u_{s_j}}}|\nabla u_{s_j}|^2\ge
\lim_{j\to \infty}\int_{B_1\setminus \Gamma_{u_{s_j}}}|\nabla u_{s_j}|^2 .
\end{equation}
Passing to the limit $j\to \infty$ also on the left and in the middle of (\ref{100thTime}) we may conclude that 
if $u_{s_j}\to u_0$ then 
$$
\frac{1}{t}\int_{B_t \setminus \Gamma_{u_{0}}}|\nabla u_{0}|^2
$$
is constant. It follows from Lemma \ref{BonnetMonFormula} that $u_0$ is homogeneous of order $1/2$. Furthermore $u_0$ cannot be identically zero since 
$(u_0,\Gamma_{u_0})$ is a minimizer (by a slight variation of Theorem 2.2 in \cite{B}) and $\Gamma_{u_0}$ is a rectifiable curve starting at  
the origin and,  as is easily seen by the assumption \ref{connectedCurve} in Definition \ref{epsCloseDef}, goes to infinity.

Since $u_0$ is a non-zero harmonic minimizer that is homogeneous of order $1/2$ it 
follows, from Lemma \ref{VariationsInTang} (and the remark thereafter), that
$$
u_0=\sqrt{\frac{2}{\pi}}r^{1/2}\sin\left(\frac{\phi+\phi_0}{2}\right).
$$
This proves (\ref{BlowupIsCT}). 

To see that $\Reg\left(\nabla \frac{u(s_j x)}{\sqrt{s_j}}\right)$ converges strongly we may use that the free discontinuity $\Gamma_{u_{s_j}}$ will 
converge in the Hausdorff metric to $\Gamma_0$. This will imply that for any $x_j\in \Gamma_{u_{s_j}}$ such that $|x_j|\ge \tau>0$ , for $j$
large enough, $(u_{s_j},\Gamma_{u_{s_j}})$ will be a local minimizer of the Mumford-Shah problem in $B_{\tau/2}(x_j)$ with $\Gamma_{u_{s_j}}\cap B_{\tau/2}(x_j)$
being a flat graph. A standard calculation shows that $|\nabla u_{s_j}|\le C|x_j|^{-1/2}$ in $B_{\tau/2}(x_j)$. Furthermore, by 
Lemma \ref{BonnetMonFormula}, for any $\delta>0$
$$
\int_{B_{\delta}(0)\setminus \Gamma_{u_{s_j}}}|\nabla u_{s_j}|^2\le C\delta.
$$
We may conclude that $\nabla u_{s_j}$ converges weakly (by weak compactness), point-wise (by harmonicity) and without concentration (by the above estimates).
Strong convergence follows. \qed

\section{Linearization.}\label{Sec:Lin}

From now on we will assume that we have a sequence of minimizers $(u^j,\Gamma_j)$ that are $\epsilon_j$ close to a crack-tip, $\epsilon_j\to 0$.
We will write
$$
u^j=\Pi(u^j,1)+\epsilon_j v^j,
$$
where $\epsilon_j$ is chosen so that $\|\nabla v^j\|_{L^2(B_1\setminus \Gamma_j)}= 1$. We remark that, by Lemma \ref{FirstBadC1},
\begin{equation}\label{JustifiedByThm11}
\Gamma_j\setminus B_{\tau_j}=\left\{(x,\epsilon_j f_j(x));\; x\in (-1,-\tau_j)\right\}\cap B_1
\end{equation}
where $\tau_j\to 0$ and $f_j\in C^1$.

The main result of this section is that $v^0=\lim_{j\to \infty}v^j$ and $f_0=\lim_{j\to \infty}f_j$ satisfies the following linear system of PDE
\begin{equation}\label{limitsystprop}
   \begin{array}{ll}
    \Delta v^0=0 & \textrm{ in } B_1\setminus\{ x_1<0, x_2=0\} \\
    -\frac{\partial v^0(x_1,0^+)}{\partial x_2}=\frac{1}{2}\sqrt{\frac{2}{\pi}}\frac{\partial}{\partial x_1}\left(\frac{1}{\sqrt{-x_1}}f_0(x_1)\right)
    & \textrm{ for }x_1<0 \\
    -\frac{\partial v^0(x_1,0^-)}{\partial x_2}=-\frac{1}{2}\sqrt{\frac{2}{\pi}}\frac{\partial}{\partial x_1}\left(\frac{1}{\sqrt{-x_1}}f_0(x_1)\right)
    & \textrm{ for }x_1<0 \\
    \frac{\partial^2 f_0(x_1)}{\partial x_1^2}=-
    \sqrt{\frac{2}{\pi}\frac{1}{r}}\left(\frac{\partial v^0(x_1,0^+)}{\partial x_1}+\frac{\partial v^0(x_1,0^-)}{\partial x_1}\right)
     & \textrm{ for }x_1<0,
    \end{array}
  \end{equation}
we will discuss this linear system in more detail in the Appendix, where we will show existence, uniqueness and the regularity of the 
solutions of (\ref{limitsystprop}), under an additional condition on the function $f$.

In order to prove the main regularity result we need to control the first two orders of the asymptotic expansion of a solution $(u,\Gamma)$. 
The first order expansion is, by Corollary \ref{BlowUpAtCarackTipCor}, $\sqrt{\frac{2}{\pi}}r^{1/2}\sin(\phi/2)$. 
Next we show that if we linearize the problem we get a solution to the system of equations analyzed in the 
appendix.

\begin{prop}\label{Linearization} {\sc [Linearization away from the crack-tip.]}
Let $(u^j,\Gamma_j)$ be a sequence of minimizers to the Mumford-Shah problem that are $\epsilon_j-$close to a crack tip for some sequence
$\epsilon_j\to 0$. Furthermore, let
\begin{equation}\label{defofVjinProp}
 v^j(x)=\frac{u^j(x)-\Pi(u^j,1)}{\epsilon_j}.
\end{equation}
Then if (as we will show in Lemma \ref{WhatShouldWeCallThis}) there exists a subsequence, which we still denote $v^j$, such that $v^j\to v^0$ strongly in 
$L^2(B_1(0))$ and $\Reg(\nabla v^j)\to \Reg(\nabla v^0)$ weakly in $L^2(B_b(0)\setminus B_a(0))$ for every $0<a<b<1$. If $f_j$ is the function 
satisfying (see (\ref{JustifiedByThm11}))
$$
\Gamma_j\setminus \{|x_1|< \tau_j\}=\left\{ (x_1,\epsilon_j f_j(x_1));\; x_1\le -\tau_j\right\}
$$
then $f_j\to f_0$ weakly for some function $f_0\in C([-1,0]) \cap W^{1,2}((-b,-a))$ for any $0<a<b<1$. Furthermore $v^0$ and $f_0$ 
satisfies (\ref{limitsystprop}).
\end{prop}
\textsl{Proof:}  For simplicity of notation we will write $\Pi(u^j)$ for $\Pi(u^j,1)$ in this proof. By (\ref{defofVjinProp}) $u^j=\Pi(u^j)+\epsilon_j v^j(x)$, and as we remarked in the beginning of this section
$\epsilon_j$ is chosen such that $\|\nabla v^j\|_{L^2(B_1(0)\setminus \Gamma_{u^j})}=1$.

Pick $\mu_0>0$ and let $j$ be large enough so that 
$\Gamma_{u^j}$ is a graph of $\epsilon_j f_j(x_1)$ in $B_1(0)\setminus B_{\mu_0}(0)$, such a $j$ always exist by Lemma \ref{FirstBadC1}. 
Let $\eta(x)= \psi(x) e_2$ with $\psi\in C^{\infty}_c(B_1(0)\setminus B_{\mu_0}(0))$ and $D_2 \psi(x)=0$ close to $\Gamma_{u^j}$. Making a domain variation 
as in (\ref{DomVar}) we can derive that
$$
0=\int_{B_1(0)\setminus \Gamma_{u^j}} \bigg( \left|\nabla (\Pi(u^j)+\epsilon_j v^j) \right|^2 \frac{\partial \psi}{\partial x_2}-
$$
$$
-2(\nabla (\Pi(u^j)+\epsilon_jv^j)\cdot e_2) (\nabla (\Pi(u^j)+\epsilon_jv^j)\cdot \nabla \psi)\bigg)+
$$
$$
+\int_{-1}^0 \frac{\epsilon_j f'_j(x_1)}{\sqrt{1+\epsilon_j^2 |f'(x_1)|^2}}\frac{\partial \psi(x)}{\partial x_1}=
$$
\begin{equation}\label{zerothorder}
=\int_{B_1(0)\setminus \Gamma_{u^j}} \left(\left|\nabla \Pi\right|^2 \frac{\partial \psi}{\partial x_2}-2(\nabla \Pi\cdot e_2)(\nabla \Pi\cdot \nabla \psi) \right)+
\end{equation}
\begin{equation}\label{firstorder1}
 +\epsilon_j\int_{B_1\setminus \Gamma_{u^j}}\Big( 2\nabla \Pi\cdot \nabla v^j \frac{\partial \psi}{\partial x_2}-2(\nabla \Pi\cdot e_2)(\nabla v^j\cdot \nabla \psi)-
\end{equation}
 \begin{equation}\label{firstorder3}
 -2(\nabla v^j \cdot e_2)(\nabla \Pi\cdot \nabla \psi)\Big)+
\end{equation}
\begin{equation}\label{firstorder2}
 +\epsilon_j\int_{-1}^0 \frac{f_j'(x_1)}{\sqrt{1+\epsilon_j^2 |f_j'(x_1)|^2}}\frac{\partial \psi(x)}{\partial x_1}+
\end{equation}
\begin{equation}\label{secondorder}
 +\epsilon_j^2\int_{B_1(0)\setminus \Gamma_{u^j}}\left(|\nabla v^j|^2\frac{\partial \psi}{\partial x_2}-2(\nabla v^j\cdot e_2)(\nabla v^j\cdot \nabla \psi) \right).
\end{equation}

If we make an integration by parts in (\ref{zerothorder}) we see that the integral in (\ref{zerothorder}) equals
\begin{equation}\label{zerothorderafterintbyparts}
\int_{\Gamma_{u^j}^\pm}\left|\nabla \Pi\right|^2\psi (\nu^\pm \cdot e_2)-2(\nabla \Pi\cdot e_2)(\nabla \Pi\cdot \nu^\pm) \psi.
\end{equation}

Notice that
\begin{equation}\label{expressionNablaPi}
\nabla \Pi(u^j)=\frac{1}{2}\sqrt{\frac{2}{\pi}}\frac{1}{r^{1/2}}\big(-\sin(\phi/2),\cos(\phi/2) \big).
\end{equation}
In particular, since the value of $\phi$ differs by $2\pi$ on $\Gamma^\pm$, we can conclude that
$\nabla \Pi(u^j)\lfloor_{\Gamma^+}=-\nabla \Pi(u^j)\lfloor_{\Gamma^-}$. Since also $\nu^+=-\nu^-$ it follows that
the value of (\ref{zerothorderafterintbyparts}), and therefore (\ref{zerothorder}), is identically zero.

Let us continue to estimate the integral in (\ref{secondorder}). Lemma \ref{FirstBadC1}, together with the normalization 
$\| \Reg(\nabla v^j)\|_{L^2(B_1)}\le 1$, implies that we may estimate the integral in (\ref{secondorder}) by 
\begin{equation}\label{TheEstOf32}
C\epsilon_j\sigma(\epsilon_j)\|\nabla \psi^2\|_{L^2},
\end{equation}
where $\sigma(\epsilon_j)\to 0$ is the modulus of continuity of Lemma \ref{FirstBadC1}.

This means that the terms in
(\ref{firstorder1})-(\ref{firstorder2}) must tend to zero as $j\to \infty$. We can thus conclude that
\begin{equation}\label{preceptions}
 \int_{B_1\setminus \Gamma_{u^j}}\Big( (2\nabla \Pi\cdot \nabla v^j)\frac{\partial \psi}{\partial x_2} -2(\nabla \Pi\cdot e_2)(\nabla v^j \cdot \nabla \psi)-
\end{equation}
\begin{equation}\label{preceptions22}
 -2(\nabla v^j \cdot e_2)(\nabla \Pi\cdot \nabla \psi)\Big)+
\end{equation}
\begin{equation}\label{decieve}
 +\int_{-1}^0 \frac{f_j'(x_1)}{\sqrt{1+\epsilon_j^2 |f_j'(x_1)|^2}}\frac{\partial \psi(x)}{\partial x_1}=
\end{equation}
\begin{equation}\label{WritingOut1}
=\int_{\Gamma_{u^j}^\pm}\Big[2(\nabla \Pi\cdot \nabla v^j)(\nu^\pm\cdot e_2)-2(\nabla \Pi\cdot e_2)(\nabla v^j\cdot \nu^\pm)-
\end{equation}
\begin{equation}\label{WritingOut2}
-2(\nabla v^j \cdot e_2)(\nabla \Pi\cdot \nu^\pm)\Big]\psi+
\end{equation}
\begin{equation}\label{sympathy}
+\int_{-1}^0 \frac{f_j'(x_1)}{\sqrt{1+\epsilon_j^2 |f'(x_1)|^2}}\frac{\partial \psi(x)}{\partial x_1}=o(1)\|\nabla \psi^j\|_{L^2}
\end{equation}

It follows from the expression (\ref{expressionNablaPi}) that on $\Gamma_{u^j}^\pm\cap (B_b\setminus B_a)$
$$
\nabla \Pi(u^j)=\frac{1}{2}\sqrt{\frac{2}{\pi}}\frac{1}{r^{1/2}}\left(\pm 1+O(\sqrt{\epsilon_j}),O(\sqrt{\epsilon_j}) \right)
.$$
Using this in the integrals (\ref{WritingOut1}) and (\ref{WritingOut2}) we may conclude from (\ref{WritingOut1})-(\ref{sympathy}) that 
\begin{equation}\label{firstequation}
o(1)=\int_{-1}^0\left( \sqrt{\frac{2}{\pi}\frac{1}{r}}\left(\frac{\partial v^j}{\partial x_1}\Big\lfloor_{\Gamma_u^+}+
\frac{\partial v^j }{\partial x_1}\Big\lfloor_{\Gamma_u^-}\right)\psi+\frac{f'_j}{\sqrt{1+\epsilon_j^2 |f_j|^2}}\frac{\partial \psi(x)}{\partial x_1}\right).
\end{equation}

By passing to the limit in (\ref{firstequation}) and using that, by Lemma \ref{FirstBadC1},
\begin{equation}\label{devil}
\sqrt{1+\epsilon_j^2|f'_j|^2}=1+o(1).
\end{equation}
we conclude that, in the weak sense,
\begin{equation}\label{SonsFindDevils1}
f_0''(x_1)=\sqrt{\frac{2}{\pi}\frac{1}{r}}\left(\frac{\partial v^0(x_1,0^+)}{\partial x_1}+\frac{\partial v^0(x_1,0^-)}{\partial x_1}\right).
\end{equation}
This proves that $f''_0$ satisfies the last equation in (\ref{limitsystprop}).

\vspace{3mm}

To derive the second and third equation for $(v^0,f_0)$ in (\ref{limitsystprop}) we use that on $\Gamma_{u^j}^\pm$
\begin{equation}\label{BoringAsHell}
0=\nabla u^j\cdot \nu= \nabla \left( \Pi(u^j)+\epsilon_j v^j\right)\cdot (\epsilon_j f_j', -1).
\end{equation}
That is
\begin{equation}\label{someBSequation}
\epsilon_j \frac{\partial v^j}{\partial x_1} f_j'-\frac{\partial v^j}{\partial x_2}=\frac{1}{\epsilon_j}\left(\nabla \Pi \cdot (-\epsilon_j f_j',1)\right).
\end{equation}
On $\Gamma_{u^j}^+$
$$
\phi=\pi-\frac{\epsilon_j f_j}{r}+O\left(\left(\frac{\epsilon_j f_j}{r}\right)^3\right)
$$
and therefore, by (\ref{expressionNablaPi}), the right side in (\ref{someBSequation}) may be written 
$$
\frac{1}{\epsilon_j}\left(\nabla \Pi \cdot (-\epsilon_j f_j',1)\right)=
$$
\begin{equation}\label{Dontever}
=\frac{1}{2\epsilon_j}\sqrt{\frac{2}{\pi}\frac{1}{r}}\Big(\epsilon_j\sin\left(\frac{\pi-\epsilon_j f_j/r}{2}\right)(f_j'+\delta_j g_j')+
\end{equation}
$$
+\cos\left(\frac{\pi-\epsilon_j f_j/r}{2}\Big) \right)=
$$
$$
=\frac{1}{2}\sqrt{\frac{2}{\pi}\frac{1}{r}}\left(f_j'(x_1)+\frac{f_j(x_1)}{2r} \right)+O(\epsilon_j).
$$

Equations (\ref{someBSequation}) and (\ref{Dontever}) together implies that, in the weak sense,
\begin{equation}\label{Nonsense}
-\epsilon_j f_j'\frac{\partial v^j(x_1, \epsilon_jf_j(x_1)^-)}{\partial x_1}-\frac{\partial v^j(x_1,\epsilon_jf_j(x_1)^-)}{\partial x_2}=
\end{equation}
$$
=\frac{1}{2}\sqrt{\frac{2}{\pi}\frac{1}{r}}\left(f_j'(x_1)+\frac{f_j(x_1)}{2r} \right)+o(\epsilon_j).
$$

Passing to the limit in (\ref{Nonsense}) we may conclude that 
\begin{equation}\label{SonsFindDevils2}
-\frac{\partial v^0(x_1,0^+)}{\partial x_2}=\frac{1}{2}\sqrt{\frac{2}{\pi}\frac{1}{r}}\left(f_0'(x_1)+\frac{f_0(x_1)}{2r} \right),
\end{equation}
which is the second equation in (\ref{limitsystprop}).

Similarly, on $\Gamma_u^-$,
where $\phi=-\pi-\frac{\epsilon_j f_j}{r}+o\left(\left(\frac{\epsilon_j f_j}{r}\right)^3\right)$,
we can conclude, after passing to the limit, that the third equation i (\ref{limitsystprop}) holds
\begin{equation}\label{SonsFindDevils2Minus}
-\frac{\partial v^0(x_1,0^-)}{\partial x_2}=-\frac{1}{2}\sqrt{\frac{2}{\pi}\frac{1}{r}}\left(f_0'(x_1)+\frac{f_0(x_1)}{2r} \right).
\end{equation}
\qed

In the next lemma we prove that the convergence needed to apply Proposition \ref{Linearization} indeed holds.

\begin{lem}\label{WhatShouldWeCallThis}
Let $(u^j,\Gamma_j)$ be a sequence of minimizers to the Mumford-Shah problem that are $\epsilon_j-$close to a crack tip for some sequence
$\epsilon_j\to 0$. Furthermore, let
\begin{equation}\label{defofVjinProp2}
 v^j(x)=\frac{u^j(x)-\Pi(u^j,1)}{\epsilon_j}.
\end{equation}
Then there exists a subsequence, which we still denote $v^j$, such that $v^j\to v^0$ strongly in $L^2$ and $\Reg(\nabla v^j)\to \Reg(\nabla v^0)$
strongly in $L^2(B_b(0)\setminus B_a(0))$ for any $0<a<b<1$. If $f_j$ is given by (see (\ref{JustifiedByThm11}))
$$
\Gamma_j\setminus \{|x_1|< \tau_j\}=\left\{ (x_1,\epsilon_j f_j(x_1));\; x_1\le -\tau_j\right\}
$$
then $f_j\to f_0$ for some function $f_0\in C([-1,0])\cap W^{1,2}((-b,-a))$ for any $0<a<b<1$.
\end{lem}
\textsl{Proof:} We will use the same idea as in the proof of Proposition \ref{Linearization} and freely referring to calculations 
made in that proof.

Let $\psi^j$ be the solution to 
\begin{equation}
 \begin{array}{ll}
  \Delta \psi^j=0 & \textrm{ in }(B_b\setminus B_a)\setminus \Gamma_{u^j} \\
  \psi^j= 0 & \textrm{ on } \partial B_b \cup \partial B_a \\
  \psi^j = f_j -\frac{f_j(b)-f_j(a)}{b-a}(-x_1-a)+f^j(a) & \textrm{ on }\Gamma_{u^j}.
 \end{array}
\end{equation}
Then, since $\Gamma_{u^j}$ is a $C^{1,\alpha}-$graph with small norm in $B_b\setminus B_a$, it follows that 
\begin{equation}\label{Sylt}
 \|\nabla \psi^j\|_{L^2(B_b\setminus B_a)}\le C\|f_j'-(f_j')_{(-a,-b)}\|_{L^2(-a,-b)},
\end{equation}
where we use the notation $(f'_j)_{(-a,-b)}$ for the average of $f'_j$ on the interval $(-a,-b)$. 
Inserting $\psi^j$ for $\psi$ in (\ref{preceptions})-(\ref{sympathy}) we see that
\begin{equation}\label{preceptions2}
 \int_{B_1\setminus \Gamma_{u^j}}\left( (2\nabla \Pi\cdot \nabla v^j)\frac{\partial \psi}{\partial x_2}-2(\nabla \Pi\cdot e_2)(\nabla v^j\cdot \nabla \psi^j)-
 2(\nabla v^j\cdot e_2)(\nabla \Pi\cdot \nabla \psi^j)\right)+
\end{equation}
\begin{equation}\label{decieve2}
 +\int_0^1 \frac{f_j'(x_1)}{\sqrt{1+\epsilon_j^2 |f_j'(x_1)|^2}}\frac{\partial \psi^j(x)}{\partial x_1}=o(1)\|\nabla \psi^j\|_{L^2(B_1\setminus \Gamma_{u^j})}.
\end{equation}
Rearranging terms in (\ref{preceptions2})-(\ref{decieve2}) and using that $\|\Reg(\nabla v^j)\|_{L^2}=1$ by normalization and that 
$\|\nabla \Pi(u^j)\|_{L^2}\le C$ (which follows from Lemma \ref{VariationsInTang} if $\epsilon_j$ is small) we may conclude that 
$$
\int_{-a}^{-b}\frac{|f_j'|^2}{\sqrt{1+\epsilon_j^2|f'_j|^2}} \le C\|f_j'-(f_j')_{(-a,-b)}\|_{L^2(-a,-b)}+
$$
\begin{equation}\label{EstOfF}
+\int_{-a}^{-b}\frac{f'_j (f'_j)_{(-a,-b)}}{\sqrt{1+\epsilon_j^2|f'_j|^2}}+C\|\nabla \psi^j\|_{L^2} +C\le
\end{equation}
$$
\le C\Big[\big\|f_j'-(f_j')_{(-a,-b)}\big\|_{L^2(-a,-b)}+|(f'_j)_{(-a,-b)}|^2+C\|f_j'\|_{L^2}\Big]+C,
$$
where we used that we can estimate $\nabla \psi$ in terms of its boundary data $f_j$ in the last equality.
In order to show that $\|f'_j\|_{L^2}$ is bounded we need to control the average $|(f'_j)_{(-a,-b)}|^2$.

To estimate the average $|(f'_j)_{(-a,-b)}|^2$ we will use (\ref{Nonsense}) which may be formulated:
for any $\zeta\in C^{\infty}_{0}(B_b\setminus B_a)$,
\begin{equation}\label{sqrtx1fprimest2}
\int_{\Gamma_{u^j}^+} \left( \epsilon_j f_j'\frac{\partial v^j}{\partial x_1}-\frac{\partial v^j(x_1,\epsilon_jf_j(x_1)^+)}{\partial x_2}\right)\zeta d\H\lfloor_{\Gamma_{u^j}}^1+
\end{equation}
\begin{equation}\label{Ivar}
+\int_{\Gamma_{u^j}^-} \left( -\epsilon_j f_j'\frac{\partial v^j}{\partial x_1}+\frac{\partial v^j(x_1,\epsilon_jf_j(x_1)^-)}{\partial x_2}\right)\zeta d\H\lfloor_{\Gamma_{u^j}}^1+
\end{equation}
\begin{equation}\label{NeedMoreSeriousNames2}
+\frac{1}{2}\sqrt{\frac{2}{\pi}}\int_{\Gamma_{u^j}}\frac{\partial}{\partial x_1}\left(\frac{1}{\sqrt{-x_1}}f(x_1)\right)\zeta d\H\lfloor_{\Gamma_{u^j}}^1=o(\epsilon_j).
\end{equation}
Notice that, by an integration by parts in (\ref{sqrtx1fprimest2})-(\ref{Ivar}), using that of $\Gamma_{u^j}$ is the graph of $\epsilon f'_j$ c.f. (\ref{BoringAsHell}),
\begin{equation}\label{thezetaintbyparts}
\int_{\Gamma_{u^j}^+}  \left( \epsilon_j f_j'\frac{\partial v^j}{\partial x_1}-\frac{\partial v^j(x_1,\epsilon_j f_j(x_1)^+)}{\partial x_2}\right)\zeta+
\end{equation}
\begin{equation}
\int_{\Gamma_{u^j}^-}  \left( \epsilon_j f_j'\frac{\partial v^j}{\partial x_1}-\frac{\partial v^j(x_1,\epsilon_j f_j(x_1)^-)}{\partial x_2}\right)\zeta=
\end{equation}
\begin{equation}\label{Ivar2}
=\int_{B_1\setminus \Gamma_u}\nabla v^j\cdot \nabla \zeta\le \|\nabla v^j\|_{L^2}\|\nabla \zeta\|_{L^2}\le \|\nabla \zeta\|_{L^2}.
\end{equation}
we can conclude from (\ref{sqrtx1fprimest2})-(\ref{NeedMoreSeriousNames2}) and (\ref{thezetaintbyparts})-(\ref{Ivar2}) that
\begin{equation}\label{ToBEMax}
\frac{1}{2}\sqrt{\frac{2}{\pi}}\int_{\Gamma_{u^j}}\frac{\partial}{\partial x_1}\left(\frac{1}{\sqrt{-x_1}}f_j(x_1)\right)\zeta \le \left(\int_{B_1\setminus \Gamma}|\nabla \zeta|^2\right)^{1/2}+o(1)
\end{equation}
Choosing the $\zeta$ that maximizes the left side in (\ref{ToBEMax}) under the constraint $\|\nabla \zeta\|_{L^2}\le 1$ we can conclude that
$$
\left\|\frac{\partial}{\partial x_1}\left(\frac{1}{\sqrt{-x_1}}f_j(x_1)\right)\right\|_{H^{-1/2}}\le 2.
$$
It follows that $\frac{1}{\sqrt{-x_1}}f_j(x_1)\in H^{1/2}$ modulo solutions, $h(x_1)$, to the ODE
$$
\frac{\partial}{\partial x_1}\left(\frac{1}{\sqrt{-x_1}}h(x_1)\right)=0.
$$
We can conclude that $f_j\in L^{2}(a,b)$ and that for some $\gamma_j\in \R$ and some constant $C$
\begin{equation}\label{AllYouWEverWanted}
\left\|f_j-\gamma_j \sqrt{-x_1}\right\|_{L^2(-a,-b)}\le C.
\end{equation}
But (\ref{AllYouWEverWanted}) clearly implies that $(f_j')_{(-a,-b)}$ is uniformly bounded. If follows, from this and (\ref{EstOfF}), 
that 
\begin{equation}\label{Shitstar2}
\|f_j\|_{W^{1,2}(-a,-b)}\le C.
\end{equation}

\vspace{2mm}

It still remains to show to show that $\Reg(\nabla v^j)$ converges strongly in $L^2(B_b\setminus B_a)$.
We notice that since $u^j$ has zero Neumann data on $\Gamma_u^\pm$ it follows that
\begin{equation}\label{Shitstar1}
\nabla v^j\cdot \nu_j^\pm= -\frac{1}{\epsilon_j} \nabla \Pi(u^j)\cdot \nu_j^\pm \textrm{ on }\Gamma_u^\pm \cap B_{(1+b)/2}(0)\setminus B_{a/2}(0).
\end{equation}
A direct calculation shows that
\begin{equation}\label{Shitstar3}
\nabla \Pi(u^j)\cdot \nu_j^\pm=
\frac{1}{2}\sqrt{\frac{2}{\pi}\frac{1}{r}}
\left(\epsilon_j\sin\left(\frac{\pi-\epsilon_jf_j/r}{2}\right)f_j'+\cos\left(\frac{\pi-\epsilon_j f_j/r}{2\sqrt{1+\epsilon_j^2|f'_j|^2}}\right) \right).
\end{equation}
It follows, from $f_j\in W^{1,2}((-(1+b)/2,-a/2))$ that
$$
\nabla \Pi(u^j)\cdot \nu_j^\pm\in L^2(\Gamma^\pm_u\cap B_{(1+b)/2}(0)\setminus B_{a/2}(0)).
$$
In particular, the non-tangential maximal function of $\nabla v^j$ is an $L^2-$function. From this it easily follows that
$$
\left\| \nabla v^j\right\|_{L^2(\{|x_2|\le \delta\})\cap B_{b}(0)\setminus B_a(0))}\le \sigma(\delta),
$$
if $j$ is large enough and $\sigma(\delta) \to 0$ as $\delta\to 0$.

By the triangle inequality we can conclude that
$$
\lim_{j\to \infty}\left\| \nabla \left(\Reg(v^j)-\Reg(v^0)\right)\right\|_{L^2(B_b\setminus B_a)}\le
$$
$$
\le \lim_{j\to \infty}\Big( \left\| \nabla \left(\Reg(v^j)-\Reg(v^0)\right)\right\|_{L^2(B_b\setminus (B_a\cup \{|x_2|\le \delta\}))}+
$$
$$
+\left\| \nabla \left(\Reg(v^j)-\Reg(v^0)\right)\right\|_{L^2(B_b\cap \{|x_2|\le \delta\}\setminus B_a)}\Big)\le
\sigma(\delta)
$$
since $v^j$ is converges uniformly in $C^1$ in the compact set $B_b\setminus (B_a\cup \{|x_2|\le \delta\})$. Since $\delta>0$ is arbitrary
we can conclude that
$$
\lim_{j\to \infty}\left\| \nabla \left(\Reg(v^j)-\Reg(v^0)\right)\right\|_{L^2(B_b\setminus B_a)}=0.
$$
This proves the claim.  \qed

\begin{cor}\label{FconvC1alpha}
 Under the assumptions of Lemma \ref{WhatShouldWeCallThis} and for any $-1<-b<-a<0$
 $$
 f_j\to f_0 \quad \textrm{ in }C^{1,\alpha}(-b,-a),
 $$
 for every $\alpha <1/2$. In particular, for $\epsilon_j$ small enough $\|f_j\|_{C^{1,\alpha}(-b,-a)}\le C$
 where $C$ is independent of $\epsilon$ (but may depend on $a$ and $b$).
\end{cor}
\textsl{Proof:} Since $f_j\in W^{1,2}$ it follows that $\left[|\nabla u^j(x_1,f_j(x_1))|^2\right]^\pm\in W^{1,2}$ and 
therefore the curvature of $\Gamma_j$ is in $L^2$ which implies that $f_j\in W^{2,2}\subset C^{1,\alpha}$ for every $\alpha<1/2$. 
The corollary follows by compactness in $C^{1,\alpha}$. \qed


\section{Strong Convergence.}\label{SceStrongConv}

In this section we prove that the linearizing sequence $\Reg(\nabla v^j)$ converges locally strongly in $L^2(B_1)$ and $f_j$ converges locally 
strongly in $W^{1,2}((-1,0))$. Throughout this section 
$u^j$, $\Gamma_{u^j}$, $v^j$ and $f_j$ will be as in Proposition \ref{Linearization}.

We begin this section with a lemma that proves strong convergence under an extra assumption.

\begin{lem}\label{EstIfThereExistsBound}
 Let $u^j$, $\Gamma_{u^j}$, $v^j$ and $f_j$ be as in Proposition \ref{Linearization}. Furthermore assume that there exists a constant $C$ such that,
 for some $0<\kappa<1/2$,
 \begin{equation}\label{ToBeUsedWithClarity}
 r^\kappa \left\|\nabla\left( \frac{u^j(rx)}{\sqrt{r}}-\Pi\left(\frac{u^j(rx)}{\sqrt{r}}\right)\right)\right\|_{L^2(B_1(0)\setminus \Gamma_{u^j_r})}\le
 C\epsilon_j
 \end{equation}
 for all $j$ and $r\in (0,1)$.

 Then
 \begin{enumerate}
  \item $\Reg(\nabla v^j)\to \Reg(\nabla v^0)$ strongly in $L^2(B_1)$.
  \item and $v^0$ and $f_0$ satisfies the following estimates
  \begin{equation}\label{v0controlandshit}
  \|\nabla v^0\|_{B_r(0)\setminus \Gamma_0}\le C(1+\ln(1/r))\sqrt{r},
  \end{equation}
  \begin{equation}\label{f0pointwise}
  |f_0(x_1)|\le C(1+\ln(1/|x_1|))|x_1|
  \end{equation}
  and
  \begin{equation}\label{fprimepointwise}
  |f_0'(x_1)|\le C(1+\ln(1/|x_1|)).
  \end{equation}
 \end{enumerate}

\end{lem}
\textsl{Proof:} Using (\ref{ToBeUsedWithClarity}) we see that 
$$
\left\|\nabla\left(\Pi(u^j,2^{-k})-\Pi(u^j,2^{-k-1}) \right)\right\|_{L^2(B_{1/2}\setminus \Gamma_{u^j(2^{-k}x)})}\le
$$
$$
\le \left\|\nabla\left(\frac{u^j(2^{-k}x)}{2^{-k/2}}-\Pi(u^j,2^{-k}) \right)\right\|_{L^2(B_{1/2}\setminus \Gamma_{u^j(2^{-k}x)})}+
$$
$$
+\left\|\nabla\left(\frac{u^j(2^{-k-1}x)}{2^{-(k+1)/2}}-\Pi(u^j,2^{-k-1}) \right)\right\|_{L^2(B_{1/2}\setminus \Gamma_{u^j(2^{-k}x)})}\le
$$
$$
\le C\epsilon_j2^{\kappa(k+1)},
$$
where $C$ only depend on the constant $C$ in (\ref{ToBeUsedWithClarity}).

It follows that
\begin{equation}\label{MutherF}
\left\| \nabla \left(\Pi(u^j, 1)-\Pi(u^j,2^{-k})\right)\right\|_{L^2(B_1)}\le
\end{equation}
$$
\le \sum_{l=0}^{k-1}\left\| \nabla \left(\Pi(u^j, 2^{-l})-\Pi(u^j,2^{-l+1})\right)\right\|_{L^2(B_1)}\le C2^{\kappa k}\epsilon_j.
$$

In particular, using (\ref{ToBeUsedWithClarity}) again,
$$
\|\Reg(\nabla v^j)\|_{L^2(B_{2^{-k}})}\le 
$$
$$
\le \frac{2^{-k/2}}{\epsilon_j}\left\|\nabla\left( \frac{u^j(2^{-k}x)}{2^{-k/2}}-\Pi\left(\frac{u^j(2^{-k}x)}{2^{-k/2}}\right)\right)\right\|_{L^2(B_{1/2}(0)\setminus \Gamma_{u^j_r})}
$$
\begin{equation}\label{ClearlyOrWhat}
+\frac{2^{-k/2}}{\epsilon_j}\left\|\nabla\left(\Pi(u^j,2^{-k})-\Pi(u^j) \right)\right\|_{L^2(B_{1/2}\setminus \Gamma_{u^j})}\le C2^{-(1/2-\kappa) k}+
\end{equation}
$$
+\sum_{i=1}^{k}\frac{1}{\epsilon_j}\left\|\nabla\left(\Pi(u^j,2^{-i})-\Pi(u^j, 2^{-(i-1)}) \right)\right\|_{L^2(B_{1/2}\setminus \Gamma_{u^j})}\le
$$
$$
\le C2^{-k(1/2-\kappa)},
$$
where we used (\ref{MutherF}) to estimate the series and scaling invariance in the integrals. We also notice that (\ref{ClearlyOrWhat})
implies that 
\begin{equation}\label{Almost61}
 \|\Reg( \nabla v^j)\|_{L^2(B_r(0)}\le Cr^{-\kappa}\sqrt{r}.
\end{equation}
The same estimate carries over to $v^0$, this is only slightly weaker than (\ref{v0controlandshit}) (the full strength of (\ref{v0controlandshit}) will 
be proved shortly).

The estimate (\ref{ClearlyOrWhat}) allows us, for any $\epsilon>0$, to find a $k$ such that 
$$
\|\Reg(\nabla v^j)\|_{L^2(B_{2^{-k}})}<\epsilon,
$$
and since $\Reg(\nabla v^j)$ converges strongly in $B_1(0)\setminus B_{2^{-k}}$ (by Lemma \ref{WhatShouldWeCallThis}) it follows that 
$\Reg(\nabla v^j)$ converges strongly in $B_1(0)$.

\vspace{3mm}

To derive the desired estimates for $f_0$ we notice that if
\begin{equation}\label{ThisWouldBeLateStar}
\Pi(u,r)=\lambda_r r^{1/2}\sin\left(\frac{\phi-\phi_r}{2}\right)
\end{equation}
then, using (\ref{MutherF}),
$$
\left\|\nabla \left(\Pi(u^j,1)-\Pi(u^j,2^{-k})\right)\right\|_{L^2(B_1\setminus \Gamma_{u^j})}\le 
$$
\begin{equation}\label{ThisWouldBeLateStar2}
\le \sum_{i=1}^{k}\frac{1}{\epsilon_j}\left\|\nabla\left(\Pi(u^j,2^{-i})-\Pi(u^j, 2^{-(i-1)}) \right)\right\|_{L^2(B_{1/2}\setminus \Gamma_{u^j})} \le
\end{equation}
$$
\le C2^{-\kappa} 2^{-k}\epsilon_j.
$$
This implies, using a Taylor expansion in (\ref{ThisWouldBeLateStar}) and (\ref{ThisWouldBeLateStar2}), that if $2^{\kappa k}\epsilon_j$ is small then
\begin{equation}\label{TwistedSister}
 \left|\phi_1-\phi_{2^{-k}}\right|\le C2^{\kappa k} \epsilon_j.
\end{equation}

We will now prove slightly weaker versions of (\ref{f0pointwise}) and (\ref{fprimepointwise}).
To prove the weaker version (\ref{f0pointwise}) we just notice that rotating the coordinate system by an angle $\phi_k$
amounts to subtracting a linear function $l(x_1)=a_k x_1$ from $f_j$, modulo lower order terms, where $ \epsilon_j a_k\approx \phi_k+
O\left(\left( \phi_k\right)^3\right)$ for $\phi_k$ small enough (that which follows from $\epsilon_j$ being small). In particular $|a_k|\le C$.

From Corollary \ref{L2boundsCor} we may conclude that
\begin{equation}\label{VemSnarkarPaMinDorr}
\sup_{B_1(0)\setminus B_{1/2}(0),\; x_1<0}\left| \frac{\epsilon_j f_j(2^{-k} x_1)}{2^{-k}}-a_k x_1\right|\le C 2^{\kappa k} \epsilon_j.
\end{equation}
But (\ref{VemSnarkarPaMinDorr}), together with $|a_k|\le C $, implies that
$$
-C|x_1|^{1-\kappa}\le f(x_1)\le C|x_1|^{1-\kappa}.
$$

To prove the weaker version of (\ref{fprimepointwise}) one argues similarly. In particular, Lemma \ref{FirstBadC1}, together with (\ref{VemSnarkarPaMinDorr})
and that $u^j$ is $C\epsilon_j r^{-\kappa}$ close to a crack-tip (by \ref{v0controlandshit}) will imply that
$|f'(x_1)-a_k|\le C |x_1|^{-\kappa}$ for any $x_1\in (-2^{-k},-2^{-k-1})$. We leave the details to the reader. 

We have now shown that $(v^0,f_0)$ satisfies the estimates (\ref{v0controlandshit}), (\ref{f0pointwise}), (\ref{fprimepointwise})
with $r^{-\kappa}$ in place of the $\ln(1/r)$ term. But, by Proposition \ref{Linearization} $(v^0,f_0)$ also satisfies the linearized
system. Therefore Corollary \ref{Sigma} implies the estimates (\ref{v0controlandshit}), (\ref{f0pointwise}) and (\ref{fprimepointwise}).
\qed

\begin{lem}\label{CloseForSmallerr}
 Let $u^j,$ $\Gamma_{u^j}$ be as in Proposition \ref{Linearization}. Then
 \begin{equation}\label{limrtozeroepsilon}
 \lim_{j\to \infty} \sup_{r\in (0,1]}\left\|\nabla \left( \frac{u^j(rx)}{\sqrt{r}}-\Pi(u^j,r)\right)\right\|_{L^2(B_1\setminus \Gamma_{u^j(rx)})}=0
 \end{equation}
 and for each $j$ the supremum is achieved at some $r_j\in (0,1]$.
\end{lem}
\textsl{Proof:} First we notice that, by Corollary \ref{BlowUpAtCarackTipCor}, since for every $j$ and any subsequence $r\to 0$ 
$$
\lim_{r\to 0}\Reg\left( \nabla \frac{u^j(rx)}{\sqrt{r}}\right)=\Reg\left( \nabla \sqrt{\frac{2}{\pi}}r^{1/2}\sin((\phi+\phi_0)/2)\right),
$$
as long as the limit exists and the angle $\phi_0$ could depend on the sub-sequence. That is, for every $j$,
$$
\lim_{r\to 0}\left\|\nabla \left( \frac{u^j(rx)}{\sqrt{r}}-\Pi(u^j,r)\right)\right\|_{L^2(B_1\setminus \Gamma_{u^j(rx)})}=0.
$$
It follows that any positive supremum of
$$
\left\|\nabla \left( \frac{u^j(rx)}{\sqrt{r}}-\Pi(u^j,r)\right)\right\|_{L^2(B_r\setminus \Gamma_{u^j})}
$$
occurs at a strictly positive $r$.

To see that (\ref{limrtozeroepsilon}) holds we argue by contradiction and assume that there exists a $\delta>0$ and $r_j>0$ such that
the supremum is achieved at $r_j$ and
$$
\frac{1}{\sqrt{r_j}}\left\|\nabla \left( u^j-\Pi(u^j,r)\right)\right\|_{L^2(B_{r_j}\setminus \Gamma_{u^j})}\ge \delta>0.
$$

We also notice that, by Lemma \ref{CPiLem} and the remark thereafter,
\begin{equation}\label{EngFuTo1}
\Pi(u^j,1)\to \sqrt{\frac{2}{\pi}}r^{1/2}\sin(\phi/2)\Rightarrow \int_{B_1\setminus \Gamma_{u^j}}|\nabla u^j|^2\to 1
\end{equation}
and, by Lemma \ref{BonnetMonFormula} and the remark after Lemma \ref{CPiLem},
\begin{equation}
\lim_{r\to 0}\int_{B_1\setminus \Gamma_{u^j}}\left|\nabla \frac{u^j(rx)}{\sqrt{r}}\right|^2\to 1.
\end{equation}

We may conclude that $\frac{u^j(r_j x)}{\sqrt{r_j}}\to u^0$ where $u^0$ is a minimizer of the Mumford-Shah functional,
\begin{equation}\label{WhereTheKnife}
\left\|\nabla \left( u^0-\Pi(u^0,1)\right)\right\|_{L^2(B_{1}\setminus \Gamma_{u^0})}\ge \delta>0
\end{equation}
and
\begin{equation}\label{AjAjAj}
1=\lim_{r\to 0}\frac{1}{r}\int_{B_r(0)\setminus \Gamma_{u^0}}|\nabla u^0|^2\le \frac{1}{r}\int_{B_r(0)\setminus \Gamma_{u^0}}|\nabla u^0|^2\le
\int_{B_1(0)\setminus \Gamma_{u^0}}|\nabla u^0|^2=1
\end{equation}
where we used (\ref{EngFuTo1}) and the Monotonicity formula (Lemma \ref{BonnetMonFormula}) in the last step. 

The monotonicity formula (Lemma \ref{BonnetMonFormula}) together with (\ref{AjAjAj}) implies that $u^0$ is homogeneous which implies that
$$
u^0=\sqrt{\frac{2}{\pi}}r^{1/2}\sin(\phi/2)
$$
which contradicts (\ref{WhereTheKnife}).

We can conclude that
$$
\lim_{j\to \infty} \sup_{r\in (0,1]}\left\|\nabla \left( \frac{u^j(rx)}{\sqrt{r}}-\Pi(u^j,r)\right)\right\|_{L^2(B_r\setminus \Gamma_{u^j})}=0.
$$
\qed

\begin{prop}\label{SatisfiesRightEstProp}
 Let $u^j$ and $\Gamma_{u^j}$ be as in Proposition \ref{Linearization}.

 For every $0<\kappa<1/2$ there exist a constant $C$ such that, for every $r\in (0,1]$,
 \begin{equation}\label{CroMangon}
 \left\|\nabla\left( \frac{u^j(rx)}{\sqrt{r}}-\Pi\left(\frac{u^j(rx)}{\sqrt{r}}\right)\right)\right\|_{L^2(B_1(0)\setminus \Gamma_{u^j_r})}\le
 Cr^{-\kappa}\epsilon_j.
 \end{equation}

In particular, $u^0$ and $f_0$ satisfy the estimates in Lemma \ref{EstIfThereExistsBound}.
\end{prop}
\textsl{Proof:} We will again argue by contradiction and assume that there is a sequence $(u^j,f_j)$ that are $\epsilon_j\to 0$
close to a crack-tip such that
\begin{equation}\label{theonewithTilde}
\frac{1}{\epsilon_j\tilde{r}_j^{-\kappa}}\left\|\nabla\left( \frac{u^j(\tilde{r}_jx)}{\sqrt{\tilde{r}_j}}-\Pi\left(\frac{u^j(\tilde{r}_jx)}{\sqrt{\tilde{r}_j}}\right)\right)\right\|_{L^2(B_1(0)\setminus \Gamma_{u^j_{\tilde{r}_j}})}= j.
\end{equation}
for some sequence $\tilde{r}_j\in (0,1]$ such that the expression on the right in (\ref{theonewithTilde}) is maximized for $\tilde{r}_j$.
By Lemma \ref{CloseForSmallerr} it follows that $j\epsilon_j\tilde{r}_j^{-\kappa}\to 0$ as $j\to \infty$. It is also easy to see that $\tilde{r}_j\to 0$.

We let $r_j=c\tilde{r}_j$ where $c>1$ is some fixed constant, to be defined later, that depends on the norm of the mapping $A-I$ from (\ref{MapL2toL2}).

We define
$$
u^j_{r_j}(x)=\frac{u^j(r_j x)}{\sqrt{r_j}}.
$$
Then $u^j_{r_j}$ satisfies the criteria in Lemma \ref{EstIfThereExistsBound} and the sequence
$$
v^j(x)=\frac{\frac{u^j(r_j x)}{\sqrt{r_j}}-\Pi(\frac{u^j(r_j x)}{\sqrt{r_j}},1)}{j\epsilon_j r_j^{-\kappa}}
$$
converges strongly to a solution $v_0$ of Proposition \ref{AnalysisOfLinear}.

In particular, $v^0$ satisfies the series expansions in (\ref{SumOfHomos})
\begin{equation}\label{morehomos}
v^0(r,\phi)=a+a_0 \mathfrak{z}(r,\phi)+\sum_{k=1}^\infty a_k r^{\alpha_k}\cos(\alpha_k \phi)+ \sum_{k=1}^\infty b_k r^{k-1/2}\sin((k-1/2)\phi)
\end{equation}
 and (\ref{SumofHetros})
$$
f_0(x_1)=a_0 \mathfrak{h}(|x_1|)+\sum_{k=1}^\infty 2a_k \sqrt{\frac{\pi}{2}}\sin(\alpha_k \pi)|x_1|^{\alpha_k+\frac{1}{2}}.
$$

But by our choice of $r_j=c\tilde{r}_j$, and that the maximum occurs at $\tilde{r}_j$, it follows that
$$
\left\|\nabla\left( \frac{u^j(r_jx)}{\sqrt{r_j}}-\Pi\left(\frac{u^j(r_jx)}{\sqrt{r}}\right)\right)\right\|_{L^2(B_1(0)\setminus \Gamma_{u^j_{r_j}})}\le
$$
$$
\le c^{-\kappa} \left\|\nabla\left( \frac{u^j(\tilde{r}_jx)}{\sqrt{\tilde{r}_j}}-\Pi\left(\frac{u^j(\tilde{r}_jx)}{\sqrt{\tilde{r}_j}}\right)\right)\right\|_{L^2(B_1(0)\setminus \Gamma_{u^j_{\tilde{r}_j}})}
$$
which implies that
$$
\|\nabla v^0(x)\|_{L^2(B_1\setminus \Gamma_0)}\le c^{-\kappa}\left\|\nabla (\sqrt{c}v^0(x/c))\right\|_{L^2(B_1\setminus \Gamma_0)}
$$
which is not true if we choose $c>1$ large enough (since all the terms in the series expansion of $v^0$ (\ref{morehomos}) have homogeneity greater 
than $1/2$). This is a contradiction. We may conclude that there exists a constant $C$ such that (\ref{CroMangon}) holds.

By the inequality (\ref{CroMangon}) and Lemma \ref{EstIfThereExistsBound} the second conclusion holds. \qed


\section{Getting Rid of the First Term in the Asymptotic Expansion.}\label{C1alphaSEC}

In this section we prove the first regularity improvement for the linearized system. We will prove a simple lemma that states that we may rotate
the coordinate system to get rid of the $\cos(\phi/2)$ term in the asymptotic expansion of $v^0$. The Lemma is a somewhat annoying technical curiosity. 
The need for this lemma arises because we define $v^j$ in such a way that $v^j$ is orthogonal to $\sin(\phi/2)$ and $\cos(\phi/2)$ - which is natural 
when we consider convergence properties of $v^j$. But for the regularity theory we would want $v^0$ to consist of terms of higher homogeneities than $1/2$. 
But the even terms, $r^{\alpha_k}\cos(\alpha_k\phi)$, in the homogeneous expansion of $v^0$ are not orthogonal to $r^{1/2}\cos(\phi/2)$. These 
non-orthogonality properties means that we get an extra $r^{1/2}\cos(\phi/2)$ term in the homogeneous expansion of the solutions to the 
linearized problem in Proposition \ref{AnalysisOfLinear}. Fortunately it is rather easy to get rid of this extra term by slightly rotating the 
coordinate system.

\begin{lem}\label{ImprovedLinearization}
Let $(u^j,\Gamma_j)$ be as in Proposition \ref{Linearization} (In particular Proposition \ref{SatisfiesRightEstProp} holds). Then there exists a sequence
of rotations of $\R^2$, $P^j:\R^2\mapsto \R^2$ such that if we express $(u^j,\Gamma_u^j)$ in these rotated coordinate systems and define $v^j\to v^0$ and
$f_j\to f_0$ in these rotated coordinate systems then
\begin{equation}\label{formofTildev0}
v^0(r,\phi)=a+a_0 \mathfrak{z}(x)+\sum_{k=2}^\infty a_k r^{\alpha_k}\cos(\alpha_k \phi)+ \sum_{k=2}^\infty b_k r^{k-1/2}\sin((k-1/2)\phi)
\end{equation}
and
\begin{equation}\label{formofTildef0}
f_0(x_1)=a_0 \mathfrak{h}(|x_1|)+\sum_{k=2}^\infty 2a_k \sqrt{\frac{\pi}{2}}\sin((k-1/2) \pi)|x_1|^{k-\frac{1}{2}},
\end{equation}
where $\alpha_k$, $\mathfrak{z}$ and $\mathfrak{h}$ are as in Proposition \ref{AnalysisOfLinear}.
\end{lem}
\textsl{Proof:} We choose $P^j$ to be the rotation of an angle $c\epsilon_j$, $P^j(r,\phi)=(r, \phi +c\epsilon_j)$.
We may apply Proposition \ref{Linearization} on the sequence of rotated solutions \linebreak $(u^j(P^j(r,\phi)), P^j(\Gamma_{u^j}))$.
By Proposition \ref{Linearization} and Proposition \ref{SatisfiesRightEstProp} it follows that the convergence $v^j\to v^0$ and $f_j\to f_0$
is strong and $(v^0,f_0)$ will satisfy the estimates needed to apply the second half of Proposition \ref{AnalysisOfLinear}. Therefore by
Proposition \ref{AnalysisOfLinear} $v^0$ and $f_0$ will be of the following form
$$
 v^0(r,\phi)=a+a_0 \mathfrak{z}(x)+\sum_{k=1}^\infty a_k r^{\alpha_k}\cos(\alpha_k \phi)+ \sum_{k=1}^\infty b_k r^{k-1/2}\sin((k-1/2)\phi)
$$
 and
$$
 f_0(x_1)=a_0 \mathfrak{h}(|x_1|)+\sum_{k=1}^\infty a_k 2\sqrt{\frac{\pi}{2}}\sin(\alpha_k \pi)|x_1|^{\alpha_k+\frac{1}{2}}
$$

We need to show that, if we choose the rotations appropriately, then $a_1=0$. It is easy to see that slightly rotate the
coordinate system by $c\epsilon_j$ amounts to adding $c x_1$ to the limit function $f_0$. Thus by choosing
$c$ appropriately we will get $a_1=0$. The formal argument is not very illustrative. But we provide the details for the sake of completeness.

Remember that $u^j(r,\phi)=\Pi(u^j)+\epsilon_j v^j(r,\phi)$ and that, for some $\lambda_j\in \R$,
$$
\Pi(u^j)(r,\phi)=\sqrt{\frac{2}{\pi}}\lambda_j \sin\left( \frac{\phi}{2}\right).
$$
This implies that
$$
u^j(r,\phi+c\epsilon_j)=\Pi(u^j)(r,\phi +c\epsilon_j)+\epsilon_j v^j(r,\phi+c\epsilon_j)=
$$
$$
=\sqrt{\frac{2}{\pi}}\lambda_j \left(\sin\left( \frac{\phi}{2}\right)\cos\left(\frac{c\epsilon_j}{2}\right)+
\cos\left( \frac{\phi}{2}\right)\sin\left(\frac{c\epsilon_j}{2}\right)\right)+\epsilon_j v^j(r,\phi+c\epsilon_j)=
$$
\begin{equation}\label{AddedExpl}
=\sqrt{\frac{2}{\pi}}\left[\lambda_j\cos\left(\frac{c\epsilon_j}{2}\right)\right]\sin\left( \frac{\phi}{2}\right)+
\epsilon_j \left[\frac{c}{2}\,\frac{\sin\left(\frac{c\epsilon_j}{2}\right)}{\frac{c\epsilon_j}{2}}\cos\left( \frac{\phi}{2}\right)+v^j(r,\phi+c\epsilon_j)\right].
\end{equation}
If we denote the first square bracket in (\ref{AddedExpl}) by $\tilde{\lambda}_j$ and the second square bracket $\tilde{v}^j$ then we see that
if $c$ is chosen so that
$$
\lim_{j\to \infty}=\frac{c}{2}\,\frac{\sin\left(\frac{c\epsilon_j}{2}\right)}{\frac{c\epsilon_j}{2}}\cos\left( \frac{\phi}{2}\right)=-a_1
$$
then it follows that $\tilde{v}^j\to \tilde{v}^0$ where $\tilde{v}^0$ has the form in (\ref{formofTildev0}). Since $(\tilde{v}^0, \tilde{f}_0)$
will satisfy the equations (\ref{theLinearSyst}) it follows that $\tilde{f}^0$ has the form of (\ref{formofTildef0}), see Proposition \ref{AnalysisOfLinear}.
\qed


\section{Variations in the Orthogonal Direction.}\label{OrthVarSec}

In this section we show that the $\mathfrak{z}$ and $\mathfrak{h}$ terms does not appear in the expansions of the limit of linearizing sequences 
of minimizers to the Mumford-Shah problem. This will directly imply the $C^{1,\alpha}$ regularity at the crack-tip.

The proofs consists mostly of rather tedious Taylor expansions that we are able to preform since we have exact information of the asymptotic expansions
of solutions.

\begin{lem}\label{LemVarInOrth}
 Assume that $(u,\Gamma_u)$ is a minimizer of the Mumford-Shah problem in $B_1(0)$ (with given boundary data) under the restriction that 
 $\Gamma_u$ is a connected rectifiable set that starts at the origin and ends at $\partial B_1(0)$.
 
 Furthermore assume that 
 \begin{equation}\label{SoDamedBoring}
 u=\sqrt{\frac{2}{\pi}}r^{1/2}\sin(\phi/2)+\epsilon\left(r^{1/2}\phi \sin(\phi/2)-r^{1/2}\ln(r)\cos(\phi/2) \right)+R_0
 \end{equation}
 and (for some small $\tau$)
 $$
 \Gamma_u\cap \left(B_1(0)\setminus B_\tau(0)\right)=\{(x_1,\epsilon f(x_1));\; x_1\in (-1,-\tau)\}
 $$
 where $f(x_1)\approx \sqrt{2\pi}x_1\ln(-x_1)$ and $\|\Reg(\nabla R_0)\|_{L^2(B_1(0))}\le \sigma \epsilon$. Then if 
 $\epsilon$ and $\sigma$ are small enough the pair $(u,\Gamma_u)$ is not a minimizer to the Mumford-Shah problem without the restriction that 
 $\Gamma_u$ starts an the origin and ends at $\partial B_1(0)$.
\end{lem}
\textsl{Proof:} The proof is entirely trivial, although it depends on significant calculations and in particular on rather messy Taylor expansions.
We will throughout the proof assume that $\sigma >0$ is a small fixed constant.

There is no loss of generality, possibly after rescaling to $\sqrt{2}u(x/2)$, to assume that $u$ is defined in $B_2(0)$ and that 
$|\nabla (u-\Pi(u))|<C\epsilon$ and $|\nabla R_0|\le C\sigma \epsilon$ in a neighborhood of $\partial B_1(0)\setminus \Gamma_u$.

The idea of the proof is to construct a competitor for minimality $(w,\Gamma_w)$ with less energy than $(u,\Gamma_u)$
in the ball $\hat{B}=B_1(\delta e_2)$ for some small $\delta$ satisfying  $\sigma\epsilon < < \delta < < \epsilon$. 

The basic construction of $(w,\Gamma_w)$ is presented in the figure below. To the left we have marked $B_1(0)$ and $\Gamma_u$. The middle picture 
shows $\hat{B}$ (dashed) $\Gamma_u$ with a cross marking where $\Gamma_u$ intersects $\partial \hat{B}$. The right picture shows $\Gamma_u$
translated by $\delta e_2$, from the picture it is clear that if we rotate the translated $\Gamma_u$ by $\tilde{\delta}\approx \delta$ radians, in the 
direction of the arrow, the rotated discontinuity set will intersect $\partial \hat{B}$ at the cross, the point, where $\Gamma_u$ enters $\hat{B}$. It follows that if we chose 
$\Gamma_w$ to be the rotated translation of $\Gamma_u$ then $\Gamma_w$ may be extended by $\Gamma_u$ into $B_2(0)\setminus \hat{B}$
to a rectifiable curve connecting the center of $\hat{B}$ to $\partial B_2(0)$. We define $w$ in $B_2$, such that 
$w=u$ in $B_2\setminus \hat{B}$.

\vspace{6mm}

\begin{center}
\includegraphics[height=4cm]{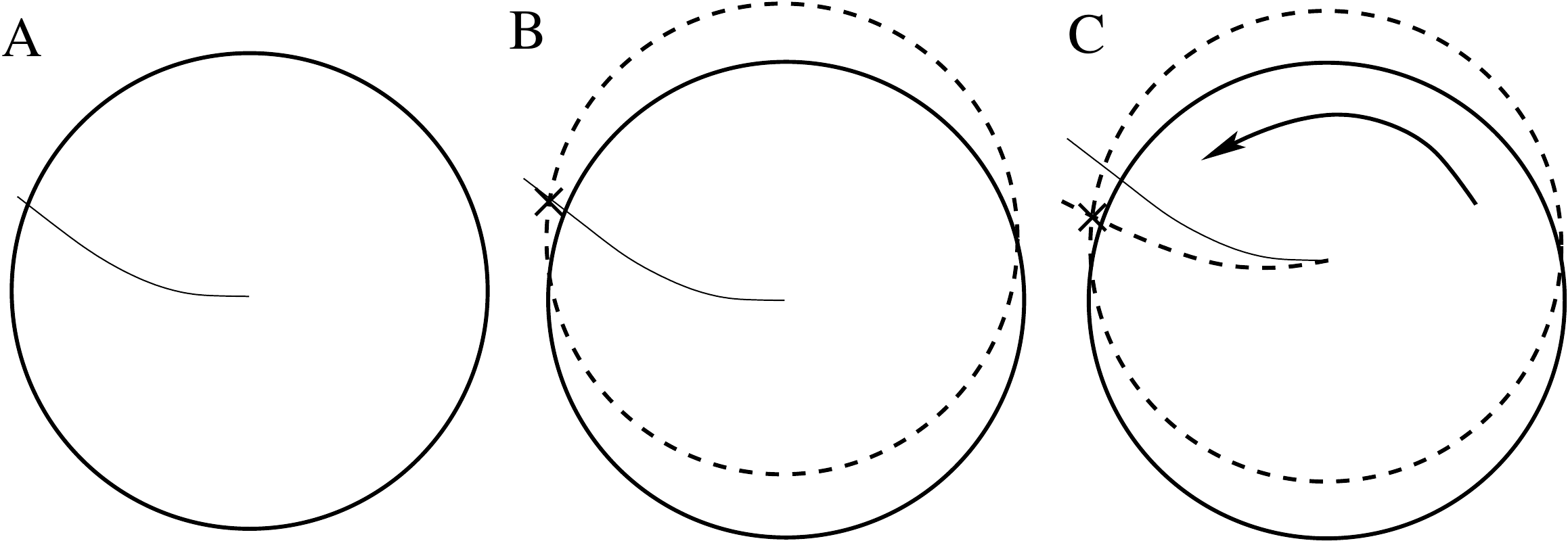}
\end{center}

\vspace{1mm}


\vspace{2mm}

Then, in order to prove the lemma, we need to compare the Mumford-Shah energies of $u$ and $w$ in $\hat{B}\setminus\Gamma_u$ and in $\hat{B}\setminus\Gamma_w$ respectively, and the main difficulty is the comparison
of the Dirichlet energies. The solution of this problem is the following: instead of comparing the Dirichlet energies with each other, we compare both with the Dirichlet energy of $u$ in $B_1(0)\setminus\Gamma_u$ (see Claim 3 and Claim 4 below). 

Below we present several claims, which prove the lemma. In Claim 1 we calculate the boundary values of $u$ on $\partial \hat{B}$. In Claim 2 we define the function $w$ in $\hat{B}$. 
In Claim 3 we compare the Dirichlet energy of $w$ in  $\hat{B}\setminus\Gamma_w$ with the Dirichlet energy of $u$ in $B_1(0)\setminus\Gamma_u$ using the fact that the domain $\hat{B}\setminus\Gamma_w$ is obtained from the domain $B_1(0)\setminus\Gamma_u$ by translation and rotation, thus $w$ and $u$ differ by a  term originating 
from the boundary values.  In Claim 4 we compare the Dirichlet energies of $u$ in the domains $B_1(0)\setminus\Gamma_u$ and $\hat{B}\setminus\Gamma_u$ with each other using direct calculations. Finally, in Claim 5 we 
compare the Hausdorff measures of $\Gamma_w$ and $\Gamma_u$ in $\hat{B}$.

We need to calculate the boundary values of $u$ on $\partial \hat{B}$.
First we need to 
express the values of $u$ on $\partial \hat{B}$ in polar coordinates 
$(\hat{r},\hat{\phi})$ that are centered at the center of $\hat{B}$ 
(note that $\hat{r}=1$ on $\partial \hat{B}$ and does not appear in the
formulas in Claim 1). We will also use $(\hat{x}_1,\hat{x}_2)$ for the
Cartesian coordinates with origin at the center on $\hat{B}$; 
$(\hat{x}_1,\hat{x}_2+\delta)=(x_1,x_2)$.

\vspace{1mm}

{\bf Claim 1.} {\sl Let $(x_1,x_2)\in \partial B_1(0)$ and then we may 
express the value of $u(x_1,x_2+\delta)$:}
$$
u(1,\hat{\phi})-R_0(x_1,x_1+\delta):=u(x_1,x_2+\delta)-R_0(x_1,x_1+\delta)= 
\sqrt{\frac{2}{\pi}}\sin(\hat{\phi}/2)+\epsilon\hat{\phi} \sin(\hat{\phi}/2) +
$$
$$
+\frac{\delta}{\sqrt{2\pi}}\cos(\hat{\phi}/2)-\frac{\epsilon \delta}{2}\sin(\hat{\phi}/2)+R_1(x)
$$
{\sl where $(\hat{r},\hat{\phi})$ are polar coordinates with origin at the center of $\hat{B}$ chooses so that the point $(x_1,x_2+\delta)$ corresponds to the point $(1,\hat{\phi})$. The rest term $R_1$ can be written $R_1=R_1^1+R_1^2$, 
where $R^1_1$  is odd in the $\hat{x}_2$ variable, orthogonal to $r^{1/2}\sin(\phi/2)$ and satisfies $|R_1^1|, |\nabla R_1^1|\le C\epsilon \delta$ in a 
neighborhood of $\partial \hat{B}\setminus \Gamma_u$, while $R_1^2$ satisfies $|R_1^2|, |\nabla R_1^2|\le C \delta^2$.}

\vspace{1mm}

{\sl Proof of Claim 1:} By making a Taylor expansion in the $x_2$ direction we may conclude that, 
having substituted $r=1$ in places where $r$ is not affected by a derivative,
\begin{equation}\label{SubtleBranchPoint1}
u(x_1,x_2+\delta)-R_0(x_1,x_2+\delta)=
\end{equation}
\begin{equation}\label{SubtleBranchPoint2}
=u(x_1,x_2)-R_0(x_1,x_2)+ \delta \sqrt{\frac{2}{\pi}}\frac{d}{dx_2}r^{1/2}\sin(\phi/2)+
\delta \epsilon \frac{d \mathfrak{z}}{d x_2}+O(\delta^2),
\end{equation}
We should remark something regarding 
the branch cut in the calculation (\ref{SubtleBranchPoint1}). Clearly the left side is discontinuous across $\Gamma_u$ wherefore the right side 
has to be discontinuous across $\Gamma_u$. Therefore we have to choose the 
branch cut in (\ref{SubtleBranchPoint2}) so that the right hand side is 
discontinuous at the point $\partial B_1(0) \cap \Gamma_u$, this is of course 
always possible.

Expressing 
$\frac{\partial }{\partial x_2}=\sin(\phi)\frac{\partial}{\partial r}+\frac{\cos(\phi)}{r}\frac{\partial}{\partial \phi}$ in
polar coordinates we can calculate
$$
\left. \frac{\partial r^{1/2}\sin(\phi/2)}{\partial x_2}\right\lfloor_{r=1} = \frac{1}{2}\sin(\phi)\sin(\phi/2) + \frac{1}{2}\cos(\phi)\cos(\phi/2)=
\frac{1}{2}\cos(\phi/2)
$$
and 
$$
\left. \frac{\partial \mathfrak{z}}{\partial x_2}\right\lfloor_{r=1} = \frac{1}{2}\phi \sin(\phi)\sin(\phi/2)-
$$
$$
-\frac{\ln(1)}{2}\sin(\phi)\cos(\phi/2)-\sin(\phi)\cos(\phi/2)+\cos(\phi)\sin(\phi/2)+
$$
$$
+\frac{1}{2}\phi\cos(\phi)\cos(\phi/2) +\frac{\ln(1)}{2}\cos(\phi)\sin(\phi/2)=
$$
$$
=\frac{\phi}{2}\cos(\phi/2)-\sin(\phi/2).
$$

It follows that on $\partial B_1(\delta e_2)$
$$
u(x_1,x_2+\delta)-R_0(x_1,x_2+\delta)=u(x_1,x_2)-R_0(x_1,x_2)+
\frac{\delta}{2}\sqrt{\frac{2}{\pi}}\cos(\phi/2)+
$$
$$
\epsilon\delta \left(\frac{\phi}{2}\cos(\phi/2)-\sin(\phi/2) \right)+O(\delta^2).
$$
We can conclude that 
$$
u(x_1,x_2+\delta)=u(x_1,x_2)+\frac{\delta}{2}\sqrt{\frac{2}{\pi}}\cos(\phi/2)+
$$
$$
\epsilon\delta \left(\frac{\phi}{2}\cos(\phi/2)-\sin(\phi/2) \right)+O(\delta^2)
+\left(R_0(x_1,x_2+\delta)-R_0(x_1,x_2) \right).
$$

It follows that on $\partial B_1(\delta e_2)$
$$
u(x_1,x_2+\delta)-R_0(x_1,x_1+\delta)=\sqrt{\frac{2}{\pi}}\sin(\phi/2)+\epsilon\phi\sin(\phi/2)+\frac{\delta}{2}\sqrt{\frac{2}{\pi}}\cos(\phi/2)+
$$
$$
\epsilon\delta \left(\frac{\phi}{2}\cos(\phi/2)-\sin(\phi/2) \right)+O(\delta^2).
$$

Noticing that 
$$
\int_{-\pi}^\pi \left(\frac{\phi}{2}\cos(\phi/2)-\frac{1}{2}\sin(\phi/2) \right)\sin(\phi/2)d\phi=0
$$
we may conclude that the claim holds with 
$$
R_1=\epsilon\delta \left(\frac{\phi}{2}\cos(\phi/2)-\frac{1}{2}\sin(\phi/2)\right)+\left(R_0(x_1,x_2+\delta)-R_0(x_1,x_2) \right)+O(\delta^2)=R^1_1+R_1^2,
$$
where we choose $R_1^1$ to be orthogonal to $\sin(\phi/2)$ and $R_1^2$ is of order 
$\delta^2$.

To finish the proof of Claim 1 we just notice that if the point $(x_1,x_2+\delta)$ corresponds to the angle 
$\hat{\phi}$ and $\phi$ corresponds to $(x_1,x_2)$ then $\hat{\phi}=\phi$

\vspace{3mm}

We want to define the comparison pair $(w,\Gamma_w)$ to $(u,\Gamma_u)$ by 
\begin{equation}\label{Euler}
w=\left\{\begin{array}{ll}
          w & \textrm{ in }\hat{B} \\
          u & \textrm{ in } B_2(0)\setminus \hat{B}
         \end{array}
\right.\quad\quad
\Gamma_{w}=\left\{\begin{array}{ll}
          \Gamma_{w} & \textrm{ in }\hat{B} \\
          \Gamma_u & \textrm{ in } B_2(0)\setminus \hat{B}
         \end{array}
\right.\quad\quad
\end{equation}
where $w$ will be chosen so that $\nabla w\in L^2(B_2\setminus \Gamma_{w})$ and $\Gamma_{w}$ is a rectifiable curve connecting 
$\partial B_2$ to the center of $\hat{B}$.

In order to get $\Gamma_{w}$ to connect $\partial B_2$ to the center of $\hat{B}$ we will use the free discontinuity set 
$\Gamma_{u(r,\phi-\tilde{\delta})}$ translated by $\delta$ in the $x_2$ 
direction, where 
$\tilde{\delta}=\delta+O(\epsilon \delta)$ is chosen so that $\Gamma_w$ connects the origin to $\partial B_1(0)$. We let $(\hat{r},\hat{\phi})$
be polar coordinates with respect to the center of $\hat{B}$, using $u(\cdot,\cdot)$ to denote the functional expression of $u$ in the coordinates $(r,\phi)$, we may write 
$$
w(\hat{r},\hat{\phi})=u(\hat{r},\hat{\phi}-\tilde{\delta})+\textrm{ error terms.}
$$
The idea is that with such a choice of $(w, \Gamma_w)$ we may explicitly calculate the error terms up to lower 
order and therefore the energy of $(w,\Gamma_w)$ in $\hat{B}$ (up to lower order)
and see that $(u,\Gamma)$ is not the minimizer. 

\vspace{1mm}

{\bf Claim 2:} {\sl If we let $w(\hat{r},\hat{\phi})=u(x_1,x_2+\delta)$ on $\partial\hat{B}$, where $(\hat{r},\hat{\phi})$
is choosen so that $(\hat{r},\hat{\phi})$ cooresponds to the point $(x_1,x_2+\delta)$ for $(x_1,x_2)\in \partial B_1(0)$, then we may (and will) define $w$ in the 
set $\hat{B}\setminus \Gamma_w$ according to}
\begin{equation}\label{CalcOfW}
w(\hat{r},\hat{\phi})=u(\hat{r},\hat{\phi}-\tilde{\delta})+
\sqrt{\frac{2}{\pi}}\delta \hat{r}^{1/2}\cos(\hat{\phi}/2)+
\end{equation}
$$
+\frac{\epsilon\delta}{2}\left(\hat{\phi}\hat{r}^{1/2}\cos(\hat{\phi}/2)+
\hat{r}^{1/2}\ln(\hat{r})\sin(\hat{\phi}/2)\right)
+\frac{\epsilon\delta}{2}\hat{r}^{1/2}\sin(\hat{\phi}/2)+R_1+R_2
$$
{\sl where $\|R_2\|_{C^{1,\alpha}}\le C(\sigma \epsilon \delta+\delta^2)$, 
$R_1$ is the error term from Claim 1, $R_1$ and $R_2$ are extended to be harmonic in $\hat{B}\setminus \Gamma_w$ with Neumann boundary condition on $\Gamma_w$.}

\vspace{1mm}

\textsl{Proof of Claim 2:} First we notice that, by (\ref{SoDamedBoring}),
$$
 u(r,\phi-\tilde{\delta})-R_0(r,\phi-\tilde{\delta})=\sqrt{\frac{2}{\pi}}r^{1/2}\sin((\phi-\tilde{\delta})/2)+
 $$
 $$
 +\epsilon\left(r^{1/2}(\phi-\tilde{\delta}) \sin((\phi-\tilde{\delta})/2)-r^{1/2}\ln(r)\cos((\phi-\tilde{\delta})/2) \right).
$$
Next we want to do a Taylor expansion on $\partial B_1(0)$. In particular it follows that
$$
u(r,\phi-\tilde{\delta})-R_0(r,\phi-\tilde{\delta})=
$$
\begin{equation}\label{TaylorTurn}
=u(r,\phi)-R_0(r,\phi)+\tilde{\delta}\frac{\partial (u(r,\phi)-R_0(r,\phi))}{\partial \phi}-\tilde{\delta}\frac{\partial R_0}{\partial \phi}+O(\delta^2)
\end{equation}
$$
= u(r,\phi)-R_0(r,\phi)-\tilde{\delta}\left( \sqrt{\frac{1}{2\pi}}r^{1/2}\cos(\phi/2)+\epsilon\sin(\phi/2)+
\frac{\epsilon}{2}\phi \cos(\phi/2)\right)+O(\delta^2),
$$
where we interpret the right hand side (which actually only contain
trigonometric and linear functions in $\phi$) so that the discontinuity in $\phi$ appears at the same $\phi$ as in the right  side. 
By choosing $R_2=R_0(x_1,x_2+\delta)-R_0(r,\phi-\tilde{\delta})$ it follows, since $R_0$ is smooth on $\partial B_1\setminus \Gamma_u$, that $\|R_2\|_{C^{1,\alpha}}<C\sigma \epsilon(\delta+\tilde{\delta}) <C\sigma \epsilon \delta$.

Using the expression in Claim 1 together with (\ref{TaylorTurn}) we may conclude that 
$$
u(x_1,x_2+\delta)-u(\hat{r},\hat{\phi}-\tilde{\delta})=
$$
$$
=\left( \frac{\tilde{\delta}+\delta}{\sqrt{2\pi}}+\frac{\epsilon\tilde{\delta}}{2}\right)\cos(\hat{\phi}/2)
+\epsilon\frac{2\tilde{\delta}-\delta}{2}\sin(\hat{\phi}/2)+R_1+R_2.
$$
Next we note that $\tilde{\delta}=\delta+O(\delta^2)$ which implies that we may,
after yet another Taylor expansion, use $\delta$ in place of $\tilde{\delta}$
with an error of at most $O(\delta^2)$ which we may include in the term $R_2$.
This yields that $u$ equals the expression in (\ref{CalcOfW}) on $\partial \hat{B}$ which means that with $w$ defined as in (\ref{CalcOfW}) in 
$\hat{B}\setminus \Gamma_w$ will be continuous across $\partial \hat{B}$.
It follows that $w\in W^{1,2}(B_2\setminus \Gamma_w)$ with the defnition 
of $w$ as in (\ref{Euler}). It follows that $w$ is an admissible variation.

\vspace{2mm}

{\bf Claim 3:} {\sl The Dirichet energy of $w$ may be calculated }
\begin{equation}\label{CalcOfWenergy}
\int_{\hat{B}\setminus \Gamma_w}|\nabla w|^2= 
\int_{B_1\setminus \Gamma_u}|\nabla u(r,\phi)|^2-
\epsilon \delta\sqrt{\frac{\pi}{2}}+ O(\epsilon^2\delta +\sigma \epsilon \delta).
\end{equation}

\textsl{Proof of Claim 3:} 
We may express $w(\hat{r},\hat{\phi}+\tilde{\delta})$, where $R_3$
includes the rest terms $R_1$ and $R_2$,
$$
w(\hat{r},\hat{\phi}+\tilde{\delta})=u(\hat{r},\hat{\phi})+
 \sqrt{\frac{2}{\pi}}\delta\cos(\hat{\phi}/2) +
$$
$$
+ \frac{\epsilon \delta}{2}\left(\hat{\phi}\cos(\hat{\phi}/2)+\hat{r}^{1/2}\ln(\hat{r})\sin(\hat{\phi}/2)\right)+
\frac{\epsilon \delta}{2}\sin(\hat{\phi}/2)+R_3,
$$
we would want to stress again that here we interpret $u(\cdot,\cdot)$ and the functional expression to the right in (\ref{SoDamedBoring}) that takes two real numbers as input and not a point in the plane; that is $u(\hat{r},\hat{\phi})$
takes the value we get when we substitute $(\hat{r},\hat{\phi})$ in (\ref{SoDamedBoring}) and not the value
of the function $u$ at the point with cartesian coordinates $(\hat{r}\cos(\hat{\phi}),\hat{r}\sin(\hat{\phi})+\delta)$ 
in the $(x_1,x_2)-$plane.

Using that $\Gamma_w$ is defined by a translation of $\Gamma_u$ by $\delta e_2$ and 
then a rotation by $\tilde{\delta}$ around the center of $\hat{B}$ we may 
calculate
$$
\int_{\hat{B}\setminus \Gamma_w}|\nabla w(\hat{r},\hat{\phi}+\tilde{\delta})|^2= 
\int_{B_1(0)\setminus \Gamma_u}|\nabla u(r,\phi)|^2+
$$
$$
+2\sqrt{\frac{2}{\pi}}\delta \int_{B_1\setminus \Gamma_u}\nabla u \cdot \nabla (r^{1/2}\cos(\phi/2)) +
$$
$$
+\frac{\epsilon \delta}{2}\int_{B_1(0)\setminus \Gamma_u}\nabla u \cdot \nabla \left(\phi r^{1/2}\cos(\phi/2)+r^{1/2}\ln(r)\sin(\phi/2)\right)+
$$
$$
+\frac{\epsilon \delta}{2}\int_{B_1\setminus \Gamma_u}\nabla u \cdot \nabla (r^{1/2}\sin(\phi/2))+O(\delta^2\sigma\epsilon+\epsilon^2\delta)=
$$
$$
= \int_{B_1\setminus \Gamma_u}|\nabla u(r,\phi)|^2+ I_1+I_2+I_3+O(\delta^2\sigma\epsilon+ \epsilon^2\delta),
$$
here we also used that $R_3=R_1^1+R_1^2+R_2+o(\epsilon\delta)$ and that $R_1^1$ is orthogonal to $r^{1/2}\sin(\phi/2)$ together with the 
estimates on $R_1^1$, $R_1^2$ and $R_2$ provided in Claim 1 and Claim 2.


First let us observe that

\vspace{3mm}

$
\nabla\left(r^{1/2}\cos(\phi/2)\right)=\dfrac{1}{2}r^{-1/2}(\cos(\phi/2),\sin(\phi/2)),
$
\vspace{3mm}

$
\nabla\left(r^{1/2}\sin(\phi/2)\right)=\dfrac{1}{2}r^{-1/2}(-\sin(\phi/2),\cos(\phi/2)),
$

\begin{multline*}
\nabla\left(r^{1/2}\ln(r)\cos(\phi/2)\right)=\\
\frac{1}{2}r^{-1/2}\ln(r)(\cos(\phi/2),\sin(\phi/2))+r^{-1/2}\cos(\phi/2)(\cos\phi,\sin\phi),
\end{multline*}

\begin{multline*}
\nabla\left(r^{1/2}\ln(r)\sin(\phi/2)\right)=\\
\frac{1}{2}r^{-1/2}\ln(r)(-\sin(\phi/2),\cos(\phi/2))+r^{-1/2}\sin(\phi/2)(\cos\phi,\sin\phi),
\end{multline*}

\begin{multline*}
\nabla\left(r^{1/2}\phi\cos(\phi/2)\right)=\\
\frac{1}{2}r^{-1/2}\phi(\cos(\phi/2),\sin(\phi/2))+r^{-1/2}\cos(\phi/2)(-\sin\phi,\cos\phi),
\end{multline*}

\begin{multline*}
\nabla\left(r^{1/2}\phi\sin(\phi/2)\right)=\\
\frac{1}{2}r^{-1/2}\phi(-\sin(\phi/2),\cos(\phi/2))+r^{-1/2}\sin(\phi/2)(-\sin\phi,\cos\phi).
\end{multline*}

\vspace{3mm}


Let us now calculate $I_1$, $I_2$ and $I_3$ in turn. We begin by calculating $I_1$

$$
I_1= 2\delta \sqrt{\frac{2}{\pi}}\int_{B_1\setminus \Gamma_u} \nabla u \cdot \nabla (r^{1/2}\cos(\phi/2))=
$$
$$
 2\delta \sqrt{\frac{2}{\pi}}\int_{B_1\setminus \Gamma_u} \nabla (r^{1/2}\sin(\phi/2))\cdot \nabla (r^{1/2}\cos(\phi/2))-
$$
$$
-2\delta\sqrt{\frac{2}{\pi}}\int_{\hat{B}_1\setminus \Gamma_w}\nabla \left( \frac{1}{2}r^{1/2}\ln(r)\cos(\phi/2)\right)\cdot 
\nabla \left( r^{1/2}\cos(\phi/2)\right)
$$
$$
+2\epsilon\delta\sqrt{\frac{2}{\pi}}
\int_{B_1\setminus \Gamma_u}\nabla (r^{1/2}\phi \sin(\phi/2))\cdot \nabla(r^{1/2}\cos(\phi/2))+C\sigma \epsilon\delta=
$$
$$
=0+0-2\epsilon\delta\sqrt{\frac{\pi}{2}}+O(\sigma \epsilon \delta),
$$
where we have used the following elementary integrals
$$
\int_{B_1\setminus \Gamma_u} \nabla (r^{1/2}\sin(\phi/2))\cdot \nabla (r^{1/2}
\cos(\phi/2))=
\int_{B_1\setminus \Gamma_u} 0 =0,
$$
$$
\int_{B_1\setminus \Gamma_u}\nabla \left(r^{1/2}\ln(r)\cos(\phi/2)\right)\cdot 
\nabla \left( r^{1/2}\cos(\phi/2)\right)=
$$
$$
= \int_{-\pi}^\pi \int_0^1 \left( \frac{\ln(r)}{4}+\frac{1}{2}\cos^2(\phi/2)\right)drd\phi=0
$$
and 
$$
\int_{B_1\setminus \Gamma_u}\nabla (r^{1/2}\phi \sin(\phi))\cdot \nabla(r^{1/2}\cos(\phi/2))=
-\frac{1}{4}\int_{-\pi}^\pi \int_0^1 (1-\cos(\phi))drd\phi=-\frac{\pi}{2}.
$$

To estimate $I_2$ we use that 
$u(r,\phi)=\sqrt{\frac{2}{\pi}}r^{1/2}\sin(\phi/2)+O(\epsilon)$ to calculate
$$
I_2=\frac{\epsilon \delta}{2}\int_{B_1\setminus \Gamma_u}\nabla u \cdot \nabla \left(\phi r^{1/2}\cos(\phi/2)+r^{1/2}\ln(r)\sin(\phi/2)\right)=
$$
$$
=\frac{\epsilon \delta}{2}\int_{B_1\setminus \Gamma_u}\nabla \left(\sqrt{\frac{2}{\pi}}r^{1/2}\sin(\phi/2) \right) \cdot \nabla \left(\phi r^{1/2}\cos(\phi/2)\right)+
$$
$$
+\frac{\epsilon \delta}{2}\int_{B_1\setminus \Gamma_u}\nabla \left(\sqrt{\frac{2}{\pi}}r^{1/2}\sin(\phi/2) \right) \cdot \nabla \left(r^{1/2}\ln(r)\sin(\phi/2)\right)+
O(\epsilon^2\delta)
$$
$$
=\frac{\epsilon \delta}{8}\sqrt{\frac{2}{\pi}}\int_{-\pi}^\pi \int_0^1 (1+\cos(\phi))drd\phi
+
\frac{\epsilon\delta}{2}\sqrt{\frac{2}{\pi}}\int_{-\pi}^\pi \int_0^1(\frac{1}{4}\ln(r)+\frac{1}{2}\sin^2(\phi/2))drd\phi
$$
$$
=\frac{\epsilon\delta}{2}\sqrt{\frac{\pi}{2}}+O(\epsilon^2\delta).
$$

For $I_3$ we use 
$$
I_3=\frac{\epsilon\delta}{2}\int_{B_1\setminus \Gamma_u}\nabla u\cdot \nabla r^{1/2}\sin(\phi/2)=
$$
$$
=\frac{\epsilon\delta}{2}\sqrt{\frac{2}{\pi}}\int_{B_1\setminus \Gamma_u}|\nabla r^{1/2}\sin(\phi/2)|^2+ O(\epsilon^2\delta)=
\frac{\epsilon\delta}{2}\sqrt{\frac{\pi}{2}}+O(\epsilon^2\delta).
$$

Putting the estimates for $I_1$, $I_2$ and $I_3$ together we may conclude that 
$$
\int_{\hat{B}\setminus \Gamma_w}|\nabla w|^2= 
\int_{B_1\setminus \Gamma_u}|\nabla u(r,\phi)|^2-\epsilon \delta \sqrt{\frac{\pi}{2}}+ O(\epsilon^2\delta +\sigma \epsilon \delta)
$$
this proves (\ref{CalcOfWenergy})


\vspace{2mm}

{\bf Claim 4:} {\sl We may calculate the Dirichlet energy of $u$}
$$
\int_{\hat{B}\setminus \Gamma_u}|\nabla u|^2= \int_{B_1(0)\setminus \Gamma_u}|\nabla u|^2+\epsilon\delta \sqrt{\frac{\pi}{2}}+O(\sigma\epsilon \delta). 
$$

\vspace{1mm}

\textsl{Proof of Claim 4:} Since $\hat{B}=B_{1}(\delta e_2)$ we may make a Taylor expansion in $\delta$ and calculate
$$
\int_{\hat{B}\setminus \Gamma_u}|\nabla u|^2=\int_{B_1(0)\setminus \Gamma_u}|\nabla u|^2+
\delta \left.\frac{d}{d\delta}\int_{B_1(\delta e_2)}|\nabla u|^2\right\lfloor_{\delta=0} +O(\delta^2)=
$$
$$
=\int_{B_1(0)\setminus \Gamma_u}|\nabla u|^2+\delta\int_{\partial B_1(0)}\sin(\phi)|\nabla u|^2+ O(\delta^2).
$$
We need to calculate the boundary integral in the final expression.

A trivial calculation, using the expression for $u$ from (\ref{SoDamedBoring}), shows that
$$
\int_{\partial B_1(0)}\sin(\phi)|\nabla u|^2=\frac{2}{\pi}\int_{\partial B_1}\sin(\phi)|\nabla r^{1/2}\sin(\phi/2)|^2+
$$
$$
+ \epsilon\int_{\partial B_1}\sin(\phi)\nabla \left( \sqrt{\frac{2}{\pi}}r^{1/2}\sin(\phi/2)\right)\cdot 
\nabla \left( r^{1/2}\phi\sin(\phi/2)-r^{1/2}\ln(r)\cos(\phi/2)\right)+
$$
$$
+O(\sigma \epsilon+ \epsilon^2)=\epsilon\sqrt{\frac{2}{\pi}}\int_{-\pi}^\pi \frac{\phi\sin(\phi)(\sin^2(\phi/2)+\cos^2(\phi/2))}{4}+
$$
$$
+O(\sigma\epsilon+\epsilon^2)=\epsilon\sqrt{\frac{\pi}{2}}+O(\sigma \epsilon+\epsilon^2).
$$
Claim 4 follows.

\vspace{2mm}

{\bf Claim 5:} {\sl We may estimate }
$$
\H^1(\Gamma_u\cap \hat{B})-\H^1(\Gamma_v\cap \hat{B})= R_\Gamma(\epsilon,\delta,\sigma),
$$
{\sl where the rest $\frac{|R_\Gamma|}{|\epsilon\delta|}\to 0$ as $\epsilon,\delta,\sigma\to 0$.}

\vspace{1mm}

{\sl Proof of Claim 5:} Since $\Gamma_w\cap \hat{B}$ is a rigid motion of $\Gamma_u\cap B_1(0)$ it clearly follows that 
$\H^1(\Gamma_u\cap B_1)=\H^1(\Gamma_v\cap \hat{B})$. Therefore we need to estimate the difference $\H^1(\Gamma_u\cap \hat{B})-\H^1(\Gamma_u\cap B_1(0))$.
The difference between $\H^1(\Gamma_u\cap B_1)$ and $\H^1(\Gamma_u\cap \hat{B})$ will depend on the part of $\Gamma_u$ contained in the symmetric difference 
$B_1(0) \Delta \hat{B}$.

\vspace{2mm}

\begin{center}

\includegraphics[height=5cm]{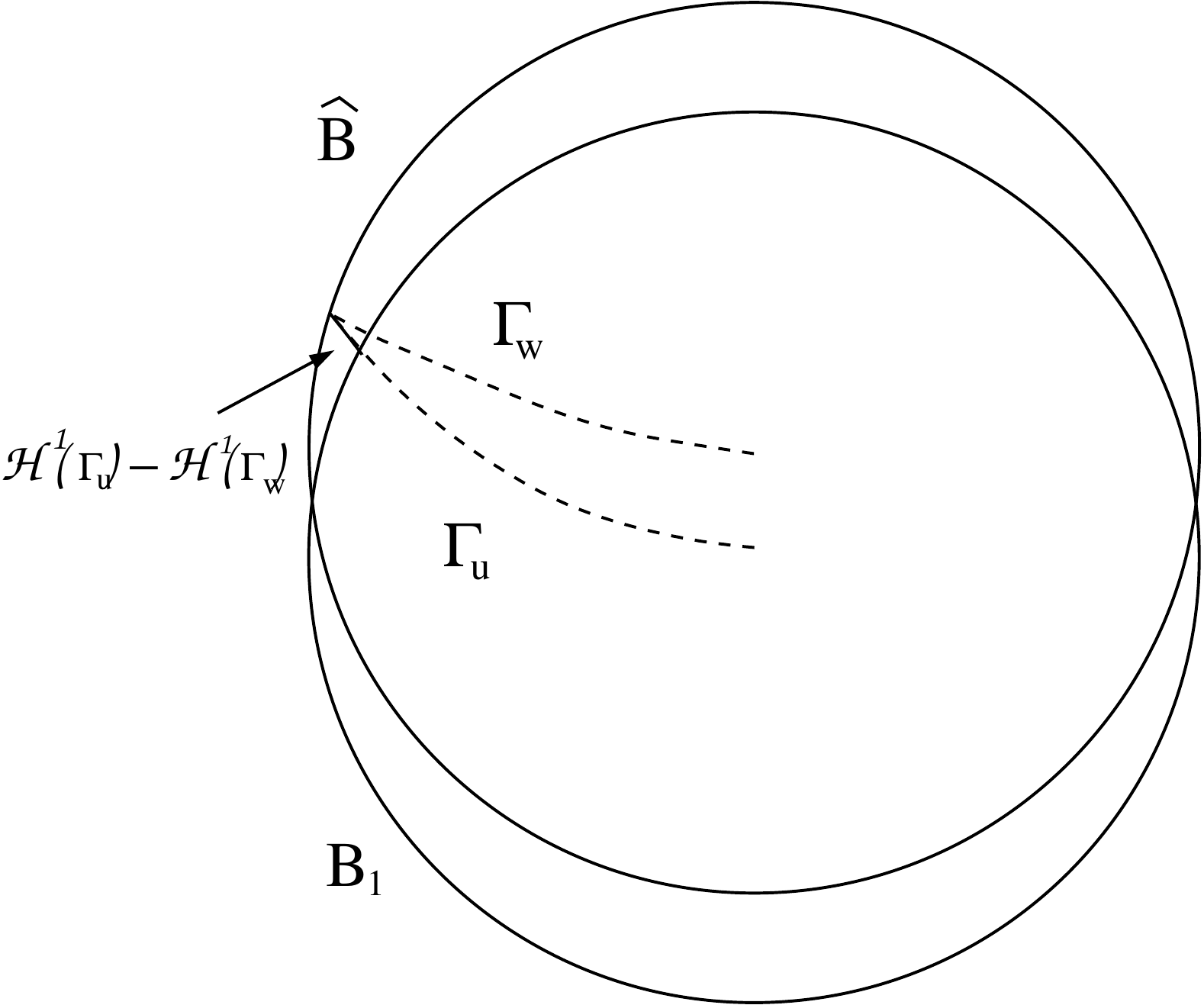}

\end{center}

\vspace{2mm}

{\bf Figure:} {\sl The above figure shows the geomerty of the situation in claim 5. The difference that we need to estimate is the part of 
$\Gamma_u$ contined in the ball $\hat{B}$ (upper ball) but not in $B_1(0)$.}

\vspace{2mm}

We let $\Sigma$ be the projection of $\Gamma_u\cap (B_1(0)\Delta \hat{B})$ onto the $x_1-$axis:
$$
\Sigma = \{x_1;\; (x_1,\epsilon f(x_1))\in \Gamma_u\cap (B_1(0)\Delta \hat{B})\} 
$$
By Corollary \ref{FconvC1alpha} it follows that $\epsilon |f'|\le C\epsilon$ for $x_1\in (-3/2, -1/2)$ if $\epsilon$ is small enough.
It follows that 
\begin{equation}\label{GoeasOn}
\left|\H^1(\Gamma_u\cap \hat{B})-\H^1(\Gamma_v \cap \hat{B})\right|=\int_\Sigma \sqrt{1+\epsilon^2 |f'|^2}dx_1\le (1+C\epsilon^2) \H^1(\Sigma).
\end{equation}

To estimate $\H^1(\Sigma)$ we let $x_0\in (-1,0)$ be the point where $(x_0,\epsilon f(x_0))\in \partial B_1(0)$. We may then approximate 
$\partial B_1(0)$ around $(x_0,\epsilon f(x_0))$ by the tangent $T$ to the circle:
\begin{equation}\label{TangentOfB1}
\frac{\sqrt{1-\epsilon^2f(x_0)^2}}{\epsilon \delta f(x_0)}(x-x_0)=y-\epsilon f(x_0).
\end{equation}
We may also approximate $\partial \hat{B}$ by the shifted tangent line $\hat{T}$
\begin{equation}\label{TangentOfBHat}
\frac{\sqrt{1-\epsilon^2f(x_0)^2}}{\epsilon \delta f(x_0)}(x-x_0)=y-\epsilon f(x_0)-\delta.
\end{equation}

\vspace{2mm}

\begin{center}

\includegraphics[height=5cm]{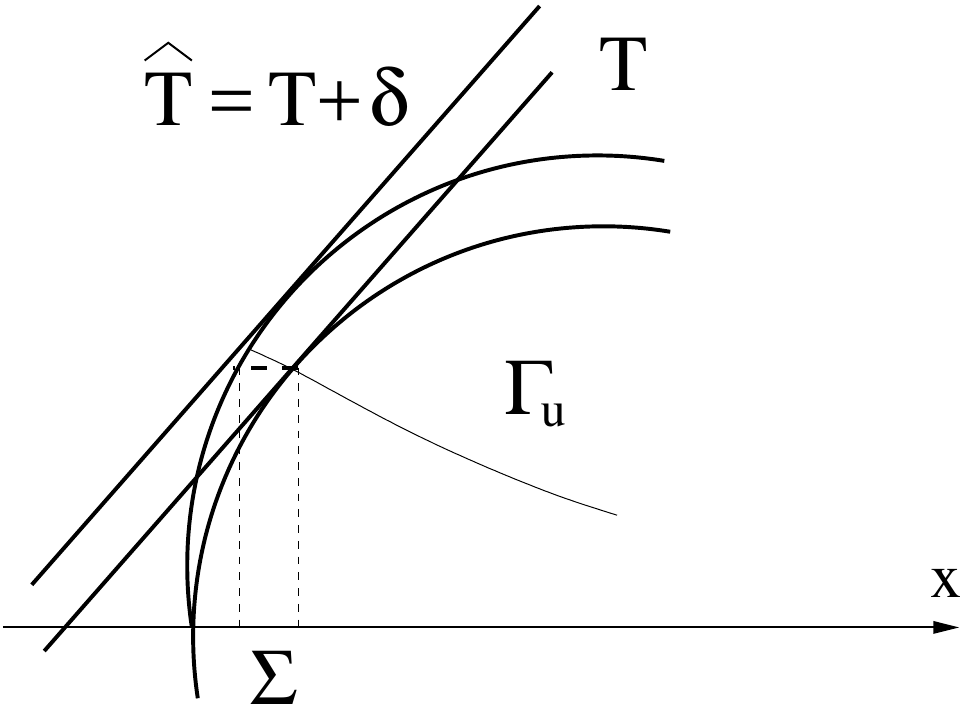}

\end{center}

\vspace{2mm}

{\bf Figure:} {\sl The above picture shows the geometry of the estimate of $\H^1(\Sigma)$ in claim 5. We want to estimate the length of $\Sigma$.
In order to do that we approximate $\partial B_1(0)$ and $\partial \hat{B}$ by the tangent lines $T$ given by (\ref{TangentOfB1})
and $\hat{T}=T+\delta$ given by (\ref{TangentOfBHat}). Since the slope of $T$ is explicitly calculable it is easy to estimate $\H^1(\Sigma)$.}

\vspace{2mm}

The tangent $\hat{T}$ is $T$ shifted by $\delta$ in the $x_2$ direction which is the same as shifting it by $\frac{-\epsilon f(x_0)}{\sqrt{1-\epsilon^2f(x_0)^2}}\delta$
in the $x_1$ direction. It follows that the length 
\begin{equation}\label{StillBoring}
\H^1(\Sigma) = \left|\frac{\epsilon f(x_0)}{\sqrt{1-\epsilon^2f(x_0)^2}}\delta\right|+O(\delta^2)
\end{equation}
with an error of lower order that comes from the Taylor expansions.

By Corollary \ref{FconvC1alpha} it follows that $f(-1)= o(1)$ as $\sigma \to 0$. Using that $x_0\approx -1$ with an error of order $\delta^2$ we may 
conclude that 
\begin{equation}\label{EstClaim5}
\left|\H^1(\Gamma_u\cap \hat{B})-\H^1(\Gamma_u\cap B_1(0))\right|\le (1+C\epsilon^2)\H^1(\Sigma)\le 
\end{equation}
$$
\le C \epsilon \delta f(-1)+O(\delta(\epsilon^3+\delta))=o(\epsilon\delta)
$$
where we also used (\ref{GoeasOn}) and (\ref{StillBoring}). Since the right side in (\ref{EstClaim5}) goes to 
zero faster than $\epsilon\delta$ as $\epsilon,\delta,\sigma \to 0$ it follows that the error term 
$R_\Gamma$ satisfies the condition in the claim.

It is easy to see that the sign of $\H^1(\Gamma_u\cap \hat{B})-\H^1(\Gamma_u\cap B_1(0))$ must be the same as the sign of $\epsilon\delta$, 
for instance from the figures in this part of the proof where $\epsilon>0$ and $\delta>0$. This proves Claim 5.

\vspace{2mm}

We are now ready to show that $(u,\Gamma_u)$ is not a minimizer in $\hat{B}$. In particular:
$$
J(u,\Gamma, \hat{B})-J(w,\Gamma_w, \hat{B})=
$$
$$
=\int_{\hat{B}\setminus \Gamma_u}|\nabla u|^2-\int_{\hat{B}\setminus \Gamma_w}|\nabla w|^2+
\H^1(\Gamma_u\cap \hat{B})-\H^1(\Gamma_w\cap \hat{B})=
$$
$$
\int_{\hat{B}\setminus \Gamma_u}|\nabla u|^2-\int_{B_1\setminus \Gamma_u}|\nabla u|^2+
\int_{B_1\setminus \Gamma_u}|\nabla u|^2-\int_{\hat{B}\setminus \Gamma_w}|\nabla w|^2+
\H^1(\Gamma_u\cap \hat{B})-\H^1(\Gamma_w\cap \hat{B})
$$
$$
=\epsilon\delta\sqrt{\frac{\pi}{2}}+ \epsilon \delta \sqrt{\frac{\pi}{2}}+R(\sigma, \epsilon, \delta)
=2\epsilon \delta\sqrt{\frac{\pi}{2}}+R(\sigma, \epsilon, \delta),
$$
where $\frac{|R|}{|\epsilon\delta|}\to 0$ as $\epsilon,\delta,\sigma\to 0$, we also used Claim 3 and Claim 4 to estimate the Dirichlet energy of $w$
and $\hat{u}$ and Claim 5 to estimate the difference in Hausdorff measure. This shows that $J(w,\Gamma_w,\hat{B})<J(u,\Gamma_u,\hat{B})$, if 
$\epsilon \delta>0$ and $\epsilon,\delta$ and $\sigma$ are small enough, contradicting that $(u,\Gamma_u)$ is a minimizer. \qed


\begin{prop}\label{NoLOTinLim}
 As in Proposition \ref{Linearization} we let $(u^j,\Gamma_j)$ be a sequence of minimizers to the Mumford-Shah problem that are $\epsilon_j-$close to a
 crack tip for some sequence $\epsilon_j\to 0$. Then $v^j=\frac{u^j-\Pi(u^j)}{\epsilon_j}\to v^0$, $f_j\to f_0$ in the same sense as in 
 Lemma \ref{EstIfThereExistsBound}. Then, by choosing the coordinate systems appropriately,
 \begin{equation}
 v^0(r,\phi)=a+\sum_{k=2}^\infty a_k r^{\alpha_k}\cos(\alpha_k \phi)+ \sum_{k=2}^\infty b_k r^{k-1/2}\sin((k-1/2)\phi)
 \end{equation}
 and
 $$
 f_0(x_1)=\sum_{k=2}^\infty 2a_k \sqrt{\frac{\pi}{2}}\sin(\alpha_k \pi)|x_1|^{\alpha_k+\frac{1}{2}}
 $$
\end{prop}

The interesting thing with this proposition is that the limit does not contain the $\mathfrak{z}$ and $\mathfrak{h}$ terms and that the summations 
start from $k=2$ (cf. Lemma \ref{ImprovedLinearization}).

\textsl{Proof of Proposition \ref{NoLOTinLim}:} By Lemma \ref{ImprovedLinearization} we know that, by choosing 
the coordinate system appropriately,
\begin{equation}
 v^0(r,\phi)=a+a_0\mathfrak{z}(r,\phi)+\sum_{k=2}^\infty a_k r^{\alpha_k}\cos(\alpha_k \phi)+ \sum_{k=2}^\infty b_k r^{k-1/2}\sin((k-1/2)\phi)
\end{equation}
 and
$$
 f_0(x_1)=a_0\mathfrak{h}(|x_1|)+\sum_{k=1}^\infty 2 a_k \sqrt{\frac{\pi}{2}}\sin(\alpha_k \pi)|x_1|^{\alpha_k+\frac{1}{2}}.
$$
We need to show that $a_0=0$ - we will therefore assume that $a_0\ne 0$ and derive a contradiction. If we rescale the functions $u^j\mapsto u^j_s(x)=\frac{u^j(sx)}{\sqrt{s}}$
then the corresponding linearilized sequence of functions 
$$
v^j_s=\frac{\frac{u^j(sx)}{\sqrt{s}}-\Pi(u^j,s)}{\epsilon_j}
$$  
and  $f^j_s$ (corresponding to  $\Gamma_{u^j_s}$) will converge to
\begin{equation}
 v^0_s(r,\phi)=a_0 v+\sum_{k=2}^\infty a_k s^{\alpha_k-1/2}r^{\alpha_k}\cos(\alpha_k \phi)+ \sum_{k=2}^\infty b_k s^{k-1}r^{k-1/2}\sin((k-1/2)\phi)=
\end{equation}
$$
=a_0\mathfrak{z}+ O\left(s\right)
$$
and
$$
 f_0(x_1)=a_0 \mathfrak{h}+\sum_{k=2}^\infty a_k s^{\alpha_k-1/2}2\sqrt{\frac{\pi}{2}}\sin(\alpha_k \pi)|x_1|^{\alpha_k+\frac{1}{2}}=
$$
$$
=a_0\mathfrak{h}+s^{\alpha_2-1/2}\left(2a_1\sqrt{\frac{\pi}{2}}\sin(\alpha_2 \pi)|x_1|^{\alpha_2+\frac{1}{2}} \right)+O(s^{\alpha_3-1/2})
$$
where the $O(s^{\alpha_3-1/2})$ are to be understood as a $W^{1,2}-$function with norm controlled by $Cs^{\alpha_2-1/2}$
for some fixed $C$.

In particular, by the strong convergence $\Reg(\nabla v^j)\to \Reg(\nabla v^0)$ Proposition \ref{SatisfiesRightEstProp}, we can deduce that
\begin{equation}\label{TheChoiceIsClear}
u_s^j(x)=\sqrt{\frac{2}{\pi}}r^{1/2}\sin(\phi/2)+\epsilon_j a_0 \mathfrak{z} + R(j,s),
\end{equation}
where $R(j,s)$ is a function with $\|\Reg(\nabla R)\|_{L^2(B_1)}=o(\epsilon_j)+O(s^{\alpha_2-1/2})$.

If we introduce the notation $\hat{\epsilon}_j=\epsilon_j a_0$ we may write (\ref{TheChoiceIsClear})
$$
u_s^j(x)=\sqrt{\frac{2}{\pi}}r^{1/2}\sin(\phi/2)+\hat{\epsilon} \mathfrak{z} + \hat{R}(j,s),
$$
where $\hat{R}(j,s)$ is a function with $\|\Reg(\nabla \hat{R})\|_{L^2(B_1)}=o(\hat{\epsilon}_j)+O(s^{\alpha_2-1/2})$.

By choosing $s$ small enough and $j$ large enough we can conclude that $u^j_s$ satisfies the 
assumptions of Lemma \ref{LemVarInOrth}. 

Similarly, we may write
$$
\Gamma_{u^j_s}=\{(x_1,\hat{\epsilon}_j f^j_s(x_1));\; x_1\in (-1,0)\}
$$
where
$$
f^j_s(x_1)=\sqrt{2\pi} x_1+\sqrt{2\pi} x_1 \ln(-x_1)+o(\hat{\epsilon}_j),
$$
as in Lemma \ref{LemVarInOrth}.

But then Lemma \ref{LemVarInOrth} implies that $(u^j_s,\Gamma_{u^j_s})$ is not a minimizer which is a contradiction. \qed



\section{$C^{1,\alpha}-$Regularity of the Crack tip.}

We are now ready to prove the first regularity theorem at of the crack-tip. The proof is entirely standard. We begin by showing regularity 
improvement of a solution close to a crack-tip.

\begin{prop}\label{FlatnessImprovement}
 For every $\alpha < 1/2$ there exists an $\epsilon_\alpha>0$ such that if
 $(u,\Gamma)$ is $\epsilon-$close to a crack-tip solution for some $\epsilon<\epsilon_\alpha$. Then, with $s_\alpha$ as in Corollary \ref{flatimpForLin},
 $(u_{s_\alpha}, \Gamma_{u_{s_\alpha}})$ is $s_\alpha^{\alpha} \epsilon-$close to a crack-tip solution:
 $$
 \left\|\nabla \left(\frac{u(s_\alpha x)}{\sqrt{s_\alpha}}-\Pi\left(\frac{u(s_\alpha x)}{\sqrt{s_\alpha}}\right)\right)\right\|_{L^2(B_1\setminus \Gamma_{u_{s_\alpha}})}\le
 s_\alpha^{\alpha}\left\|\nabla \left( u-\Pi(u)\right)\right\|_{L^2(B_1\setminus \Gamma_u)}.
 $$
\end{prop}
\textsl{Proof:} The proof is almost trivial. We argue by contradiction and assume that there exists a
sequence of minimizers $(u^j, \Gamma_{u^j})$ that are $\epsilon_j\to 0$ close to a crack-tip solution but
$(u^j_{s_\alpha}, \Gamma_{u^j_{s_\alpha}})$ is not $s_\alpha^{\alpha} \epsilon-$close to any crack-tip solution.

Denoting, as in Proposition \ref{Linearization},
$$
v^j=\frac{u^j-\Pi(u^j)}{\epsilon_j}
$$
we see that the assumption in the previous paragraph implies that
\begin{equation}\label{ArgumentofIvar}
\|\nabla v^j\|_{L^2(B_{s_\alpha}\setminus \Gamma_{u_j})}\ge s_\alpha^{\alpha}.
\end{equation}

If we use appropriately rotated coordinates, then by Lemma \ref{ImprovedLinearization}, $v^j\to v^0$ strongly (by Proposition \ref{SatisfiesRightEstProp})
where $v^0$ solves the linear system in Proposition \ref{AnalysisOfLinear}. By Proposition \ref{NoLOTinLim} $a_0=0$ and thus by Corollary \ref{flatimpForLin}
$$
\|\nabla v^0\|_{L^2(B_{s_\alpha}\setminus \Gamma_{0})}< s_\alpha^\alpha\|\nabla v^0\|_{L^2(B_{1}\setminus \Gamma_{0})}=s_\alpha^\alpha.
$$
This contradicts the strong convergence and (\ref{ArgumentofIvar}). This finishes the proof. \qed

\begin{thm}\label{C1alpha}
For every $\alpha< 1$ there exists an $\epsilon_\alpha>0$ such that if $(u,\Gamma)$ is $\epsilon$-close to a crack-tip solution for 
some $\epsilon\le \epsilon_\alpha$ then $\Gamma_u$ is $C^{1,\alpha}$ at the crack-tip.

 This in the sense that the tangent at the crack-tip is a well defined line, which we may, upon rotation of the coordinates, assume to be $\{(x_1,0);\; x_1\in \R\}$,
 and there exists a constant $C_\alpha$ such that
 $$
 \Gamma_u\subset \left\{(x_1,x_2);\; |x_2|<C_\alpha \epsilon |x_1|^{1+\alpha},\; x_1<0\right\}.
 $$
 Here $C_\alpha$ may depend on $\alpha$ but not on $\epsilon<\epsilon_0$.
\end{thm}

\textsl{Proof:} Denoting $t=s_\alpha$ it follows directly from an iteration of Theorem \ref{FlatnessImprovement}
that if $\epsilon$ is small enough then
\begin{equation}\label{Idiot}
\left\|\nabla \left(u_{t^k}-\Pi(u_{t^k})\right)\right\|_{L^2(B_1\setminus \Gamma_{u_{t^k}})}\le t^{(\alpha-1/2) k} \epsilon.
\end{equation}

Since, by the triangle inequality and the fact that the projection decreases the norm
\begin{equation}\label{StupidTriangle}
\left\|\nabla \left( \Pi(u_{t^{k+1}})-\Pi(u_{t^{k}})\right)\right\|_{L^2(B_1\setminus \Gamma_0)}\le
\end{equation}
$$
\le C \left\|\nabla \left(u_{t^k}-\Pi(u_{t^{k}})\right)\right\|_{L^2(B_1\setminus \Gamma_{u_{t^k}})}\le C t^{(\alpha-1/2) k} \epsilon
$$
for some dimensional constant. We see that $\Pi(u_{t^k})$ forms a Cauchy sequence in 
$W^{1,2}(B_1(0)\setminus \Gamma_0)$ and thus
$\lim_{k\to \infty}\Pi(u_{t^k})=\Pi_0$ exists. By a choice of coordinate system we may assume that
\begin{equation}\label{RotatedLimit}
\Pi_0=\sqrt{\frac{2}{\pi}}r^{1/2}\sin(\phi/2).
\end{equation}
Also, from the triangle inequality and (\ref{StupidTriangle})
$$
\left\|\nabla \left( \Pi(u_{t^k})-\Pi_0\right)\right\|_{L^2}\le C\epsilon \sum_{j=k}^\infty t^{\alpha j}\le C\epsilon t^{(\alpha-1/2) k}.
$$

From Corollary \ref{L2boundsCor} we can also conclude that the free discontinuity set $\Gamma_{u_{t^k}}$ is within distance
$C\epsilon t^{(\alpha-1/2) k}$ from a line $l_{t^k}=\{(x_1,x_2);\; x_2=h_k x_1\}$.

The coefficients $h_k$ and $h_{k+1}$ will differ by a less than $C \epsilon t^{(\alpha-1) k}$ for some constant $C$. This since
$$
\left\{ x;\, x\in B_{t^{k+1}},\, \textrm{dist}(x,l_{t^{k+1}})\le C\epsilon t^{(\alpha-\frac{1}{2}) (k+1)} \right\}\subset
$$
$$
\subset \left\{ x;\, x\in B_{t^{k}},\, \textrm{dist}(x,l_{t^{k}})\le C\epsilon t^{(\alpha-\frac{1}{2}) k} \right\}
$$
only if
\begin{equation}\label{242}
|h_k-h_{k+1}|\le C\epsilon t^{(\alpha-1/2) k}
\end{equation}

The coefficients  $h_{k}$ therefore forms a Cauchy sequence and therefore converges.
By the rotation made in (\ref{RotatedLimit}) we may conclude that $h_k\to 0$. Moreover, by the triangle inequality and (\ref{242})
\begin{equation}\label{AKBOUND}
|h_k|\le \sum_{j=k}^\infty|h_j-h_{j+1}|\le C \epsilon t^{(\alpha-1/2) k}.
\end{equation}

It follows from (\ref{AKBOUND}) and Corollary \ref{L2boundsCor} that
\begin{equation}\label{conclusion}
\Gamma_u \cap \left(B_{t^{k}}(0)\setminus B_{t^{k+1}}(0)  \right)\subset \{(x_1,x_2);\; |x_2|\le C \epsilon t^{(\alpha+1/2) k}\},
\end{equation}
where the change from $\alpha-1/2$ to $\alpha+1/2$ in the exponent is due to a scaling factor. The inclusion
(\ref{conclusion}) concludes the proof. \qed


\section{$C^{2,\alpha}$-Regularity of the Crack Tip.}\label{SecC2alpha}

In this section we prove the $C^{2,\alpha}$ regularity at the crack-tip. The proof is again rather long and we will split it into several 
minor results. However, most of the proof consists in estimating the error term in the linearization. The techniques that 
where used to prove strong convergence, with only minor modifications, can also be used to estimate the error term. Therefore this section
will mostly be a repetition of previous parts of the article - at times almost direct copies of previous proofs. Since the calculations are 
so similar to other calculations in this paper we hope that the reader can forgive us for freely referring to previous proofs. Even though this
puts more demands on the reader at times it also allows us to refrain from too much repetition of previous calculations. 

We begin with a definition of a projection operator $\P$ that will have the same role in this section as $\Pi$ had in previous sections.

\begin{definition}\label{DefProjOnLin}
 If $(u,\Gamma)$ is a minimizer of the Mumford-Shah energy,
 we will denote by $\P(u,\Gamma)$, or at times just $\P(u)$, the projection, w.r.t. $\|\nabla \cdot \|_{L^2(B_1\setminus \Gamma)}$, of the function 
 $u-\Pi(u,1)$ to the set of solutions to 
 the linearized system (\ref{theLinearSyst}). That is if $\Pi(u,1)=\lambda(0)r^{1/2}\sin(\phi+\phi_0)/2$, then $\P(u)$ is the function of the form
 $$
 \P(u)= a+a_0 \mathfrak{z}(r,\phi+\phi_0)+\sum_{k=1}^\infty a_k r^{\alpha_k}\cos\alpha_k( \phi+\phi_0)+ \sum_{k=2}^\infty b_k r^{k-1/2}\sin((k-\tfrac{1}{2})(\phi+\phi_0))
 $$
 that minimizes the following integral
 $$
 \int_{B_1\setminus \Gamma} \left|\nabla (u-\Pi(u,1))-\P(u)\right|^2.
 $$
\end{definition}

Throughout this section we will be using the notations which we would like to 
introduce in the following remark.

\begin{rem}\label{Rem:proj-fj}
 Let $(u,\Gamma)$ be a minimizer of the Mumford-Shah energy and let us assume that $\phi_0=0$, that is 
 we assume that there exists a constant $\lambda(u)\in \R$ such that
$$
\Pi(u,1)=\lambda(u)r^{1/2}\sin(\phi/2).
$$ 
We will also denote
$$
\epsilon=\|\nabla \P(u)\|_{L^2(B_1\setminus \Gamma)}<<1.
$$ 
Using asymptotic expansion (\ref{SumofHetros}) we can construct a function $f$ such that 
$(\P(u), \epsilon f)$ turns into a solution of the system (\ref{theLinearSyst}).
Moreover, we can define the constant $\delta>0$ by the equation 
$$
u=\Pi(u)+\P(u)+\epsilon\delta w,
$$
where $\|\nabla w\|_{L^2(B_1\setminus \Gamma_{u})}=1$, and the function $g$ by the representation
$$
\Gamma_{u}=\{(x_1,\epsilon(f+\delta g)),\;\textrm{ for }x_1<-\tau \},
$$
where $\tau$ and $\epsilon$ are related by Lemma \ref{FirstBadC1}.
\end{rem}

The next proposition is very similar to Proposition \ref{Linearization}.

\begin{prop}\label{NLinearization}
Let   
$$
u^j=\Pi(u^j)+\P(u^j)+\epsilon_j\delta_j w^j
$$ 
and
$$
\Gamma_{u^j}=\{(x_1,\epsilon_j(f_j+\delta_j g_j)),\;\textrm{ for }x_1<-\tau_j \}
$$
be a sequence of Mumford-Shah minimizers, 
where $w_j$, $\epsilon_j$, $\delta_j$, $f_j$ and $g_j$ are introduced in the Remark \ref{Rem:proj-fj}.

Assume furthermore that $\Reg(\nabla \P(u^j))\to \Reg(\nabla v^0)$ strongly in $L^2(B_1(0))$ where 
$(v^0,f_0)$, for some (unique) $f_0$, is a solution to (\ref{theLinearSyst}) that satisfies (\ref{snurrfotolin})-(\ref{fprimeestlin}) of 
Proposition \ref{AnalysisOfLinear} and that $\Reg(\nabla w^j)\to \Reg(\nabla w^0)$ strongly in $L^2(B_1)$. Furthermore we assume that there exists a 
function $g_0$ such that $g_j\rightharpoonup g_0$ locally in $W^{1,2}_{loc}(-1,0)$,
and $(w^0,g_0)$ satisfies 
(\ref{snurrfotolin})-(\ref{fprimeestlin}). 

Then there exists a constant $C$ such that $\delta_j\le C\epsilon_j$.
\end{prop}
\textsl{Proof:} The proof follows almost line for line the proof of Proposition \ref{Linearization} with slight and obvious changes.
We will include most details partly for completeness but also to convince the reader that to control the higher order asymptotics in the linearization
is done in exactly the same way as we control the second order asymptotics.

The idea of the proof is to assume that $\epsilon_j/\delta_j\to 0$ and then show that $w^j$ converges to a solution to 
(\ref{theLinearSyst}) which would contradict that $\P(u^j)$ is the projection of $u^j-\Pi(u^j,1)$ on the space of solutions to (\ref{theLinearSyst}).

For simplicity of notation we will denote $ v^j:=\epsilon_j^{-1}\P(u^j)$.

We begin by doing a domain variation (see (\ref{DomVar})) of the Mumford-Shah energy, with $\eta(x)= \psi(x) e_2$ with
$\psi\in C^{\infty}_c(B_1(0)\setminus B_{\mu_0}(0))$ and $D_2 \psi(x)=0$ close to $\Gamma_{u^j}$, and derive that
$$
0=\int_{B_1(0)\setminus \Gamma_{u^j}} \bigg[ \left|\nabla (\Pi(u^j)+\epsilon_j(v^j +\delta_j w^j)) \right|^2 \frac{\partial \psi}{\partial x_2}-
$$
$$
-2(\nabla (\Pi(u^j)+\epsilon_j(v^j +\delta_j w^j))\cdot e_2) (\nabla (\Pi(u^j)+\epsilon_j(v^j +\delta_j w^j))\cdot \nabla \psi)\bigg]+
$$
$$
+\int_{-1}^0 \frac{\epsilon_j (f'_j(x_1)+\delta_j g_j(x_1))}{\sqrt{1+\epsilon_j^2 |f'(x_1)|^2}}\frac{\partial \psi(x)}{\partial x_1}=
$$
\begin{equation}\label{Nzerothorder}
=\left[ \int_{B_1(0)\setminus \Gamma_{u^j}} \left(\left|\nabla \Pi\right|^2 \frac{\partial \psi}{\partial x_2}-2(\nabla \Pi\cdot e_2)(\nabla \Pi\cdot \nabla \psi) \right)
\right]+
\end{equation}
\begin{equation}\label{Nfirstorder1}
 +\epsilon_j\int_{B_1\setminus \Gamma_{u^j}}\Big[ 2\nabla \Pi\cdot \nabla (v^j +\delta_j w^j)\frac{\partial \psi}{\partial x_2}-2(\nabla \Pi\cdot e_2)(\nabla (v^j +\delta_j w^j)\cdot \nabla \psi)-
\end{equation}
 \begin{equation}\label{Nfirstorder3}
 -2(\nabla (v^j +\delta_j w^j)\cdot e_2)(\nabla \Pi\cdot \nabla \psi)\Big]+
\end{equation}
\begin{equation}\label{Nfirstorder2}
 +\epsilon_j\int_{-1}^0 \frac{f_j'(x_1)+\delta_j g'_j(x_1)}{\sqrt{1+\epsilon_j^2 |f_j'(x_1)|^2}}\frac{\partial \psi(x)}{\partial x_1}+
\end{equation}
\begin{equation}\label{Nsecondorder}
 +\epsilon_j^2\int_{B_1(0)\setminus \Gamma_{u^j}}\left(|\nabla (v^j +\delta_j w^j)|^2\frac{\partial \psi}{\partial x_2}-2(\nabla (v^j +\delta_j w^j)\cdot e_2)(\nabla (v^j +\delta_j w^j)\cdot \nabla \psi) \right).
\end{equation}

If we make an integration by parts in (\ref{Nzerothorder}) we arrive at
\begin{equation}\label{Nzerothorderafterintbyparts}
\int_{\Gamma_{u^j}^\pm}\left|\nabla \Pi\right|^2\psi (\nu^\pm \cdot e_2)-2(\nabla \Pi\cdot e_2)(\nabla \Pi\cdot \nu^\pm) \psi.
\end{equation}

Without loss of generality we can assume 
\begin{equation}\label{NexpressionNablaPi}
\nabla \Pi(u^j)=\frac{\lambda(j)}{2}\sqrt{\frac{2}{\pi}}\frac{1}{r^{1/2}}\big(-\sin(\phi/2),\cos(\phi/2) \big),
\end{equation}
where we, for simplicity of notation write $\lambda(j)=\lambda(u^j)\in \R$ for the constant in Remark
\ref{Rem:proj-fj}.
In particular, since the value of $\phi$ differs by $2\pi$ on $\Gamma^\pm$, we can conclude that
$\nabla \Pi(u^j)\lfloor_{\Gamma^+}=-\nabla \Pi(u^j)\lfloor_{\Gamma^-}$. Since also $\nu^+=-\nu^-$ it follows that
the value of (\ref{Nzerothorderafterintbyparts}), and therefore (\ref{Nzerothorder}), is identically zero.

Let us continue to estimate the integral in (\ref{Nsecondorder}). Lemma \ref{FirstBadC1}, together with the normalization 
$\|\nabla \Reg(v^j+\delta_j w^j)\|_{L^2(B_1)}\le 1+\delta_j\le 2$, implies that we may estimate the integral in (\ref{Nsecondorder}) by 
\begin{equation}\label{NTheEstOf32}
C\epsilon_j\sigma(\epsilon_j)\|\nabla \psi^2\|_{L^2},
\end{equation}
where $\sigma(\epsilon_j)\to 0$ is the modulus of continuity of Lemma \ref{FirstBadC1}.

This means that the terms in
(\ref{Nfirstorder1})-(\ref{Nfirstorder2}) must tend to zero as $j\to \infty$. We can thus conclude that
\begin{equation}\label{Npreceptions}
 \int_{B_1\setminus \Gamma_{u^j}}\Big[ 2\nabla \Pi\cdot \nabla (v^j +\delta_j w^j)-2(\nabla \Pi\cdot e_2)(\nabla (v^j +\delta_j w^j)\cdot \nabla \psi)-
\end{equation}
\begin{equation}\label{Npreceptions2}
 2(\nabla (v^j +\delta_j w^j)\cdot e_2)(\nabla \Pi\cdot \nabla \psi)\Big]+
\end{equation}
\begin{equation}\label{Ndecieve}
 +\int_{-1}^0 \frac{f_j'(x_1)+\delta_j g_j'(x_1)}{\sqrt{1+\epsilon_j^2 |f_j'(x_1)|^2}}\frac{\partial \psi(x)}{\partial x_1}=
\end{equation}
\begin{equation}\label{NWritingOut1}
=\int_{\Gamma_{u^j}^\pm}\Big[2(\nabla \Pi\cdot \nabla (v^j +\delta_j w^j))(\nu^\pm\cdot e_2)-2(\nabla \Pi\cdot e_2)(\nabla (v^j +\delta_j w^j)\cdot \nu^\pm)-
\end{equation}
\begin{equation}\label{NWritingOut2}
-2(\nabla (v^j +\delta_j w^j)\cdot e_2)(\nabla \Pi\cdot \nu^\pm)\Big]\psi+
\end{equation}
\begin{equation}\label{Nsympathy}
+\int_{-1}^0 \frac{f_j'(x_1)+\delta_j g_j'(x_1)}{\sqrt{1+\epsilon_j^2 |f'(x_1)|^2}}\frac{\partial \psi(x)}{\partial x_1}=o(1)\|\nabla \psi^j\|_{L^2}
\end{equation}

Remember that, by construction, $(v^j, f_j)$ solves the linearized system (\ref{theLinearSyst}). First of 
all, this implies, by the series expansion of $(v^j,f_j)$ in Proposition \ref{AnalysisOfLinear}, that $|f_j(x_1)|, |f_j'(x_1)|<C$ for $x_1\in (-b,-a)$. We may therefore 
use (\ref{NexpressionNablaPi}) to deduce that $\sin(\phi/2)=\pm 1+O(\epsilon_j)$ and $\cos(\phi/2)=O(\epsilon_j)$ on
$\Gamma_u^\pm\cap B_b \setminus B_a$ for any $a,b\in (0,1)$. Using that $(v^j,f_j)$ solves the linearized system in equation (\ref{Nsympathy}) we may also 
conclude that
\begin{equation}\label{NfirstequationPart2}
O(\epsilon_j)=\delta_j \int_0^1\left( \sqrt{\frac{2}{\pi}\frac{1}{r}}\left(\frac{\partial w^j}{\partial x_1}\Big\lfloor_{\Gamma_u^+}+
\frac{\partial w^j }{\partial x_1}\Big\lfloor_{\Gamma_u^-}\right)\psi+\frac{g'_j}{\sqrt{1+\epsilon_j^2 |f_j|^2}}\frac{\partial \psi(x)}{\partial x_1}\right),
\end{equation}
where we used Corollary \ref{CPiLem} to estimate $|\lambda(j)-1|\le C\epsilon_j$.

We will now use the assumption that will lead to a contradiction: 
\begin{equation}\label{NToC}
\frac{\epsilon_j}{\delta_j}\to 0.
\end{equation}
From (\ref{NfirstequationPart2}) and the assumption (\ref{NToC}) we may conclude that 
\begin{equation}\label{NBdrywg1}
g_0''(x_1)=\sqrt{\frac{2}{\pi}\frac{1}{r}}\left(\frac{\partial w^0(x_1,0^+)}{\partial x_1}+\frac{\partial w^0(x_1,0^-)}{\partial x_1}\right).
\end{equation}
Our aim is to show that $(w^0,g_0)$ also satisfies the other boundary condition in (\ref{theLinearSyst}) this will lead to a contradiction to
$v^j$ being the projection on the solutions to the linearized system. But in order to do this we need to analyze the second boundary condition.

\vspace{3mm}

To derive the second boundary condition for $(w^0,g_0)$ we use that on $\Gamma_{u^j}^\pm$
\begin{equation}\label{NBoringAsHell}
0=\nabla u^j\cdot \nu= \nabla \left( \Pi(u^j)+\epsilon_j (v^j+\delta_j w^j)\right)\cdot (\epsilon_j (f_j'+\delta_j g_j'), -1).
\end{equation}
That is
\begin{equation}\label{NsomeBSequation}
\epsilon_j \frac{\partial (v^j+\delta_j w^j)}{\partial x_1} (f_j'+\delta_j g_j')-\frac{\partial (v^j+\delta_j w^j)}{\partial x_2}=
\frac{1}{\epsilon_j}\left(\nabla \Pi \cdot (-\epsilon_j f_j',1)\right).
\end{equation}
If we use (\ref{NexpressionNablaPi}) in (\ref{NsomeBSequation}) we see that on $\Gamma_{u^j}^+$
$$
\phi=\pi-\frac{\epsilon_j (f_j+\delta_j g_j)}{r}+O\left(\left(\frac{\epsilon_j (f_j+\delta_j g_j)}{r}\right)^3\right)
$$
and therefore
$$
\frac{1}{\epsilon_j}\left(\nabla \Pi \cdot (-\epsilon_j (f_j'+\delta_j g_j'),1)\right)=
$$
\begin{equation}\label{NDontever}
=\frac{\lambda(j)}{2\epsilon_j}\sqrt{\frac{2}{\pi}\frac{1}{r}}\Big(\epsilon_j\sin\left(\frac{\pi-\epsilon_j(f_j+\delta_jg_j)/r}{2}\right)(f_j'+\delta_j g_j')+
\end{equation}
$$
\cos\left(\frac{\pi-\epsilon_j (f_j+\delta_j g_j)/r}{2}\Big) \right)=
$$
$$
=\frac{1}{2}\sqrt{\frac{2}{\pi}\frac{1}{r}}\left((f_j'(x_1)+\delta_jg_j')+\frac{f_j(x_1)+\delta_j g_j}{2r} \right)+O(\epsilon_j),
$$
where we again used Corollary \ref{CPiLem} to estimate $|\lambda(j)-1|\le C\epsilon_j$.

Equations (\ref{NBoringAsHell}), (\ref{NsomeBSequation}) and (\ref{NDontever}) 
together imply that, in the weak sense,
\begin{equation}\label{NNonsense}
-\epsilon_j f_j'\frac{\partial (v^j+\delta_j w^j)}{\partial x_1}-\frac{\partial (v^j+\delta_j w^j)(x_1,0^-)}{\partial x_2}=
\end{equation}
$$
=\frac{1}{2}\sqrt{\frac{2}{\pi}\frac{1}{r}}\left(f_j'+\delta_j g_j'+\frac{f_j+\delta_j g_j}{2r} \right)+O(\epsilon_j).
$$
But since $(v^j,f_j)$ solves to linearized system (and also that $\|\Gamma_{u^j}\|_{C^{1,\alpha}}\le C\epsilon_j$) we can conclude that 
\begin{equation}\label{NNonsense78}
-\frac{\partial w^j(x_1,0^-)}{\partial x_2}=\frac{1}{2}\sqrt{\frac{2}{\pi}\frac{1}{r}}\left(g_j'+\frac{g_j}{2r} \right)+O(\epsilon_j/\delta_j).
\end{equation}

Passing to the limit in (\ref{NNonsense78}) and using (\ref{NToC}) we may conclude that
\begin{equation}\label{NSonsFindDevils2}
-\frac{\partial w^0(x_1,0^+)}{\partial x_2}=\frac{1}{2}\sqrt{\frac{2}{\pi}\frac{1}{r}}\left(g_0'(x_1)+\frac{g_0(x_1)}{2r} \right).
\end{equation}

We may argue similarly on $\Gamma_u^-$, where $\phi=-\pi-\frac{\epsilon_j f_j}{r}+O\left(\left(\frac{\epsilon_j f_j}{r}\right)^3\right)$,
we can conclude, after passing to the limit, that
\begin{equation}\label{NSonsFindDevils2Minus}
-\frac{\partial w^0(x_1,0^-)}{\partial x_2}=-\frac{1}{2}\sqrt{\frac{2}{\pi}\frac{1}{r}}\left(g_0'(x_1)+\frac{g_0(x_1)}{2r} \right).
\end{equation}

We are now ready to prove the proposition. If no constant $C$ exists such that $\delta_j<C\epsilon_j$ then (\ref{NBdrywg1}), (\ref{NSonsFindDevils2}) 
and (\ref{NSonsFindDevils2Minus}) holds. This implies that $(w^0,g_0)$ solves the linearized system. Also, by strong convergence 
$\Reg(\nabla w^j)\to \Reg(\nabla w^0)$, and $w^0$ is not identically zero. However, since $\epsilon_j v^j=\P(u^j)$ it follows that $w^0$ is 
orthogonal to all solutions to the linearized system. This is clearly a contradiction. It follows that $\delta_j<C\epsilon_j$ for some constant $C$. \qed

Next we show that the two dominating terms in the asymptotic expansion of a minimizer $u^j$
are $\Pi(u^j)$ and $\P(u^j)$. This should be very intuitively clear.

\begin{lem}\label{NdeltagIsZerothOrder}Let   
$$
u^j=\Pi(u^j)+\P(u^j)+\epsilon_j\delta_j w^j
$$ 
and
$$
\Gamma_{u^j}=\{(x_1,\epsilon_j(f_j+\delta_j g_j)),\;\textrm{ for }x_1<-\tau_j \}
$$
be a sequence of Mumford-Shah minimizers, 
where $w_j$, $\epsilon_j$, $\delta_j$, $f_j$ and $g_j$ are introduced in the Remark \ref{Rem:proj-fj}.

Assume furthermore that $\epsilon_j\to 0$,
$$
\frac{\Reg(\nabla \P(u^j))}{\epsilon_j}\to \Reg(\nabla v^0) \textrm{ strongly in }L^2(B_1),
$$
where $(v^0,f_0)$, for some (unique) $f_0$, is a solution to the linearized system (\ref{theLinearSyst}) and satisfies (\ref{snurrfotolin})-(\ref{fprimeestlin})
of Proposition \ref{AnalysisOfLinear}.

Then for any $0<a<b<1$ it follows that $\delta_j g_j \to 0$ in $C^{1,\alpha}((-b,-a))$ and that 
$\|\delta_j w^j\|_{C^{1,\alpha}((B_b\setminus B_a)\setminus \Gamma_{u^j})}\to 0$.
\end{lem}
\textsl{Proof:} Since, by Proposition \ref{Linearization}, Proposition \ref{SatisfiesRightEstProp} and Corollary \ref{FconvC1alpha}, $(f_j +\delta_j g_j)$ 
converges in $C^{1,\alpha}$ to a solution to linear system and the limit of $f_j$ converges to that solution, too, it follows that $\delta_jg_j\to 0$ 
and similarly for $\delta_j w^j$.\qed

We will now prove a convergence result that refines Lemma \ref{WhatShouldWeCallThis}. The proof of the next lemma is very similar to the proof of its
sister Lemma \ref{WhatShouldWeCallThis}.

\begin{lem}\label{NLocalConvW}
Let $(u^j,\Gamma_j)$ be a sequence of minimizers to the Mumford-Shah problem and assume that 
$u^j=\Pi(u^j)+\P(u^j)+\delta_j \epsilon_j w^j$ where $\|\nabla \P(u^j)\|_{L^2(B_1\setminus \Gamma_{j})}=\epsilon_j\to 0$ and 
$\|\nabla w^j\|_{L^2(B_1\setminus \Gamma_j)}\le 1$. Assume furthermore that 
$$
\frac{\Reg(\nabla \P(u^j))}{\epsilon_j}\to \Reg(\nabla v^0) \textrm{ strongly in }L^2(B_1),
$$
where $(v^0,f_0)$, for some (unique) $f_0$, is a solution to the linearized system (\ref{theLinearSyst}) and satisfies (\ref{snurrfotolin})-(\ref{fprimeestlin})
of Proposition \ref{AnalysisOfLinear}.

Then for $0<a<b<1$ if there exists a constant $C_0>0$ such that if $\delta_j>C_0\epsilon_j$ then
$$
\Reg(\nabla w^j)\to \Reg(\nabla w^0) \textrm{ strongly in }L^2(B_b\setminus B_a).
$$
Furthermore if $f_j$ is defined so that $(\P(u^j),\epsilon_j f_j)$ solves (\ref{theLinearSyst}), 
and $g_j$ is defined so that the free discontinuity 
$$
\Gamma_j=\{(x_1,\epsilon_j (f_j+\delta_j g_j));\; x_1< -a/2 )\}
$$
then $g_j\to g_0$ for in $C(-b,-a)$ and weakly in $W^{1,2}(-b,-a)$ for some function $g_0$.
\end{lem}
\textsl{Proof:} We argue as in (\ref{Shitstar1}) and use that $u^j$ satisfies the Neumann condition on $\Gamma_{u^j}$ which implies that 
\begin{equation}\label{NExprForwjOnGamma}
\nabla w^j\cdot \nu_j^\pm =-\frac{1}{\epsilon_j\delta_j}\nabla \Pi(u^j)\cdot \nu_j^\pm-\frac{1}{\delta_j}\nabla v^j\cdot \nu_j^\pm.
\end{equation}
Continuing, as in (\ref{Shitstar3}) and using (\ref{NexpressionNablaPi}), we may calculate that
$$
\nabla \Pi(u^j)\cdot \nu_j^\pm =
$$
\begin{equation}\label{NSecondOrderExpOfPi}
=\frac{1}{2}\sqrt{\frac{2}{\pi r}}\left( 
\epsilon_j\sin\left( \frac{\pi-\epsilon_j(f_j+\delta_jg_j)/r}{2}\right)+
\cos\left(\frac{\pi-\epsilon_j (f_j+\delta_j g_j)/r}{2\sqrt{1+\epsilon_j^2(f_j'+\delta_jg_j')^2}}\right)\right)=
\end{equation}
$$
=\frac{1}{2}\sqrt{\frac{2}{\pi r}}\left( -\epsilon_j (f'_j+\delta_j g'_j)-\frac{f_j+\delta_j g_j}{2r \sqrt{1+\epsilon_j^2(f_j'+\delta_jg_j')^2}}\right)+
$$
$$
+O\left(\epsilon_j^3 (f_j+\delta_jg_j)^3+\epsilon_j^3 (f_j'+\delta_jg_j')^3\right),
$$
where we used trigonometric identities and first order Taylor expansions of trigonometric functions.

Next we use that $(v^j,f_j)$ are solutions to the linearized system (\ref{theLinearSyst}) and therefore
$$
\nabla v^j \cdot \nu_j^\pm =-\frac{1}{2}\sqrt{\frac{2}{\pi r}}\left( f'_j+\frac{1}{2r}f_j\right),
$$
on $\Gamma_0$ (clearly we need to place the branch cut of $(v^j,f_j)$ differently for this equality to hold but that does not affect the argument).
But since $v^j,f_j,g_j \in C^{1,\alpha}$ and $\Gamma_{u^j}$ is given by a graph of $\epsilon_j(f_j+\delta_j g_j)$
which has $C^{1,\alpha}-$norm bounded by $C\epsilon$, here we use that $f_j$ is smooth and $\delta_jg_j$ negligible by Lemma \ref{NdeltagIsZerothOrder},
we may conclude that
\begin{equation}\label{NSomeOtherCalc}
\nabla v^j \cdot \nu_j^\pm =-\frac{1}{2}\sqrt{\frac{2}{\pi r}}\left( f'_j+\frac{1}{2r}f_j\right)+O(\epsilon_j)
\end{equation}
on $(\Gamma_{u^j}\cap B_b(0))\setminus B_a(0) $.

If we insert (\ref{NSecondOrderExpOfPi}) and (\ref{NSomeOtherCalc}) into (\ref{NExprForwjOnGamma}), and also use Lemma \ref{NdeltagIsZerothOrder}
to disregard $\delta_j g_j'$ terms, it will follow that
$$
\nabla w^j\cdot \nu_j^\pm =\frac{1}{2}\sqrt{\frac{2}{\pi r}}\left( -g_j'-\frac{g_j}{2r}\right)+O(\epsilon_j/\delta_j).
$$
Since we assume that $\delta_j>C_0 \epsilon_j$ the term $O(\epsilon_j/\delta_j)$ is bounded by a constant.

It follows that $w^j$ satisfies 
\begin{equation}\label{NOneStarApr6}
\begin{array}{ll}
 \Delta w^j =0 & \textrm{ in }B_1\setminus \Gamma_{u^j} \\
 \frac{\partial w^j}{\partial \nu^j}=\frac{1}{2}\sqrt{\frac{2}{\pi r}}\left( -g_j'-\frac{g_j}{2r}\right)+O(\epsilon_j/\delta_j) & \textrm{ on }\Gamma_{u^j}.
\end{array}
\end{equation}

In order to show that (\ref{NOneStarApr6}) implies strong convergence of $\Reg(\nabla w^j)$ it is enough to show that $\|g_j\|_{W^{1,2}(-b,-a)}\le C$. 
The proof of that is very similar to the proof that $f_j\in W^{1,2}$ in Lemma \ref{WhatShouldWeCallThis}. 

We choose $\psi^j$ as the solution to 
\begin{equation}
 \begin{array}{ll}
  \Delta \psi^j=0 & \textrm{ in }(B_b\setminus B_a)\setminus \Gamma_{u^j} \\
  \psi^j= 0 & \textrm{ on } \partial B_b \cup \partial B_a \\
  \psi^j = g_j -\frac{g_j(b)-g_j(a)}{b-a}(-x_1-a)+g_j(a) & \textrm{ on }\Gamma_{u^j}.
 \end{array}
\end{equation}
Then using $\psi^j$ as a test function in (\ref{Npreceptions})-(\ref{Nsympathy}) we may deduce that 
\begin{equation}\label{NFuckingWoi0}
\int_{B_1\setminus\Gamma_{u^j}}\Big[ 2\nabla \Pi(u^j)\cdot \nabla w^j\frac{\partial \psi^j}{\partial x_2}-2\left(\nabla \Pi(u^j)\cdot e_2 \right)\nabla w^j\cdot \nabla \psi^j-
\end{equation}
\begin{equation}\label{NFuckingWoi}
-2(\nabla w^j\cdot e_2)(\nabla \Pi(u^j)\cdot \psi^j)\Big]+\int_0^1 \frac{g_j'}{\sqrt{1+\epsilon_j^2(f_j'+\delta_jg_j')^2}}\frac{\partial \psi^j}{\partial x_1}=
\end{equation}
\begin{equation}\label{NFuckingWoi2}
=\frac{-1}{\delta_j}\bigg[\int_{B_1\setminus \Gamma_{u^j}}\Big(2\nabla \Pi(u^j)\cdot \nabla v^j-2(\nabla \Pi(u^j)\cdot \nabla v^j)\frac{\partial \psi^j}{\partial x_2}-
\end{equation}
\begin{equation}\label{NFuckingWoi3}
-2(\nabla \Pi(u^j)\cdot e_2)\nabla v^j\cdot \nabla \psi^j-2(\nabla v^j\cdot e_2)(\nabla \Pi(u^j)\cdot \nabla \psi^j)\Big)\Big]+
\end{equation}
\begin{equation}\label{NFuckingWoi4}
+\int_0^1 \frac{f_j'}{\sqrt{1+\epsilon_j^2(f_j'+\delta_jg_j')^2}}\frac{\partial \psi^j}{\partial x_1}+J(j,\psi^j),
\end{equation}
where we have written $J(j,\psi^j)$ for the $o(1)\|\nabla \psi^j\|_{L^2}/\delta_j$ term in (\ref{Nsympathy}). 
For this lemma we need a better estimate for that term so we will write it out more explicitly.
The term $o(1)\|\nabla \psi^j\|_{L^2}/\delta_j$ in (\ref{Nsympathy}) came from the integral in (\ref{Nsecondorder}), which we have divided by $\epsilon_j\delta_j$ 
in (\ref{NFuckingWoi}). This leads to 
\begin{equation}\label{NJjPsij}
J(j,\psi^j)=\frac{\epsilon_j}{\delta_j}\int_{B_1(0)\setminus \Gamma_{u^j}}\Big(|\nabla (v^j +\delta_j w^j)|^2\frac{\partial \psi}{\partial x_2}
\end{equation}
$$
-2(\nabla (v^j +\delta_j w^j)\cdot e_2)(\nabla (v^j +\delta_j w^j)\cdot \nabla \psi) \Big).
$$

So far we have only estimated $J(j,\psi^j)$, or equivalently the integral in (\ref{Nsecondorder}) which appears in (\ref{NJjPsij}), 
by $\frac{\sigma(\epsilon_j)\|\nabla \psi^j\|_{L^2}}{\delta_j}$ (which is the estimate in (\ref{NTheEstOf32})
divided by $\epsilon_j\delta_j$) and we need to provide a stronger estimate of $J(j,\psi^j)$ in order to get the desired bound on $g_j'$. 

In order to estimate the integral in (\ref{NJjPsij}) we will use that $v^j$ and $\nabla v^j$ are bounded, 
since $v^j$ solves the linearized system, and that $\delta_j \Reg(\nabla w^j)\to 0$ in 
$B_b\setminus B_a$ (by Lemma \ref{NdeltagIsZerothOrder}):
$$
\left|J(j,\psi^j)\right|= 
$$
$$
\frac{\epsilon_j}{\delta_j}\int_{B_1(0)\setminus \Gamma_{u^j}}\left||\nabla (v^j +\delta_j w^j)|^2\frac{\partial \psi}{\partial x_2}-2(\nabla (v^j +\delta_j w^j)\cdot e_2)(\nabla (v^j +\delta_j w^j)\cdot \nabla \psi) \right|
$$
\begin{equation}\label{Aristoteles}
\le C\frac{\epsilon_j}{\delta_j} \|\nabla(v^j+\delta_j w^j)\|_{L^\infty}^2\int_{B_1(0)\setminus \Gamma_{u^j}}|\nabla \psi^j|\le C\frac{\epsilon_j}{\delta_j}\|\nabla \psi^j\|_{L^2}.
\end{equation}

Using that $(v^j,f_j)$ solves the linearized system together with the estimate of $J(j,\psi^j)$ in (\ref{Aristoteles}) we may conclude that the right side in 
(\ref{NFuckingWoi0})-(\ref{NFuckingWoi4}) is of order $O(\epsilon_j/\delta_j)\|\nabla \psi^j\|_{L^2(B_1\setminus \Gamma_{u^j})}$.
We have therefore shown that 
\begin{equation}\label{NFuckingWoi5}
\int_{B_1\setminus\Gamma_{u^j}}\Big[ 2\nabla \Pi(u^j)\cdot \nabla w^j\frac{\partial \psi^j}{\partial x_2}-2\left(\nabla \Pi(u^j)\cdot e_2 \right)\nabla w^j\cdot \nabla \psi^j-
\end{equation}
\begin{equation}
 -2(\nabla w^j\cdot e_2)(\nabla \Pi(u^j)\cdot \psi^j)\Big]+\int_0^1 \frac{g_j'}{\sqrt{1+\epsilon_j^2(f_j'+\delta_jg_j')^2}}\frac{\partial \psi^j}{\partial x_1}=
\end{equation}
\begin{equation}\label{NFuckingWoi6}
=O(\epsilon_j/\delta_j)\|\nabla \psi^j\|_{L^2(B_1\setminus \Gamma_{u^j})}.
\end{equation}
Since, by assumption, $\delta_j\le C\epsilon_j$ it follows that (\ref{NFuckingWoi5})-(\ref{NFuckingWoi6}) are 
exactly the same estimates 
as  the estimates (\ref{preceptions2})-(\ref{decieve2}) we derived in Lemma \ref{WhatShouldWeCallThis}.

The rest of the proof follows from (\ref{NFuckingWoi5})-(\ref{NFuckingWoi6}) by the same calculations as in (\ref{EstOfF})-(\ref{Shitstar2}). \qed

\begin{prop}\label{PropContOfRest}
 Let $(u^j,\Gamma_{u^j})$ be as in Proposition \ref{Linearization}. Then there exists a constant $C$ such that for every $r\in (0,3/4)$
 \begin{equation}\label{ControllOnAllScales}
 \left\| \nabla \left(\frac{u^j(rx)}{\sqrt{r}}-\Pi\left(\frac{u^j(rx)}{\sqrt{r}} \right)-\P\left(\frac{u^j(rx)}{\sqrt{r}} \right)\right)\right\|_{L^2(B_1\setminus \Gamma_{u^j})}\le C_0r^{-2\kappa}\epsilon_j^2.
 \end{equation}
It follows in particular that there exist a subsequence in $j$ such that 
\begin{equation}\label{NStrongConvWj}
\frac{1}{\epsilon_j^2}\Reg\left( \nabla \left[\frac{u^j(rx)}{\sqrt{r}}-\Pi\left(\frac{u^j(rx)}{\sqrt{r}} \right)-\P\left(\frac{u^j(rx)}{\sqrt{r}} \right)\right]\right)\to \nabla \Reg(w^0)
\end{equation}
strongly in $L^2(B_{3/4})$ for some $w^0$.
\end{prop}
\textsl{Proof:} The proof is very similar to the proof of Proposition \ref{SatisfiesRightEstProp} and 
Lemma \ref{EstIfThereExistsBound} and \ref{CloseForSmallerr}. We prove the proposition in three claims.

\vspace{2mm}

{\bf Claim 1:} {\sl (Similar to Lemma \ref{CloseForSmallerr}.) The following limit holds}
\begin{equation}\label{NHereWeNeedSupremum}
\lim_{j\to \infty} \sup_{r\in (0,1]}\left\|\Reg\left(\nabla \left(\frac{u^j(rx)}{\sqrt{r}}-\Pi(u^j,r)-\P(u^j,r)\right)\right)\right\|_{L^2(B_1)}=0
\end{equation}
{\sl and for each $j$ the supremum is achieved for some $r_j$.}

\vspace{1mm}

{\sl Proof of Claim 1:} Clearly, since $\P$ is defined as a projection operator, 
$$
\sup_{r\in (0,1]}\left\|\Reg\left(\nabla \left(\frac{u^j(rx)}{\sqrt{r}}-\Pi(u^j,r)-\P(u^j,r)\right)\right)\right\|_{L^2(B_1)}\le
$$
$$
\le \sup_{r\in (0,1]}\left\|\Reg\left(\nabla \left(\frac{u^j(rx)}{\sqrt{r}}-\Pi(u^j,r)\right)\right)\right\|_{L^2(B_1)}.
$$
But the right side in the previous expression tends to zero as $j\to \infty$ by Lemma \ref{CloseForSmallerr}. This proves claim 1.

\vspace{2mm}

{\bf Claim 2:} {\sl (Similar to Lemma \ref{EstIfThereExistsBound}.)} {\sl Under the assumption that (\ref{ControllOnAllScales}) holds for every $r\in (0,3/4]$ the strong convergence 
(\ref{NStrongConvWj}) holds.}

\vspace{1mm}

{\sl Proof of Claim 2:} The proof of claim 2 is exactly the same as the proof of Lemma \ref{EstIfThereExistsBound}. The only difference is that 
one has to change $\Pi(u^j)$ to $\Pi(u^j)-\P(u^j)$ in all calculations and refer to Lemma \ref{NLocalConvW} instead of Lemma \ref{WhatShouldWeCallThis}.

\vspace{2mm}

{\bf Claim 3:} {\sl (Similar to Proposition \ref{SatisfiesRightEstProp}.)} {\sl The estimate (\ref{ControllOnAllScales}) holds for every $r\in (0,3/4]$.}

\vspace{1mm}

{\sl Proof of Claim 3:} We will use similar notation as Proposition \ref{NLinearization}. 
To be specific, we will define $u_r^j$, $\epsilon_j(r)$, $\delta_j(r)$ and $w_r^j$ by
$$
u_r^j(x)=\frac{u^j(rx)}{\sqrt{r}}=\Pi(u_r^j,1)+\P(u_r^j,1)+\epsilon_j(r)\delta_j(r)w_r^j,
$$
where 
\begin{equation}\label{DefOfepsjr}
\epsilon_j(r)=\|\Reg(\nabla \P(u_r^j,1))\|_{L^2(B_1(0))}
\end{equation}
and $\delta_j(r)$ is defined by the condition that 
$$
\|\Reg(\nabla w_r^j)\|_{L^2(B_1(0))}=1.
$$

We can estimate
\begin{equation}\label{Eq1Claim3}
\left\| \nabla \left(\frac{u^j(rx)}{\sqrt{r}}-\Pi\left(\frac{u^j(rx)}{\sqrt{r}} \right)-\P\left(\frac{u^j(rx)}{\sqrt{r}} \right)\right)\right\|_{L^2(B_1\setminus \Gamma_{u_r^j})}=
\end{equation}
$$
=\left\|\Reg(\epsilon_j(r)\delta_j(r) \nabla w^j_r)\right\|_{L^2}=
\epsilon_j(r) \delta_j(r)\le C\epsilon_j(r)\left\|\Reg(\nabla \P(u^j_r))\right\|_{L^2}=C \epsilon_j(r)^2,
$$
where we used Proposition \ref{NLinearization} in the first equality, that 
$\|\Reg(\nabla w^j_r)\|_{L^2}=1$ by 
construction in the first equality and the definition of $\epsilon_j(r)$ from 
(\ref{DefOfepsjr}) in the last equality.

From Proposition \ref{SatisfiesRightEstProp} we know that for some constant $C$
\begin{equation}\label{Eq2Claim3}
\left\|\Reg\left( \nabla \left( u^j_r-\Pi\left(u^j_r\right)\right)\right)\right\|_{L^2(B_1(0)}\le  
\end{equation}
$$
\le Cr^{-\kappa}\left\|\Reg\left(\nabla \left(u^j-\Pi\left(u^j\right)\right)\right)\right\|_{L^2(B_1(0)}. 
$$
Next we remark that 
\begin{equation}\label{Eq3Claim3}
 \left\|\Reg\left( \nabla \left(u^j_r-\Pi\left(u^j_r\right)\right)\right)\right\|_{L^2(B_1(0)}=
\end{equation}
$$
= \left\|\Reg(\nabla ( \P(u^j_r)+ \epsilon_j(r) \delta_j(r) w^j_r))\right\|_{L^2}\le \epsilon_j(r) +\epsilon_j(r)\delta_j(r) \le 2\epsilon_j(r),
$$
where we used that $|\delta_j(r)|< < 1$ in the last inequality. And similarly
\begin{equation}\label{Eq4Claim3}
\left\|\Reg\left( \nabla \left(u^j_r-\Pi\left(u^j_r\right)\right)\right)\right\|_{L^2(B_1(0)}
\ge \frac{\epsilon_j(r)}{2}.
\end{equation}

Using  (\ref{Eq3Claim3}) and (\ref{Eq4Claim3}) in (\ref{Eq2Claim3}) we may deduce that 
\begin{equation}\label{Eq5Claim3}
\epsilon_j(r)\le 4C r^{-\kappa} \epsilon_j(1)=4C r^{-\kappa} \epsilon_j.
\end{equation}

Using (\ref{Eq5Claim3}) in (\ref{Eq1Claim3}) we can conclude that 
$$
\left\| \nabla \left(\frac{u^j(rx)}{\sqrt{r}}-\Pi\left(\frac{u^j(rx)}{\sqrt{r}} \right)
-\P\left(\frac{u^j(rx)}{\sqrt{r}} \right)\right)\right\|_{L^2(B_1\setminus \Gamma_{u^j})}\le C_0r^{-2\kappa}\epsilon_j^2
$$
for some constant $C_0$. This proves Claim 3. 

\vspace{2mm}

The proposition follows from Claim 2 and Claim 3. \qed

We are now ready to prove the main regularity theorem. The proof is very similar to the proof of Theorem \ref{C1alpha}. However some 
technicalities differ wherefore we will provide some details of the proof.

\begin{thm}\label{C2alpha}
For every $\alpha< \alpha_2-3/2\approx 0.389$ there exists an $\epsilon_\alpha>0$ such that if $(u,\Gamma)$ is $\epsilon$-close to a crack-tip solution for 
some $\epsilon\le \epsilon_\alpha$ then $\Gamma_u$ is $C^{2,\alpha}$ at the crack-tip.

 This in the sense that the tangent at the crack-tip is a well defined line, which we may assume to be $\{(x_1,0);\; x_1\in \R\}$,
 and there exists a constant $C_\alpha$ such that
 $$
 \Gamma_u\subset \left\{(x_1,x_2);\; |x_2|<C_\alpha \epsilon |x_1|^{2+\alpha},\; x_1<0\right\}.
 $$
 Here $C_\alpha$ may depend on $\alpha$ but not on $\epsilon<\epsilon_\alpha$.
\end{thm}
\textsl{Proof:} We assume that $(u,\Gamma)$ is a solution that is $\epsilon-$close to a crack-tip for some very small $\epsilon$.
Then we may write, for some small constant $s$ and $j\in \N$, $u_j(x)=\frac{u(s^j x)}{\sqrt{s^j}}$. For each $j$ we will rotate the coordinate system by a small constant  
$\phi_j$ so that we may express $u_j$ in the new coordinates according to
\begin{equation}\label{ExpressionForuj}
u_j(x)=\Pi(u_j,1)+S_j(x)+C_j(x)+R_j(x)
\end{equation}
where 
$$
S_j(x)=S_j(r,\phi)=\sum_{k=2}^\infty b_k(j) r^{k-1/2}\sin((k-1/2)\phi),
$$
$$
C_j(x)=C_j(r,\phi)=\sum_{k=2}^\infty a_k(j) r^{\alpha_k}\cos(\alpha_k \phi)
$$
and 
\begin{equation}\label{StarIsTheNew}
R_j(x)=R_j(r,\phi)=a_0(j)\mathfrak{z}(r,\phi)+\hat{R}_j(r,\phi),
\end{equation}
where $\mathfrak{z}$ is defined in (\ref{DefOfW}).
Since 
$$
\P(u_j,1)=S_j(x)+C_j(x)+a_0(j)\mathfrak{z}
$$
it follows from Proposition \ref{NLinearization} that
\begin{equation}\label{InteBlind}
\|\Reg(\nabla \hat{R}_j)\|_{L^2(B_1)}\le C_R\left\|\Reg(\nabla (S_j+C_j+a_0(j)\mathfrak{z})\right\|_{L^2(B_1)}^2.
\end{equation}
We also notice that, by Proposition \ref{NoLOTinLim}, 
 there exists a modulus of continuity $\sigma(\epsilon)$ such that 
\begin{equation}\label{EstOna00}
 |a_0(0)|\le \sigma(\epsilon)\epsilon.
\end{equation}

The proof will consist in showing that $|\phi_j-\phi_{j+1}|$ will only depend on $C_j$ and $R_j$.
Since $C_j(r,\phi)$ consists of terms that decay faster than $r^{\alpha_2}$ it will follow that $r^{\alpha_2}$ determines the rotation of 
the coordinate system. Our first aim is to 
connect the functions $C_1$ and $S_1$ to the functions $C_0$ and $S_0$. By scaling invariance 
analogous estimates hold for $C_{j+1}$ and $S_{j+1}$ and $C_j$ and $S_j$ by the same proof.

\vspace{2mm}

{\bf Claim 1:} {\sl There exists constants $K_1$ and  $K_2$ such that 
if $\epsilon$ and $s$ are small enough}
\begin{equation}\label{EstFromAbovee}
 \left\|\Reg(\nabla (C_1+S_1))\right\|_{L^2(B_1(0))} \le K_1s\left\|\Reg(\nabla(S_0+C_0))\right\|_{L^2(B_1(0)}+K_2\epsilon^2.
\end{equation}

\vspace{2mm}

{\sl Proof of Claim 1:} The proof is a straightforward estimate using that 
$u_1(x)=\frac{u_0(sx)}{\sqrt{s}}$. Using this we see that  
$$
u_1(x)=\frac{u_0(sx)}{\sqrt{s}}=
$$
\begin{equation}\label{ExpForu1}
=a_0(0)\frac{\mathfrak{z}(sx)}{\sqrt{s}}+\frac{\Pi(u_0)(sx)+S_0(sx)+C_0(sx)}{\sqrt{s}}
+\frac{\hat{R}_0(sx)}{\sqrt{s}}.
\end{equation}
Furthermore, by the definition of $\mathfrak{z}$ (see (\ref{DefOfW})),
\begin{equation}\label{RescalingOffrakz}
\frac{\mathfrak{z}(sr,\phi)}{\sqrt{s}}=\mathfrak{z}(r,\phi)-\ln(s)r^{1/2}\cos(\phi/2).
\end{equation}
From the expression for $u_1$ in (\ref{ExpressionForuj}) and the definition of the 
projection operator $\P$ it follows that 
\begin{equation}\label{NeedToBeNormed}
a_0(1)\mathfrak{z}(x)+ S_1(x)+C_1(x)=\P(u_1)=\P(u_1-\Pi(u_1))=
\end{equation}
$$
=a_0(0)\mathfrak{z}(x)+
\frac{\P(S_0(sx)+C_0(sx))}{\sqrt{s}}+\P(\hat{R}_0(sx)/\sqrt{s}).
$$
Taking the norm of the gardients on both sides in (\ref{NeedToBeNormed}), and using  (\ref{EstOna00}) (and the analogous estimate for $a_0(1)$) to get rid of the $a_0$ terms, it 
follows that 
\begin{equation}\label{LeftAlone}
(1-\sigma(\epsilon))\left\|\Reg(\nabla (S_1+C_1))\right\|_{L^2(B_1)}\le
\end{equation}
$$
\le (1+\sigma(\epsilon))
\left\|\Reg\left(\nabla \frac{S_0(sx)+C_0(sx)}{\sqrt{s}}\right)\right\|_{L^2(B_1)}+
\left\|\Reg(\nabla \P(\hat{R}_0(sx)/\sqrt{s}))\right\|_{L^2(B_1)}.
$$
Using that $\P$ is a projection it follows that 
\begin{equation}\label{LeftAlone2}
\left\|\Reg(\nabla \P(\hat{R}_0(sx)/\sqrt{s}))\right\|_{L^2(B_1)}\le 
\left\|\Reg(\nabla (\hat{R}_0(sx)/\sqrt{s}))\right\|_{L^2(B_1)}\le K_2 \epsilon^2,
\end{equation}
where we also used (\ref{InteBlind}) in the last estimate. The constant $K_2$ in 
(\ref{LeftAlone2}) may depend on $s$, but we will always choose $s$ before $\epsilon$ 
and therefore the $s$ dependence will not affect the argument.

Using the expression for $S_0$ and $C_0$ we see that 
$$
\frac{S_0(sx)}{\sqrt{s}}=s\sum_{k=2}^\infty b_k(0) s^{k-2}r^{k-1/2}\sin((k-1/2)\phi)
$$
and 
$$
\frac{C_0(sx)}{\sqrt{s}}=s\sum_{k=2}^\infty a_k(0) s^{\alpha_k-3/2}r^{\alpha_k}\cos(\alpha_k\phi).
$$
Note that this implies, if $s$ is small enough, that 
\begin{equation}\label{Serious}
\left\|\Reg\left(\nabla \left(\frac{S_0(sx)+C_0(sx)}{\sqrt{s}}\right)\right)\right\|_{L^2(B_1)}
\le K_1 s \left\|\Reg(\nabla (S_0+C_0))\right\|_{L^2(B_1)}.
\end{equation}

Using (\ref{Serious}) and (\ref{LeftAlone2}) in (\ref{LeftAlone}) we can conclude that 
$$
\left\|\Reg(\nabla (S_1+C_1))\right\|_{L^2(B_1)} \le
K_1s\left\|\Reg(\nabla (S_0+C_0))\right\|_{L^2(B_1)}+K_2\epsilon^2.
$$
This proves (\ref{EstFromAbovee}).

\vspace{2mm}

{\bf Claim 2:} {\sl There is a constant $K$ such that 
\begin{equation}\label{Claim1Thm81}
|a_0(0)|\le K \epsilon^2. 
\end{equation}
It follows in particular, from (\ref{StarIsTheNew}) and (\ref{InteBlind}), 
that }
\begin{equation}\label{BetterRestEst}
\|\Reg(\nabla R_0)\|_{L^2(B_1)}\le  C\epsilon^2.
\end{equation}

{\sl For $a_0(j)$ and $\|\Reg(\nabla R_j)\|_{L^2(B_1)}$ the corresponding estimates }
$$
|a_0(j)|\le K \|\Reg(\nabla (C_j+S_j))\|_{L^2(B_1)}^2
$$
{\sl and }
$$
\|\Reg(\nabla R_j)\|_{L^2(B_1)}\le  C\|\Reg(\nabla (C_j+S_j))\|_{L^2(B_1)}^2
$$
{\sl hold.}

\vspace{2mm}

{\sl Proof of Claim 2:} First we notice that from (\ref{EstOna00}) and  (\ref{InteBlind})
\begin{equation}\label{CSLessThanEps}
(1-\sigma(\epsilon)-C\epsilon)\epsilon\le \|\Reg(\nabla (C_0(x)+S_0(x))\|_{L^2(B_1)}\le (1+\sigma(\epsilon)+C\epsilon)\epsilon.
\end{equation}
In particular, this implies that if $\epsilon$ is small enough then
\begin{equation}\label{CSLessThanEps2}
\|\Reg(\nabla (C_0(x)+S_0(x))\|_{L^2(B_1)}\approx \epsilon.
\end{equation}

By (\ref{CSLessThanEps}) it is enough to prove that for some $K$
\begin{equation}\label{ItFora0}
\left|\frac{a_0(0)}{\|\Reg(\nabla (S_0+C_0))\|_{L^2(B_1)}}\right|\le 
K\left\|\Reg(\nabla (S_0+C_0)\right\|_{L^2(B_1)}.
\end{equation}

The idea is to show that if (\ref{ItFora0}) is violated then 
\begin{equation}\label{AnotherThingToShow}
 \left|\frac{a_0(1)}{\|\Reg(\nabla (S_1+C_1))\|_{L^2(B_1)}}\right|\ge 
 \frac{K\|\Reg(\nabla (S_1+C_1))\|_{L^2(B_1)}}{s^{3/2}},
\end{equation}
by iterating this procedure we would, after finitely many iterations, get a contradiction to 
Lemma \ref{LemVarInOrth}.

We may estimate, using the relation $u_1(x)=\frac{u_0(sx)}{\sqrt{s}}$ as in Claim 1,
\begin{equation}\label{EstOfa10Below}
 |a_0(1)|\ge |a_0(0)|
 -\left\|\Reg\left( \nabla\frac{\hat{R}_0(sx)}{\sqrt{s}}\right)\right\|_{L^2(B_1(0))} \ge
 \end{equation}
 $$
 \ge  |a_0(0)|- 
 C_Rs\left\|\Reg(\nabla (S_0+C_0+a_0(0)\mathfrak{z})\right\|_{L^2(B_1)}^2 \ge
$$
$$
 \frac{1}{2}|a_0(0)|- C\left\|\Reg(\nabla (S_0+C_0)\right\|_{L^2(B_1)}^2\ge 
 \frac{1}{2}|a_0(0)|-C\epsilon^2,
$$
where we used (\ref{CSLessThanEps}) (and (\ref{CSLessThanEps2})) in the last inequality.

From Claim 1 we know that 
\begin{equation}\label{SmokeAndCoffee}
\|\Reg(\nabla (S_1(r,\phi)+C_1(r,\phi)))\|_{L^2(B_1(0))}\le Cs\|\Reg(\nabla(S_0+C_0))\|_{L^2(B_1)}+C\epsilon^2.
\end{equation}

Using the estimates (\ref{EstOfa10Below}) and (\ref{SmokeAndCoffee}) we can conclude that 
$$
\left|\frac{a_0(1)}{\|\Reg(\nabla (S_1(r,\phi)+C_1(r,\phi)))\|_{L^2(B_1)}}\right|\ge 
$$
\begin{equation}\label{EstOfA10}
\ge \left| \frac{\frac{1}{2}|a_0(0)|- 
 C\left\|\Reg(\nabla (S_0+C_0+a_0(0)\mathfrak{z})\right\|_{L^2(B_1)}^2 }{Cs\epsilon -C\epsilon^2}\right|\ge
\end{equation}
$$
\ge \left| \frac{\frac{1}{2}|a_0(0)|- C\epsilon^2}{Cs\epsilon -C\epsilon^2}\right|.
$$
By possibly decreasing $s$ and then $\epsilon$ it follows from (\ref{EstOfA10})
that if $|a_0(0)|\ge K \epsilon^2$ then
\begin{equation}\label{AlmostThere}
\left|\frac{a_0(1)}{\|\Reg(\nabla (S_1(r,\phi)+C_1(r,\phi)))\|_{L^2(B_1)}}\right|\ge 
\frac{K \epsilon}{s^{3/4}}.
\end{equation}
Also if $s$ and $\epsilon$ are small enough then, by Claim 1 and (\ref{CSLessThanEps2})
\begin{equation}\label{Rumpenstilskin}
\frac{\left\|\Reg(\nabla(C_1(x)+S_1(x)))\right\|_{L^2(B_1)}}{s^{3/4}}\le \epsilon.
\end{equation}
From (\ref{AlmostThere}) and (\ref{Rumpenstilskin}) we can conclude that 
\begin{equation}\label{SetUpForIterations}
\left|\frac{a_0(1)}{\|\Reg(\nabla (S_1+C_1))\|_{L^2(B_1)}}\right|\ge 
\frac{K \left\|\Reg(\nabla (C_1+S_1))\right\|_{L^2(B_1)}}{s^{3/2}}.
\end{equation}

The inequality (\ref{SetUpForIterations}) may be iterated for $j=1,2,...$ which leads to 
$$
\left|\frac{a_0(j)}{\|\Reg(\nabla (S_j+C_j))\|_{L^2(B_1)}}\right|\ge 
\frac{K \left\|\Reg(\nabla(C_j+S_j))\right\|_{L^2(B_1)}}{s^{3j/2}}\to \infty
$$
which would contradict Proposition \ref{NoLOTinLim}.

The estimate (\ref{BetterRestEst}) follows from the definition that 
$R_0=\hat{R}_0+a_0(0)\mathfrak{z}$, the estimate (\ref{InteBlind}) and the 
estimate on $a_0(0)$.

To that similar estimates hold for every $j$ follows from the scaling invariance
of the problem which implies that the above argument can be repeated for any $u_j$.

%
%

\vspace{2mm}

Next we need to show that the rotation $|\phi_{j+1}-\phi_j|$ only depends on $C_j$ and $R_j$.

\vspace{1mm}

{\bf Claim 3:} {\sl The following estimate holds}
\begin{equation}\label{RotPhiEst}
|\phi_{j+1}-\phi_j|\le  C\left( \|\Reg(\nabla C_j(r,\phi))\|_{L^2(B_1)}+\|\Reg(\nabla R_j(r,\phi))\|_{L^2(B_1)}\right).
\end{equation}

\vspace{1mm}

{\sl Proof of Claim 3:} By scaling invariance there is no loss of generality to prove this only for $j=0$, that is to
show (\ref{RotPhiEst}) for $\phi_0=0$ and $\phi_1$.

By trigonometric identities it follows that, where $\lambda(j)$ is used to denote the 
constant in $\Pi(u_j,1)=\lambda(j) r^{1/2}\sin((\phi+\phi_j)/2)$,
\begin{equation}\label{SomeTrigIdForPi}
\Pi(u_1,1)-\Pi(u,1)=\lambda(1)r^{1/2}\sin((\phi+\phi_1)/2)-\lambda(0)r^{1/2}\sin(\phi/2)=
\end{equation}
$$
=(\lambda(1)\cos(\phi_1/2)-\lambda(0))r^{1/2}\sin(\phi/2)+\lambda(1)r^{1/2}\sin(\phi_1/2)\cos(\phi/2).
$$
Also by the linearity of the projection
\begin{equation}\label{LinOfProj}
\Pi(u_1,1)=\Pi(u+(u_1-u),1)=\Pi(u,1)+\Pi(u_1-u,1).
\end{equation}

Comparing (\ref{SomeTrigIdForPi}) and (\ref{LinOfProj}) we see that the angle $\phi_1$ is determined by the projection of $u_1-u$
onto expressions of the form $a\cos(\phi/2)+b\sin(\phi/2)$. But using the explicit expression of $u$ and $u_1$ (equation (\ref{ExpressionForuj})) we see that 
\begin{equation}\label{SrthSeries}
u_1-u=\sum_{k=2}^\infty a_k(1)(s^{\alpha_k-1/2}-1)r^{\alpha_k}\cos(\alpha_k \phi)+
\end{equation}
$$
+\sum_{k=2}^\infty b_k(1)(s^{k-1}-1)r^{k-1/2}\sin((k-1/2)\phi)+\frac{R_1(s x)}{\sqrt{s}}-R_1(x).
$$
Since the sine series in (\ref{SrthSeries}) is orthogonal to every expression $a\cos(\phi/2)+b\sin(\phi/2)$ it clearly follows that 
the projection of $u_1-u$ onto expressions of the form $a\cos(\phi/2)+b\sin(\phi/2)$ only depend on the the cosine series and rest terms in 
(\ref{SrthSeries}). Since the projection decreases the norm the claim follows.

\vspace{2mm}

Next we state a simple technical claim. 

\vspace{1mm}

{\bf Claim 4:} {\sl If we denote }
$$
\xi_j=\|\Reg(\nabla S_j)\|_{L^{2}(B_1)},
$$
$$
\eta_j = \|\Reg(\nabla C_j)\|_{L^{2}(B_1)}
$$
{\sl and }
$$
\zeta_j =\|\Reg(\nabla R_j)\|_{L^{2}(B_1)}
$$
{\sl then (by Claim 2) }
\begin{equation}\label{IneqTau}
\zeta_j \le C(\xi_j^2+\eta_j^2)
\end{equation}
{\sl and (by Claim 1)}
\begin{equation}\label{IneqEps}
\xi_{j+1}\le s\xi_j +C\zeta_j
\end{equation}
{\sl and for any $\alpha <\alpha_2-1/2$ provided that $s$ is small enough 
(by an argument directly analogous to Claim 1)}
\begin{equation}\label{IneqDel}
\eta_{j+1}\le s^{\alpha}\eta_j +C\zeta_j.
\end{equation}
{\sl In particular, given an $\alpha <\alpha_2-1/2$, then provided that $s$ and $\epsilon$ are small enough, }
\begin{equation}\label{SimpleOm1}
\eta_{j}\le s^{j\alpha}\eta_0\le s^{j\alpha}\epsilon.
\end{equation}
{\sl and for any $\beta<2$}
\begin{equation}\label{SimpleOm2}
\|\Reg(\nabla R_j)\|_{L^2(B_1)}\le Cs^{j\beta}\epsilon^2.
\end{equation}

\vspace{2mm}

The proofs of (\ref{SimpleOm1}) and (\ref{SimpleOm2}) follows from simple iterations of the 
inequalities (\ref{IneqTau}), (\ref{IneqEps}) and (\ref{IneqDel}) and therefore omitted.
The proof of (\ref{IneqDel}) is exactly the same as the proof of Claim 1 and also omitted.

\vspace{2mm}

The proof of the Theorem follows in the same way as the proof of Theorem \ref{C1alpha}. It follows in particular, from Claim 3 and Claim 4, that 
\begin{equation}\label{BetterImprovement}
 \|\Reg(\nabla (\Pi(u_j,1)-\Pi(u_{j+1},1)))\|_{L^2(B_1)}\le C(\eta_j+\zeta_j)\le
\end{equation}
$$
\le  C(s^{j\alpha}\epsilon+s^{j\beta}\epsilon^2)\le Cs^{j\alpha}\epsilon.
$$
But the estimate (\ref{BetterImprovement}) is the same estimate as (\ref{StupidTriangle}) with $\alpha<\alpha_2-1/2$
in place of $\alpha<1$. The rest of the proof of Theorem \ref{C2alpha} follows exactly as the proof of Theorem \ref{C1alpha}
using (\ref{BetterImprovement}) in place of (\ref{StupidTriangle}). \qed


\appendix

\section{Analysis of the Linearized System.}\label{Sec:AnalysisLin}

In order to derive any information from Proposition \ref{Linearization} we need to understand the linear system that $(v^0,f_0)$
solves. The aim of this appendix is to prove the following proposition and a simple corollary stated at the 
end of the appendix (Corollary \ref{flatimpForLin}).

\begin{prop}\label{AnalysisOfLinear}
 For each $g\in H^{1/2}(\partial B_1(0)\setminus \{(-1,0)\})$ and $t\in \R$ there exists a weak solution
 $(v(x_1,x_2), f(x_1))\in W^{1,2}(B_1(0)\setminus \{(x_1,0);\; x_1\le 0\})\times W^{1,2}((-1,0))$
 to the following boundary value problem
 \begin{equation}\label{theLinearSyst}
 \begin{array}{ll}
  \Delta v=0 & \textrm{ in } B_1\setminus\{ x_1<0, x_2=0\} \\
    -\frac{\partial v(x_1,0^+)}{\partial x_2}=\frac{1}{2}\sqrt{\frac{2}{\pi}}\frac{\partial}{\partial x_1}\left(\frac{1}{\sqrt{-x_1}}f(x_1)\right)
    & \textrm{ for }x_1<0 \\
    -\frac{\partial v(x_1,0^-)}{\partial x_2}=-\frac{1}{2}\sqrt{\frac{2}{\pi}}\frac{\partial}{\partial x_1}\left(\frac{1}{\sqrt{-x_1}}f(x_1)\right)
    & \textrm{ for }x_1<0 \\
    \frac{\partial^2 f(x_1)}{\partial x_1^2}=
    \sqrt{\frac{2}{\pi}\frac{1}{r}}\left(\frac{\partial v(x_1,0^+)}{\partial x_1}+\frac{\partial v(x_1,0^-)}{\partial x_1}\right)
     & \textrm{ for }x_1<0 \\
     v=g & \textrm{ on } \partial B_1(0)\setminus \{(-1,0)\} \textrm{ and} \\
     f(-1)=t. &
 \end{array}
 \end{equation}

 Furthermore, if $(v,f)$ satisfies the following estimates, for any constants $C_v$ and $C_f$ and $\kappa<1/2$,
 \begin{equation}\label{snurrfotolin}
 \frac{1}{r}\int_{B_r(0)\setminus \{(x_1,0);\; x_1\le 0\}}|\nabla v|^2\le C_v |r|^{-2\kappa},
\end{equation}
\begin{equation}\label{fpointwiselin}
 |f(x_1)|\le C_f |x_1|^{1-\kappa}  \quad \textrm{ for } x_1<0
 \end{equation}
 and
 \begin{equation}\label{fprimeestlin}
 |f'(x_1)|\le C_f|x_1|^{-\kappa}  \quad \textrm{ for } x_1<0
 \end{equation}
 then the solution is unique and we may express $v$, in polar coordinates,
 \begin{equation}\label{SumOfHomos}
 v(r,\phi)=a+a_0 \mathfrak{z}(r,\phi)+\sum_{k=1}^\infty a_k r^{\alpha_k}\cos(\alpha_k \phi)+ \sum_{k=1}^\infty b_k r^{k-1/2}\sin((k-1/2)\phi)
 \end{equation}
 and
 \begin{equation}\label{SumofHetros}
 f(x_1)=a_0 \mathfrak{h}(|x_1|)+\sum_{k=1}^\infty a_k 2\sqrt{\frac{\pi}{2}}\sin(\alpha_k \pi)|x_1|^{\alpha_k+\frac{1}{2}}
 \end{equation}
 for some constants $a_k$ and $b_k$. Here $\alpha_1=1/2$ and $\alpha_k$, for $k\ge 2$, are the positive solutions, $\alpha_k\in (k-1,k)$, to the following equation
 \begin{equation}\label{EqForAlpha}
 \tan\left( \alpha \pi\right)=-\frac{2}{\pi}\frac{\alpha}{\alpha^2-\frac{1}{4}},
 \end{equation}
 the function $\mathfrak{z}$ is defined by 
 \begin{equation}\label{DefOfW}
 \mathfrak{z}(r,\phi)=r^{1/2}\phi\sin(\phi/2)-r^{1/2} \ln(r)\cos(\phi/2)
 \end{equation}
and $\mathfrak{h}$ is defined by
 \begin{equation}\label{DefOfh}
 \mathfrak{h}(r)=-\sqrt{2\pi}r\ln(r).
 \end{equation}
\end{prop}

{\bf Remark:} {\sl The values of $\kappa$, $C_v$ and $C_f$ are not important in this Proposition -- as long as the solution $v$ does not blow up too fast
at the origin, and the derivatives of $f$ do not blow up too fast, the solution is unique. This is analogous 
with the fact that uniqueness for, say, the Dirichlet problem in a smooth domain is only assured if we restrict the underlying space 
to, say, $W^{1,2}$.}

The series expansions (\ref{SumOfHomos}) and (\ref{SumofHetros}) immediately implies the following estimates.

\begin{cor}\label{Sigma}
 If $(v,f)$ is a solution as in Proposition \ref{AnalysisOfLinear} that satisfies (\ref{snurrfotolin}), (\ref{fpointwiselin}) and (\ref{fprimeestlin})
 then $(v,f)$ satisfies the following estimates
  \begin{equation}\label{snurrfotolin2}
 \frac{1}{r}\int_{B_r(0)\setminus \{(x_1,0);\; x_1\le 0\}}|\nabla v|^2\le C_v \ln(1/|r|)^{2},
\end{equation}
\begin{equation}\label{fpointwiselin2}
 |f(x_1)|\le C_f \ln(1/|x_1|)  \quad \textrm{ for } x_1<0
 \end{equation}
 and
 \begin{equation}\label{fprimeestlin2}
 |f'(x_1)|\le C_f\ln(1/|x_1|)  \quad \textrm{ for } x_1<0.
 \end{equation}
\end{cor}

Again the proof is quite long and we will therefore split it into several Lemmata. Our goal is to show that the
limit $v^0$ in Proposition \ref{Linearization} can be expressed as a series of homogeneous functions and 
$\mathfrak{z}$ as in (\ref{SumOfHomos}).
To that end we begin to derive an expression of all homogeneous solutions to (\ref{theLinearSyst}) in the next sub-section. In sub-section
\ref{SectSpan} we will show that these homogeneous functions span $L^2(\partial B_1(0)\setminus (-1,0))$ and thus $H^{1/2}(\partial B_1(0)\setminus (-1,0))$.
It follows that, for any boundary data $g$, we can find a solution $(u,f)$ such that $u=g$ on $\partial B_1(0)$. However, we also need to
specify the boundary data of $f(x_1)$ at $x_1=-1$. We show that that is possible in sub-section \ref{SectExt}. This shows that
we may find a solution for each boundary data $g$ on $\partial B_1(0)$ and $t=f(-1)\in \R$. We also need to show the uniqueness of these solutions
in order to conclude that the particular solution $(u^0,f_0)$, we get from the linearization (in Proposition \ref{Linearization}), has the desired form.
We show uniqueness in sub-section \ref{SEctUnique}. The proof of proposition \ref{AnalysisOfLinear} is then a simple consequence of the sub-sections
\ref{SectOfHomos}-\ref{SectExt}.

\subsection{Homogeneous solutions to (\ref{theLinearSyst}).}\label{SectOfHomos}

In this sub-section we classify the homogeneous solutions to the system (\ref{theLinearSyst}).

\begin{lem}
 Let $(v,f)$ be a homogeneous solution to (\ref{theLinearSyst}) then
$$
v= a r^{\alpha_k} \cos(\alpha_k \phi) \quad\textrm{ and } f=2a\sqrt{\frac{\pi}{2}}\sin(\alpha_k \pi)r^{\alpha_k+\frac{1}{2}}
$$
where $\alpha_k>1$ is a solution to
$$
\tan\left( \alpha \pi\right)=-\frac{2}{\pi}\frac{\alpha}{\alpha^2-\frac{1}{4}}
$$
or
$$
v(r,\phi)= a r^{1/2}\cos(\phi/2) \textrm{ and } f(x_1)=2\sqrt{\frac{\pi}{2}}a x_1
$$
or
$$
v=br^{k-1/2}\sin((k-1/2)\phi) \textrm{ and }f=0
$$
where $k\in \mathbb{N}$.
\end{lem}
\textsl{Proof:} We aim to derive expressions for the homogeneous solutions to the system (\ref{theLinearSyst}).
Let us therefore assume that $v(r,\phi)=r^\alpha \Phi(\phi)$ and $f(r)=c r^\beta$ for some constants $\alpha,\beta$ and $c$.
Since $v$ is harmonic it follows that $\Phi(\phi)=a\cos(\alpha \phi)+b\sin(\alpha \phi)$ for some constants $a$ and $b$.

If we first consider the part of the homogeneous solution that is even in $x_2$. We see that the equations involving $f$ in (\ref{theLinearSyst})
reduces to
\begin{equation}\label{mixingterms}
\begin{array}{ll}
    a\alpha r^{\alpha-1}\sin(\alpha\pi)=c\frac{1}{2}\sqrt{\frac{2}{\pi}}\left( \beta-\frac{1}{2}\right)r^{\beta-\frac{3}{2}}
    & \textrm{ for }x_1<0 \\
    -2a\sqrt{\frac{2}{\pi}}\alpha r^{\alpha-\frac{3}{2}}\cos(\alpha \pi)=c\beta\left( \beta-1\right)r^{\beta-2}
    & \textrm{ for }x_1<0.
 \end{array}
\end{equation}
Since the exponents in $r$ in the first equation in (\ref{mixingterms}) must agree we see that $\beta=\alpha+\frac{1}{2}$. A simple calculation
shows that (\ref{mixingterms}) is only solvable if 
\begin{equation}
 \alpha=\alpha_1=\frac{1}{2}
\end{equation}
or if $\alpha$ is a solution to
\begin{equation}\label{alphaCond}
\tan\left( \alpha \pi\right)=-\frac{2}{\pi}\frac{\alpha}{\alpha^2-\frac{1}{4}},
\end{equation}
we will denote the solutions to (\ref{alphaCond}) $\alpha_2,\alpha_3,...$

\vspace{2mm}

\includegraphics[width=10cm,height=6.5cm]{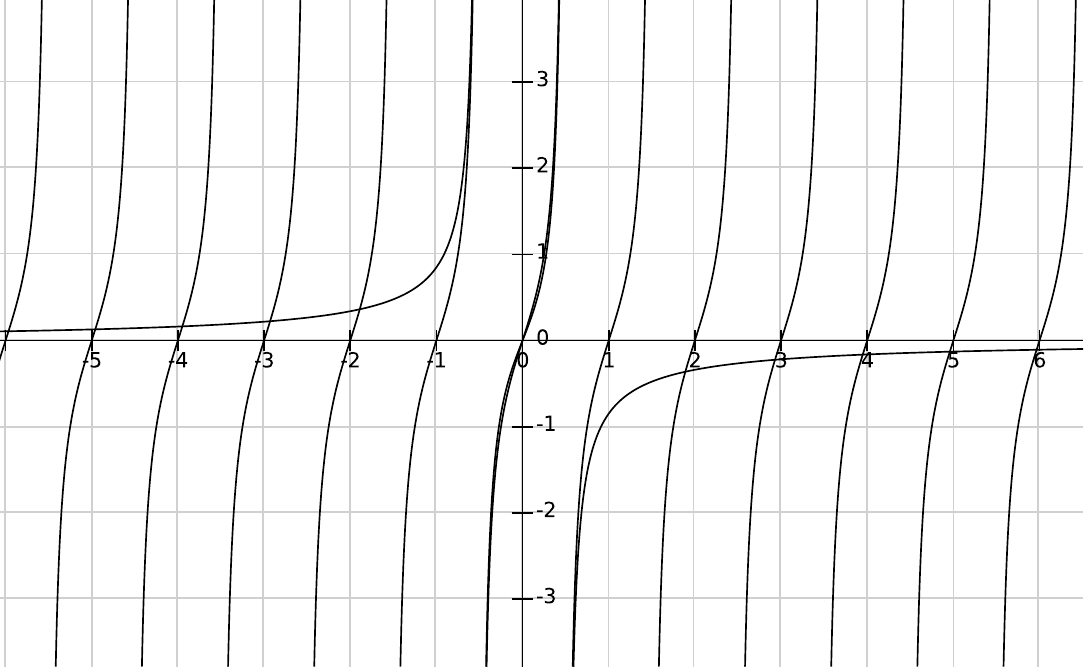}

{\bf Figure 3:} {\sl Graph that shows the values of $\alpha_k$ for $k=1,2,3,4, 5$ and $6$.}

\vspace{2mm}

Given this we can conclude that if $v$ and $f$ are homogeneous solutions, $v$ is even, then up to a multiplicative constant
$$
(v,f)=\left( r^{\alpha_k} \cos(\alpha_k \phi), 2\sqrt{\frac{\pi}{2}}\sin(\alpha_k \pi)r^{\alpha_k+\frac{1}{2}} \right)
$$
where $\alpha_k$ is a solution to (\ref{alphaCond}) or if $\alpha=1/2$ in which case 
$$
v(r,\phi)= a r^{1/2}\cos(\phi/2) \textrm{ and } f(x_1)=2\sqrt{\frac{\pi}{2}}a x_1.
$$

\vspace{2mm}

The argument is similar if $v$ is odd in $x_2$. In this case the equations involving $f$ in (\ref{theLinearSyst})
with $(v,f)=(br^\alpha \sin(\alpha \phi), cr^\beta )$ reduces to
\begin{equation}\label{mixingtermsOdd}
\begin{array}{ll}
    -b\alpha r^{\alpha-1}\cos(\alpha\pi)=c\frac{1}{2}\sqrt{\frac{2}{\pi}}\left( \beta-\frac{1}{2}\right)r^{\beta-\frac{3}{2}}
    & \textrm{ for }x_1<0 \\
    0=c\beta\left( \beta-1\right)r^{\beta-2}
    & \textrm{ for }x_1<0.
 \end{array}
\end{equation}
This implies that $c=0$ and thus that either $a=0$ or $\cos(\alpha \pi)=0$. We may conclude that $f=0$ and that $\alpha=k-\frac{1}{2}$, for $k\in \mathbb{N}$,
in case $v$ is odd in $x_2$.\qed

We end this subsection by remarking that $(\mathfrak{z},\mathfrak{h})$ is a solution to (\ref{theLinearSyst}), the proof is a simple calculation.

\begin{lem}
 The pair $(\mathfrak{z},\mathfrak{h})$, defined in (\ref{DefOfW}) and (\ref{DefOfh}) in Proposition \ref{AnalysisOfLinear}, is a solution to (\ref{theLinearSyst}).
\end{lem}


\subsection{The homogeneous solutions span $L^2(\partial B_1(0)\setminus \{(-1,0)\})$.}\label{SectSpan}

In this subsection we show that the set
$$
\left\{ 1, \cos(\alpha_1 \phi), \cos(\alpha_2\phi), \cos(\alpha_2\phi),...\right\}
$$
spans $L^2_{\textbf{even}}(\partial B_1)$, that is all the even functions in $L^2((-\pi,\pi))$, where $\alpha_k>0$ are the solutions to 
(\ref{EqForAlpha}) for $k\ge 2$ and $\alpha_1=1/2$.
Since $\sin((k-1/2)\phi)$ spans all the odd functions in $L^2((-\pi,\pi))$ it follows that for any $g\in H^{1/2}(\partial B_1\setminus (-1,0))$
we can find a pair of solutions $(v,f)$ to (\ref{theLinearSyst}) where $v=g$ on $\partial B_1$ and is on the form
\begin{equation}\label{vreducedform}
 v(r,\phi)=a+\sum_{k=1}^\infty a_k r^{\alpha_k}\cos(\alpha_k \phi)+ \sum_{k=1}^\infty b_k r^{k-1/2}\sin((k-1/2)\phi)
\end{equation}
Notice that the expression of $v$ in (\ref{vreducedform}) does not contain the term $\mathfrak{z}(r,\phi)$ that appears in
the form of $v$ in Proposition \ref{AnalysisOfLinear}. Later we will see that this missing term allows us to specify the values of $f(-1)$.

To show that $1$ and $\{\cos(\alpha_k \phi);\; k=1,2,...\}$ spans the set of even functions on the sphere we will show that
this set has the same span as $\{\cos(j\phi);\; j=0,1,2,...\}$.

Consider the linear map $A$ defined, for some $t_k\in \R$ to be specified shortly, by
\begin{equation*}
A\left(\cos(k\phi)\right)= t_k\cos(\alpha_k \phi).
\end{equation*}
Since 
$\left\{\cos(\alpha_1 \phi), \cos(\alpha_2\phi), \cos(\alpha_2\phi),...\right\}$ spans $L^2_{\textbf{even}}((-\pi,\pi))$,
$A$ defines a linear map 
$$
A:L^2_{\textbf{even}}((-\pi,\pi))\mapsto L^2_{\textbf{even}}((-\pi,\pi))
$$
where $L^2_{\textbf{even}}((-\pi,\pi))$ consists of the even functions of $L^2((-\pi,\pi))$, and
\begin{equation}\label{MapL2toL2}
A\left( \sum_{k=0}^\infty a_k\cos(k\phi))\right)= \sum_{k=0}^\infty a_k t_k\cos(\alpha_k \phi)).
\end{equation}

We claim that we can choose the numbers $t_k$ such that 
$$
\|A-I\|_{L^2_{\textbf{even}}((-\pi,\pi))\mapsto L^2_{\textbf{even}}((-\pi,\pi))}<1,
$$ 
where $\|\cdot \|_{L^2\mapsto L^2}$ denotes the operator norm.
From this it clearly follows that
$A$ is invertible and that $\left\{1, \cos(\alpha_1 \phi), \cos(\alpha_2\phi), \cos(\alpha_2\phi),...\right\}$ spans $L^2_{\textbf{even}}((-\pi,\pi))$.

Notice that, for any even function $u=\sum_{k=0}^\infty a_k\frac{\cos(k\phi))}{\sqrt{\pi}}$ the following estimate holds
$$
\|(A-I)u\|_{L^2((-\pi,\pi))}=\left\|\sum_{k=1}^\infty a_k \frac{ t_k\cos(\alpha_k \phi)-\cos( k\phi)}{\sqrt{\pi}}\right\|_{L^2}\le
$$
$$
\le \sum_{k=1}^\infty |a_k| \left\|\frac{ t_k \cos(\alpha_k \phi)-\cos( k\phi)}{\sqrt{\pi}}\right\|_{L^2}\le
$$
$$
\le \left(\sum_{k=1}^\infty |a_k|^2 \right)^{1/2}\left( \sum_{k=1}^\infty \left\|\frac{ t_k \cos(\alpha_k \phi)-\cos( k\phi)}{\sqrt{\pi}}\right\|_{L^2}^2\right)^{1/2}=
$$
$$
=\|u\|_{L^2((-\pi,\pi))}\left( \sum_{k=1}^\infty \left\|\frac{ t_k\cos(\alpha_k \phi)-\cos( k\phi)}{\sqrt{\pi}}\right\|_{L^2}^2\right)^{1/2},
$$
where we used the triangle inequality and H\"older's inequality in as well as Paresval's equality. It is therefore enough to show that
\begin{equation}\label{OperatorNormEst}
\sum_{k=1}^\infty \left\|\frac{ t_k\cos(\alpha_k \phi)-\cos( k\phi)}{\sqrt{\pi}}\right\|_{L^2}^2<1
\end{equation}
for some choice of numbers $t_k$.

It is easy to see that $\left\| t_k\cos(\alpha_k \phi)-\cos( k\phi) \right\|_{L^2}^2$ is minimized for the 
value 
\begin{equation}\label{Exprtk}
t_k=\frac{\int_{-\pi}^\pi  \cos(\alpha_k\phi)\cos(k \phi)d\phi}{\int_{-\pi}^\pi  \cos^2(\alpha_k\phi)d\phi}=
\frac{\frac{2\alpha_k}{\alpha_k+k} \frac{\sin\pi(k-\alpha_k)}{\pi(k-\alpha_k)} }{1+ \frac{\sin(2\alpha_k \pi)}{2\alpha_k\pi}}
\end{equation}

We will need to estimate $t_k$, $k\ge 2$, from below. To that end we begin by estimating the denominator in (\ref{Exprtk}).
For that we will use the number $\gamma_k$ defined by $k-\gamma_k=\alpha_k$ and the following calculation
\begin{equation}\label{Denom}
 \frac{1}{1+\frac{\sin(2\alpha_k \pi)}{2\alpha_k \pi}}\ge 1-\frac{\sin(2\pi \alpha_k)}{2\pi \alpha_k}=
\end{equation}
$$
=1-\frac{\sin\left( 2\pi k-2\pi\gamma_k\right)}{2\pi\alpha_k}=1+\frac{\sin(2\pi \gamma_k)}{2\pi\alpha_k}\ge 
1+\frac{\gamma_k}{\alpha_k}-\frac{4\pi^2 \gamma_k^3}{2\pi \alpha_k},
$$
where we used the elementary inequality $\sin(t)\ge t-t^3/6$ for $t\ge 0$ in the last inequality.

To estimate the numerator in (\ref{Exprtk}) we use the following calculation
\begin{equation}\label{Numerator}
 \frac{2\alpha_k}{\alpha_k+k}\frac{\sin(\pi (k-\alpha_k))}{\pi(k-\alpha_k)}=\frac{2\alpha_k+\gamma_k}{\alpha_k+k}\frac{\sin(\pi \gamma_k)}{\pi\gamma_k}-
\end{equation}
$$
-\frac{\gamma_k}{\alpha_k+k}\frac{\sin(\pi \gamma_k)}{\pi \gamma_k}\ge 1-\frac{\pi^2 \gamma_k^2}{6}-\frac{\gamma_k}{\alpha_k+k},
$$
where we used the elementary inequalities $x-\frac{x^3}{6}\le \sin(x)\le 1$ in the last estimate.

Next we calculate each term in the series in (\ref{OperatorNormEst})
\begin{equation}\label{CalculusEstimate1}
\frac{1}{\pi}\int_{-\pi}^\pi\left( t_k \cos(\alpha_k \phi)-\cos( k\phi)\right)^2d\phi=
\end{equation}
$$
=\frac{1}{\pi}\int_{-\pi}^\pi  t_k^2 \cos^2(\alpha_k\phi)-2t_k\cos(\alpha_k \phi)\cos(k \phi)+\cos^2(k \phi) d\phi=
$$
\begin{equation}\label{EstEachTerm}
=\frac{1}{\pi}\int_{-\pi}^\pi \cos^2(k \phi) -t_k^2 \cos^2(\alpha_k\phi)d\phi=1-t_k^2 \frac{1}{\pi}\int_{-\pi}^\pi\cos^2(\alpha_k\phi)d\phi=
\end{equation}
$$
=1-t_k \frac{1}{\pi}\int_{-\pi}^\pi\cos(\alpha_k\phi)\cos(k\phi)d\phi
=1-t_k\frac{2\alpha_k}{\alpha_k+k} \frac{\sin\pi(k-\alpha_k)}{\pi(k-\alpha_k)}. 
$$

Using (\ref{Denom}) and (\ref{Numerator}) to estimate the right side in (\ref{EstEachTerm}) we get 
$$
1-t_k\frac{2\alpha_k}{\alpha_k+k} \frac{\sin\pi(k-\alpha_k)}{\pi(k-\alpha_k)}\le 1-\left( 1-\frac{\pi^2 \gamma_k^2}{6}-\frac{\gamma_k}{\alpha_k+k}\right)^2
\left(1+\frac{\gamma_k}{\alpha_k}-\frac{4\pi^2 \gamma_k^3}{2\pi \alpha_k}\right)\le
$$
\begin{equation}\label{NeedNastyEst}
\le \frac{\pi^2 \gamma_k^2}{3}+\frac{2\gamma_k}{2k-\gamma_k}-\frac{\gamma_k}{\alpha_k}+\frac{4\pi^2\gamma_k^3}{6\alpha_k}+
\left( \frac{\pi^2\gamma_k^2}{3}+\frac{2 \gamma_k}{2k-\gamma_k}\right)\frac{\gamma_k}{\alpha_k}.
\end{equation}

We need to estimate $\gamma_k=k-\alpha_k$ in order to continue. Since $\alpha_k$ is the solution of (\ref{EqForAlpha}) that satisfies $\alpha_k\in (k-1,k]$
it follows that
\begin{equation}\label{EgForGammak}
\tan(\pi\gamma_k)+\frac{2}{\pi}\frac{k-\gamma_k}{(k-\gamma_k)^2-\frac{1}{4}}=0.
\end{equation}
We aim to show that $\gamma_k$ satisfies $0<\gamma_k<\frac{1}{4k}$ for $k\ge 2$. 

Since $-\pi \gamma_k\le \tan(\pi\alpha_k)$ when $\alpha_k=k-\gamma_k$ and $0\le \gamma_k < 1$ it follows that we may estimate $\gamma_k$
from above by the solution $\gamma$ to 
$$
0=-\pi \gamma+\frac{2}{\pi}\frac{2-\gamma}{(k-\gamma)^2-1/4}=
$$
$$
=\frac{-\pi^2 \gamma^3+2\pi^2 k \gamma^2+(-\pi^2k^2-2+\pi^2/4)\gamma+2k}{\pi(k-\gamma)^2-1/4}=\frac{q(\gamma)}{\pi((k-\gamma)^2-1/4)},
$$
where $q(\gamma)$ is defined by the last expression. To show that $0<\gamma_k<1/(4k)$ it is enough, by the mean value property, to show
that $q(0)>0$, which is trivial, and that $q(1/4k)<0$ which follows from
$$
k^3q(1/4k)=\left( 2-\frac{\pi^2}{4}\right)k^4+\left(\frac{\pi^2}{2}+\frac{\pi^2}{4}-2 \right)\frac{k^2}{4}-\frac{\pi^3}{64}\le
$$
$$
\le \left( 8-\pi^2+\frac{3\pi^2}{16}\right)k^2-\frac{\pi^3}{64}<0,
$$
where we used that $k\ge 2$ in the first inequality. It follows that 
\begin{equation}\label{EstingGamma}
 0<\gamma_k <\frac{1}{4k}\quad\textrm{ for }\quad k\ge 2.
\end{equation}

The estimate on $\gamma_k$ allows us to estimate (\ref{NeedNastyEst}) from above by 
$$
\frac{\pi^2}{48}\frac{1}{k^2}+2\frac{1}{8k^2-1}-\frac{1}{4k^2-1}+\frac{\pi^2}{96(k^4-k^2/4)}+\left(\frac{\pi^2}{48 k^2}+\frac{2}{8k^2-1}\right)\frac{1}{4k^2-1}=
$$
\begin{equation}\label{IntroducingT}
=T_1(k)+T_2(k)+T_3(k)+T_4(k)+T_5(k).
\end{equation}

Let us recapitulate what we have done and what we are aiming to do. We are aiming to prove (\ref{OperatorNormEst}). 
In (\ref{EstEachTerm}) we rewrote each term in 
(\ref{OperatorNormEst}) and then estimated the terms, for $k\ge 2$ by (\ref{IntroducingT}), it
follows that we have to 
show that the following expression is less than one
\begin{equation}\label{NightEyes}
1-t_1\frac{2\alpha_1}{\alpha_1+1} \frac{\sin\pi(1-\alpha_1)}{\pi(1-\alpha_1)}+\sum_{k=2}^\infty \left( T_1(k)+T_2(k)+T_3(k)+T_4(k)+T_5(k)\right).
\end{equation}
Using that $\alpha_1=1/2$ we may evaluate 
$$
1-t_1\frac{2\alpha_1}{\alpha_1+1} \frac{\sin\pi(1-\alpha_1)}{\pi(1-\alpha_1)}=1-\frac{16}{9\pi^2}< 0.820.
$$
Next we evaluate the series in (\ref{NightEyes}) one at the time. First we use that $\sum_{k=1}^\infty \frac{1}{k^2}=\frac{\pi^2}{6}$ to evaluate
$$
\sum_{k=2}^\infty T_1(k)=\frac{\pi^2}{48}\sum_{k=2}^\infty \frac{1}{k^2}=\frac{\pi^2}{48}\left(\frac{\pi^2}{6}-1 \right)< 0.133.
$$
To estimate $\sum_{k=2}^\infty T_2(k)$ we use a summation formula due to Euler, $\sum_{k=1}^\infty \frac{1}{k^2-a^2}=\frac{1}{2a^2}-\frac{\pi}{2a\tan(a\pi)}$,
$$
\sum_{k=2}^\infty T_2(k)=\frac{1}{4}\sum_{k=2}^\infty\left(\frac{1}{k^2-(1/2\sqrt{2})^2} \right)=\frac{1}{4}\left( 4-\frac{\sqrt{2}\pi}{\tan(\pi/2\sqrt{2})}-\frac{8}{7}\right)< 0.164.
$$
For $\sum_{k=2}^\infty T_3(k)$ we use that $\sum_{k=1}^\infty \frac{1}{4k^2-1}=\frac{1}{2}$, which may be seen by evaluating the Fourier-series of $|\sin(x)|$ on 
$(-\pi,\pi)$ at $x=0$,
$$
\sum_{k=2}^\infty T_3(k)=-\sum_{k=2}^\infty \frac{1}{4k^2-1}=-\frac{1}{6}<-0.166.
$$
For $\sum_{k=2}^\infty T_4(k)$ we use that $\sum_{k=1}^\infty\frac{1}{k^4}=\frac{\pi^4}{90}$:
$$
\sum_{k=2}^\infty T_4(k)\le \sum_{k=2}^\infty \frac{\pi^2}{90k^4}=\frac{\pi^2}{90}\left( \frac{\pi^4}{90}-1\right)< 0.010.
$$
Finally we use $\sum_{k=1}^\infty \frac{1}{k^4}=\frac{\pi^4}{90}$ to estimate 
$$
\sum_{k=2}^\infty T_5(k)=\sum_{k=2}^\infty \left(\frac{\pi^2}{48 k^2}+\frac{2}{8k^2-1}\right)\frac{1}{4k^2-1}\le \sum_{k=2}^\infty\left(\frac{10}{48}+\frac{8}{31}\right)\frac{4}{15}\frac{1}{k^4}< 0.011.
$$

Using the estimates of the of $\sum_{k=2}^\infty T_j(k)$, $j=1,...,5$, in (\ref{NightEyes}) leads to 
$$
\sum_{k=1}^\infty \left\|\frac{ t_k\cos(\alpha_k \phi)-\cos( k\phi)}{\sqrt{\pi}}\right\|_{L^2}^2\le 0.820+0.133+0.164-0.166+0.010+0.011<1.
$$
It follows that the operator $A$ is invertable. We have therefore shown the following lemma.

%

\begin{lem}
 The set
 $$
\left\{ 1, \cos(\alpha_1 \phi), \cos(\alpha_2\phi), \cos(\alpha_2\phi),...\right\}
$$
where $\alpha_1=1/2$ and $\alpha_k$, for $k\ge 2$, are the positive solutions to
$$
\tan\left( \alpha \pi\right)=-\frac{2}{\pi}\frac{\alpha}{\alpha^2-\frac{1}{4}}
$$
forms a basis of all even $L^2(-\pi,\pi)$ functions.
\end{lem}

\subsection{Uniqueness of solutions.}\label{SEctUnique}

\begin{lem}
 The solutions to (\ref{theLinearSyst}) that satisfy (\ref{snurrfotolin}), (\ref{fpointwiselin}) and (\ref{fprimeestlin})
 are unique.
\end{lem}
\textsl{Proof:} The proof is standard. Therefore we will only provide an outline.

Let us denote by $\mathcal{L}_s(v,f)$, $s\in [0,1]$, the mapping that takes $(v,f)$ to $(a_1,a_2,a_3,g,t)$ defined by 
\begin{equation}\label{theLinearSysts}
 \begin{array}{ll}
  \Delta v=0 & \textrm{ in } B_1\setminus\{ x_1<0, x_2=0\} \\
    -\frac{\partial v(x_1,0^+)}{\partial x_2}-\frac{s}{2}\sqrt{\frac{2}{\pi}}\frac{\partial}{\partial x_1}\left(\frac{1}{\sqrt{-x_1}}f(x_1)\right)=a_1(x_1)
    & \textrm{ for }x_1<0 \\
    -\frac{\partial v(x_1,0^-)}{\partial x_2}+\frac{s}{2}\sqrt{\frac{2}{\pi}}\frac{\partial}{\partial x_1}\left(\frac{1}{\sqrt{-x_1}}f(x_1)\right)=a_2(x_1)
    & \textrm{ for }x_1<0 \\
    \frac{\partial^2 f(x_1)}{\partial x_1^2}-
    s\sqrt{\frac{2}{\pi}\frac{1}{r}}\left(\frac{\partial v(x_1,0^+)}{\partial x_1}+\frac{\partial v(x_1,0^-)}{\partial x_1}\right)=a_3(x_1)
     & \textrm{ for }x_1<0 \\
     v=g & \textrm{ on } \partial B_1(0)\setminus \{(-1,0)\} \textrm{ and} \\
     f(-1)=t & f(0)=0.
 \end{array}
 \end{equation}
We will consider $\mathcal{L}_s$ to be a mapping defined on the space $X$ where $(v,f)\in X$ if 
$v\in W^{1,2}(B_1(0)\setminus \Gamma_0)$, $\Delta v=0$ in $B_1(0)\setminus \Gamma_0$ and 
$v, |x_1|\frac{\partial v}{\partial x_2}\in C^{1/2}(\Gamma_0^\pm)$ and $f,x_1f', x_1^2 f''\in C^{1/2}(\Gamma_0)$.
The mapping $\mathcal{L}_s: X\mapsto Y$ where $Y$ is the space consisting quintuples $(a_1,a_2,a_3,g,t)$ where $t\in \R$, 
$g\in H^{1/2}(\partial B_1\setminus (-1,0))$ and $x_i a_i\in C^{1/2}(\Gamma_0)$ for $i=1,2,3$. We may equip $X$ and $Y$
with their natural norms, the sum of the norms of the spaces that defines $X$ and $Y$.

First we show that $\mathcal{L}_s$ is onto. This can be done by first noticing that linear
combinations of $(r^{\beta}\sin(\beta\phi),0)\in X$ and $(r^{\beta}\cos(\beta \phi),0)\in X$, where $\beta=k$ or $\beta=k+\frac{1}{2}$, and the functions 
$(0,r^k)\in X$ is mapped onto a linear subspace in $Y$ whose restriction to the first three components is dense in the subspace consisting of 
$(a_1,a_2,a_3)$. By linearity of $\mathcal{L}_s$ we may therefore reduce the problem 
of showing that $\mathcal{L}_s$ is onto to showing that it is onto the subspace $(0,0,0,g,t)\in Y$, this can be proved as in Appendix \ref{SectSpan}.

Clearly each mapping $\mathcal{L}_s$ is bounded. The bounded inverse Theorem (see Corollary 4.30 in \cite{Bowers}) implies that 
$\mathcal{L}_s^{-1}$ is bounded for each $s\in [0,1]$. By a routine application of the Banach fixed point theorem it follows that 
$\mathcal{L}_s^{-1}$ is uniformly bounded in a neighborhood of each $s$ and since $[0,1]$ is compact we may find a uniform bound on the 
inverse in the entire interval $[0,1]$.

Arguing as in the method of continuity (see Theorem 5.2 in \cite{GT}) we may show that if $\mathcal{L}_1$ has nontrivial kernel
then so will $\mathcal{L}_0$. The mapping $\mathcal{L}_0$ is just
\begin{equation}\label{theLinearSyst0}
 \begin{array}{ll}
  \Delta v=0 & \textrm{ in } B_1\setminus\{ x_1<0, x_2=0\} \\
    -\frac{\partial v(x_1,0^+)}{\partial x_2}=a_1 & \textrm{ for }x_1<0 \\
    -\frac{\partial v(x_1,0^-)}{\partial x_2}=a_2 & \textrm{ for }x_1<0 \\
    \frac{\partial^2 f(x_1)}{\partial x_1^2}=a_3 & \textrm{ for }x_1<0 \\
     v=g & \textrm{ on } \partial B_1(0)\setminus \{(-1,0)\} \textrm{ and} \\ 
     f(-1)=t & f(0)=0.
 \end{array}
\end{equation}
which is no more than the classical Dirichlet/Neumann problem in $v$ and a trivial ode in $f$. This implies that $\mathcal{L}_0$ is injective. We may 
conclude that $\mathcal{L}_1$ is injective. This finishes the proof. \qed

\subsection{Existence of solutions.}\label{SectExt}

\begin{lem}
 Given $g\in H^{1/2}(\partial B_1(0)\setminus \{(-1,0)\})$ and $t\in \R$ there exists a solution
 to (\ref{theLinearSyst}) satisfying (\ref{snurrfotolin}), (\ref{fpointwiselin}) and (\ref{fprimeestlin}).
\end{lem}
\textsl{Proof:}
Since
\begin{equation}\label{basisAgain}
\left\{ 1, \cos(\alpha_1 \phi), \cos(\alpha_2\phi), \cos(\alpha_2\phi),...\right\}
\end{equation}
span the even functions on $L^2(-\pi,\pi)$ and
$$
\left\{ \sin(\phi/2), \sin(3\phi/2), \sin(5\phi/2), ...\right\}
$$
span the odd functions on $L^2(-\pi,\pi)$ we may express any function \linebreak$g\in H^{1/2}(\partial B_1(0)\setminus (-1,0))$
according to
$$
g=a+\sum_{k=1}^\infty a_k \cos(\alpha_k \phi)+ \sum_{k=1}^\infty b_k \sin((k-1/2)\phi)
$$
clearly
$$
u(r,\phi)=a+\sum_{k=1}^\infty a_k r^{\alpha_k}\cos(\alpha_k \phi)+ \sum_{k=1}^\infty b_k r^{k-1/2}\sin((k-1/2)\phi)
$$
and
$$
 f(x_1)= \sum_{k=1}^\infty 2 a_k \sqrt{\frac{\pi}{2}}\sin(\alpha_k \pi)|x_1|^{\alpha_k+\frac{1}{2}}
$$
will solve (\ref{theLinearSyst}) except possible the condition that $f(-1)=t$.

However, since we may expand $\mathfrak{z}$ in the basis (\ref{basisAgain}) on the boundary $\partial B_1(0)$ we can find
a solution $(s,l)$
$$
s(r,\phi)=b_0+\sum_{k=1}^\infty b_k r^{\alpha_k}\cos(\alpha_k \phi)
$$
$$
l(x_1)=\sum_{k=1}^\infty b_k 2\sqrt{\frac{\pi}{2}}\sin(\alpha_k \pi)|x_1|^{\alpha_k+\frac{1}{2}}
$$
such that $s=\mathfrak{z}$ on $\partial B_1(0)$. By the uniqueness of solutions $l(-1)\ne h(-1)$ since otherwise $(\mathfrak{z},\mathfrak{h})=(s,l)$. It follows that 
the solution $(\mathfrak{z}-s,\mathfrak{h}-l)$ vanishes on $\partial B_1$ and $\mathfrak{h}(-1)-l(-1)\ne 0$. It follows that
$$
\left( u+\frac{t-f(-1)}{\mathfrak{h}(-1)-l(-1)}(\mathfrak{z}(x)-s(x)), f(x_1)+\frac{t-f(-1)}{\mathfrak{h}(-1)-l(-1)}(\mathfrak{h}(x_1)-l(x_1)\right)
$$
is a solution that satisfies (\ref{theLinearSyst}).\qed

\subsection{Regularity of the Solution to the Linearized problem.}

The following is a simple Corollary to Proposition \ref{AnalysisOfLinear}.

\begin{cor}\label{flatimpForLin}
 Let $(v,f)$ be a solution to (\ref{theLinearSyst}) of the following form
 $$
 v(r,\phi)=\sum_{k=2}^\infty a_k r^{\alpha_k}\cos(\alpha_k \phi)+ \sum_{k=2}^\infty b_k r^{k-1/2}\sin((k-1/2)\phi)
 $$
 and
 $$
 f(x_1)=\sum_{k=2}^\infty a_k 2\sqrt{\frac{\pi}{2}}\sin(\alpha_k \pi)|x_1|^{\alpha_k+\frac{1}{2}},
 $$
 where $w$ and $h$ are as in Proposition \ref{AnalysisOfLinear}. Then for each $\alpha<3/2$ there exists an 
 $s_\alpha>0$, depending only on $\alpha$, such that
 \begin{equation}\label{obviousImp}
 \|\nabla v\|_{L^2(B_{s_\alpha}(0)\setminus \Gamma_0)}< s_\alpha^{\alpha} \|\nabla v\|_{L^2(B_1\setminus \Gamma_0)}
 \end{equation}
 and
 \begin{equation}\label{obviousImpf}
 \|f'(x_1)\|_{L^2(-s_\alpha,0)}\le s_\alpha^{\alpha-1/2}\|f'\|_{L^2(-1,0)}.
 \end{equation}
\end{cor}
\textsl{Proof:} We will only show (\ref{obviousImp}), the proof of (\ref{obviousImpf}) is similar (and somewhat simpler).
Notice that
$$
\|\nabla v\|_{L^2(B_{s_\alpha}(0)\setminus \Gamma_0)}\le
\sum_{k=2}^\infty|a_k|\left\| \nabla r^{\alpha_k}\cos(\alpha_k \phi)\right\|_{L^2(B_{s_\alpha}(0)\setminus \Gamma_0)}+
$$
\begin{equation}\label{sineTerms}
\sum_{k=2}^\infty |b_k| \left\| \nabla r^{k-1/2}\sin((k-1/2)\phi)\right\|_{L^2(B_{s_\alpha}(0)\setminus \Gamma_0)}.
\end{equation}

A change of variables $r \mapsto s_\alpha r$ shows that
$$
\left\| \nabla r^{\alpha_k}\cos(\alpha_k \phi)\right\|_{L^2(B_{s_\alpha}(0)\setminus \Gamma_0)}=
\left( \int_{B_{s_\alpha}\setminus \Gamma_0} |\nabla r^{\alpha_k}\cos(\alpha_k \phi)|^2\right)^{1/2}=
$$
$$
=s_\alpha^{\alpha_k}\left( \int_{B_{1}\setminus \Gamma_0} |\nabla r^{\alpha_k}\cos(\alpha_k \phi)|^2\right)^{1/2},
$$
a similar calculation obviously works for the $r^{k-1/2}\sin((k-1/2)\phi)$ terms in (\ref{sineTerms}).

We may thus write (\ref{sineTerms})
$$
\|\nabla v\|_{L^2(B_{s_\alpha}(0)\setminus \Gamma_0)}\le
\sum_{k=2}s_\alpha^{\alpha_k}|a_k|\left\| \nabla r^{\alpha_k}\cos(\alpha_k \phi)\right\|_{L^2(B_{1}(0)\setminus \Gamma_0)}+
$$
\begin{equation}\label{sineTerms2}
\sum_{k=2}^\infty |b_k| s_\alpha^{k-1/2}\left\| \nabla r^{k-1/2}\sin((k-1/2)\phi)\right\|_{L^2(B_{1}(0)\setminus \Gamma_0)}\le
\end{equation}
$$
\le s_\alpha^{3/2}\bigg(
\sum_{k=1}|a_k|\left\| \nabla r^{\alpha_k}\cos(\alpha_k \phi)\right\|_{L^2(B_{1}(0)\setminus \Gamma_0)}+
$$
$$
+\sum_{k=2}^\infty |b_k| \left\| \nabla r^{k-1/2}\sin((k-1/2)\phi)\right\|_{L^2(B_{1}(0)\setminus \Gamma_0)}
\bigg).
$$
It is easy to see that we may estimate (\ref{sineTerms2}) by
$$
\le \left( C s_\alpha^{3/2-\alpha}\right)s_\alpha^\alpha\bigg(\left\| \sum_{k=2}a_k\nabla r^{\alpha_k}\cos(\alpha_k \phi)\right\|_{L^2(B_{1}(0)\setminus \Gamma_0)}+
$$
$$
+\left\| \sum_{k=2}^\infty b_k\nabla r^{k-1/2}\sin((k-1/2)\phi)\right\|_{L^2(B_{1}(0)\setminus \Gamma_0)}
\bigg),
$$
for some constant $C$ depending on the almost orthogonality of the basis $\cos(\alpha_k \phi)$. Choosing $s_\alpha$
small enough so that
$$
C s_\alpha^{3/2-\alpha}< 1
$$
finishes the proof. \qed

\bibliographystyle{plain}
\bibliography{ZeroCurv.bib}

\begin{thebibliography}{10}

\bibitem{AmbExist}
L.~Ambrosio.
\newblock Existence theory for a new class of variational problems.
\newblock {\em Arch. Rational Mech. Anal.}, 111(4):291--322, 1990.

\bibitem{AFP}
Luigi Ambrosio, Nicola Fusco, and Diego Pallara.
\newblock Partial regularity of free discontinuity sets. {II}.
\newblock {\em Ann. Scuola Norm. Sup. Pisa Cl. Sci. (4)}, 24(1):39--62, 1997.

\bibitem{AFPBook}
Luigi Ambrosio, Nicola Fusco, and Diego Pallara.
\newblock {\em Functions of bounded variation and free discontinuity problems}.
\newblock Oxford Mathematical Monographs. The Clarendon Press, Oxford
  University Press, New York, 2000.

\bibitem{AP}
Luigi Ambrosio and Diego Pallara.
\newblock Partial regularity of free discontinuity sets. {I}.
\newblock {\em Ann. Scuola Norm. Sup. Pisa Cl. Sci. (4)}, 24(1):1--38, 1997.

\bibitem{ASign}
John Andersson.
\newblock Optimal regularity for the {S}ignorini problem and its free boundary.
\newblock {\em Inventiones mathematicae}, pages 1--82, 2015.

\bibitem{ASW}
John Andersson, Henrik Shahgholian, and Georg~S. Weiss.
\newblock On the singularities of a free boundary through {F}ourier expansion.
\newblock {\em Invent. Math.}, 187(3):535--587, 2012.

\bibitem{B}
A.~Bonnet.
\newblock On the regularity of edges in image segmentation.
\newblock {\em Ann. Inst. H. Poincar\'e Anal. Non Lin\'eaire}, 13(4):485--528,
  1996.

\bibitem{BD}
Alexis Bonnet and Guy David.
\newblock Cracktip is a global {M}umford-{S}hah minimizer.
\newblock {\em Ast\'erisque}, (274):vi+259, 2001.

\bibitem{Bowers}
Adam Bowers and Nigel~J. Kalton.
\newblock {\em An introductory course in functional analysis}.
\newblock Universitext. Springer, New York, 2014.
\newblock With a foreword by Gilles Godefroy.

\bibitem{CL}
Antonin Chambolle and Antoine Lemenant.
\newblock The stress intensity factor for non-smooth fractures in antiplane
  elasticity.
\newblock {\em Calc. Var. Partial Differential Equations}, 47(3-4):589--610,
  2013.

\bibitem{DMSActa}
G.~Dal~Maso, J.-M. Morel, and S.~Solimini.
\newblock A variational method in image segmentation: existence and
  approximation results.
\newblock {\em Acta Math.}, 168(1-2):89--151, 1992.

\bibitem{DT}
Gianni Dal~Maso and Rodica Toader.
\newblock A model for the quasi-static growth of brittle fractures: existence
  and approximation results.
\newblock {\em Arch. Ration. Mech. Anal.}, 162(2):101--135, 2002.

\bibitem{DCone}
Guy David.
\newblock {$C^1$}-arcs for minimizers of the {M}umford-{S}hah functional.
\newblock {\em SIAM J. Appl. Math.}, 56(3):783--888, 1996.

\bibitem{DBook}
Guy David.
\newblock {\em Singular sets of minimizers for the {M}umford-{S}hah
  functional}, volume 233 of {\em Progress in Mathematics}.
\newblock Birkh\"auser Verlag, Basel, 2005.

\bibitem{DCL}
E.~De~Giorgi, M.~Carriero, and A.~Leaci.
\newblock Existence theorem for a minimum problem with free discontinuity set.
\newblock {\em Arch. Rational Mech. Anal.}, 108(3):195--218, 1989.

\bibitem{DFo}
C.~De~Lellis and M.~Focardi.
\newblock Higher integrability of the gradient for minimizers of the {$2d$}
  {M}umford-{S}hah energy.
\newblock {\em J. Math. Pures Appl. (9)}, 100(3):391--409, 2013.

\bibitem{DFi}
Guido De~Philippis and Alessio Figalli.
\newblock Higher integrability for minimizers of the {M}umford-{S}hah
  functional.
\newblock {\em Arch. Ration. Mech. Anal.}, 213(2):491--502, 2014.

\bibitem{FM}
G.~A. Francfort and J.-J. Marigo.
\newblock Revisiting brittle fracture as an energy minimization problem.
\newblock {\em J. Mech. Phys. Solids}, 46(8):1319--1342, 1998.

\bibitem{GT}
David Gilbarg and Neil~S. Trudinger.
\newblock {\em Elliptic partial differential equations of second order}.
\newblock Classics in Mathematics. Springer-Verlag, Berlin, 2001.
\newblock Reprint of the 1998 edition.

\bibitem{G}
Alan~Arnold Griffith.
\newblock The phenomena of rupture and flow in solids.
\newblock {\em Philos. Trans. R. Soc. Ser. A Math. Phys. Eng. Sci.},
  221(582-593):163--198, 1921.

\bibitem{KLM}
Herbert Koch, Giovanni Leoni, and Massimiliano Morini.
\newblock On optimal regularity of free boundary problems and a conjecture of
  {D}e {G}iorgi.
\newblock {\em Comm. Pure Appl. Math.}, 58(8):1051--1076, 2005.

\bibitem{LemMik}
Antoine Lemenant and Hayk Mikayelyan.
\newblock Stationarity of the crack-front for the {M}umford-{S}hah problem in
  3{D}.
\newblock {\em J. Math. Anal. Appl.}, 462(2):1555--1569, 2018.

\bibitem{MS}
David Mumford and Jayant Shah.
\newblock Optimal approximations by piecewise smooth functions and associated
  variational problems.
\newblock {\em Comm. Pure Appl. Math.}, 42(5):577--685, 1989.

\end{thebibliography}

\end{document}